\documentclass[a4paper,11pt]{article}

\usepackage[colorlinks=false,pdfborderstyle={/W 1}]{hyperref}

\def\arXiv#1#2{\href{http://arxiv.org/abs/#1}{{\tt arXiv:#1 [#2]}}}

\usepackage{import}
\usepackage{amsmath}
\usepackage{amssymb}
\usepackage{amsthm}
\usepackage{amsfonts}
\usepackage{graphicx}
\usepackage{mathrsfs}
\usepackage[margin=1in]{geometry}
\usepackage{tabularx}
\usepackage{makeidx}
\makeindex
\usepackage{graphicx}
\usepackage[utf8]{inputenc}

\newcommand{\ignore}[1]{}

\usepackage[usenames,dvipsnames]{xcolor}
\def\alert#1{\textcolor{Magenta}{#1}}

\newtheorem{thm}{Theorem}[section]
\newtheorem{quest}[thm]{Question}
\newtheorem{conj}[thm]{Conjecture}
\newtheorem{lemma}[thm]{Lemma}
\newtheorem{proposition}[thm]{Proposition}
\newtheorem{corollary}[thm]{Corollary}
\theoremstyle{remark}\newtheorem{rem}[thm]{Remark}
\theoremstyle{definition}\newtheorem{definition}{Definition}[section]
\theoremstyle{definition}\newtheorem{ex}[thm]{Example}

\newcommand{\empha}[1]{{\it #1}}
\newcommand{\genN}{\langle N \rangle}
\newcommand{\gencN}{\langle \cN \rangle}

\newcommand{\R}{\mathbb{R}}
\newcommand{\p}{\mathbb{P}}
\newcommand{\E}{\mathbb{E}}
\newcommand{\Z}{\mathbb{Z}}
\newcommand{\N}{\mathbb{N}}
\newcommand{\D}{\mathbb{D}}
\newcommand{\Cov}{\mathrm{Cov}}
\newcommand{\Corr}{\mathrm{Corr}}
\newcommand{\Var}{\mathrm{Var}}

\newcommand{\B}{\mathcal{B}}

\newcommand{\cC}{\mathcal{C}}
\newcommand{\cG}{\mathcal{G}}
\newcommand{\cF}{\mathcal{F}}
\newcommand{\cU}{\mathcal{U}}
\newcommand{\cE}{\mathcal{E}}
\newcommand{\cN}{\mathcal{N}}

\newcommand{\T}{\mathbb{T}}

\newcommand{\floor}[1]{\lfloor #1 \rfloor}
\def\1{1\!\! 1}

\def\lora{\longrightarrow}

\def\clue{\mathsf{clue}}

\def\Corr{\mathrm{Corr}}
\def\Spec{\mathscr{S}}

\def\eps{\epsilon}

\def \Maj{\mathsf{Maj}}

\newcommand{\Vol}{\mathsf{Vol}}

\def\sign{\mathrm{sign}}

\def\Su{\chi}

\def\HH{\mathsf{H}}

\numberwithin{equation}{section}
\numberwithin{figure}{section}
\def\bl{\begin{lemma}}
\def\el{\end{lemma}}
\def\bth{\begin{thm}}
\def\eth{\end{thm}}
\def\bc{\begin{corollary}}
\def\ec{\end{corollary}}
\def\bcj{\begin{conj}}
\def\ecj{\end{conj}}
\def\bpr{\begin{proposition}}
\def\epr{\end{proposition}}
\def\bde{\begin{definition}}
\def\ede{\end{definition}}
\def\beq{\begin{equation}}
\def\eeq{\end{equation}}
\def\bpf{\begin{proof}}
\def\epf{\end{proof}}
\newcommand{\comm}[1]{}

\makeatletter 
\newcommand\mynobreakpar{\par\nobreak\@afterheading} 
\makeatother
\newenvironment{myitemize}{\mynobreakpar\begin{itemize}}{\end{itemize}}
\def\bit#1{
\ \vskip 4 pt 
\begin{myitemize}
\setlength{\parskip}{-2pt}
#1
\end{myitemize}
\vskip -7 pt
\nobreak
}

\begin{document}

\title{Sparse reconstruction in spin systems II:\\
Ising and other factor of IID measures}
\author{P\'al Galicza
\and 
G\'abor Pete
}
\date{}
\maketitle
\begin{abstract}
For a sequence of Boolean functions $f_n : \{-1,1\}^{V_n} \longrightarrow \{-1,1\}$, with random input given by some probability measure $\p_n$, we say that there is sparse reconstruction for $f_n$ if there is a sequence of subsets $U_n \subseteq V_n$ of coordinates satisfying $|U_n| = o(|V_n|)$  such that knowing the spins in $U_n$ gives us a non-vanishing amount of information about the value of $f_n$.

In the first part of this work, we showed that if the $\p_n$s are product measures, then no sparse reconstruction is possible for any sequence of transitive functions. In this sequel, we consider spin systems that are relatives of IID measures in one way or another, with our main focus being on the Ising model on finite transitive graphs or exhaustions of lattices. We prove that no sparse reconstruction is possible for the entire high temperature regime on Euclidean boxes and the Curie-Weiss model, while sparse reconstruction for the majority function of the spins is possible in the critical and low temperature regimes. We give quantitative bounds for two-dimensional boxes and the Curie-Weiss model, sharp in the latter case. 

The proofs employ several different methods, including factor of IID and FK random cluster representations, strong spatial mixing, a generalization of discrete Fourier analysis to Divide-and-Color models, and entropy inequalities.



\end{abstract}

\tableofcontents

\section{Introduction and main results}

\subsection{Background and motivation}

In \cite{GaPe} we introduced the concepts of \empha{clue} and \empha{sparse reconstruction}. Central to this work as well, we start by defining these notions. 
Consider the product space  $(\Omega^V,\,\p)$ where $V$ is some large finite set, and $\p$ is not necessarily a product measure, and a function $f: \Omega^V \lora \R$, often the indicator function of an event. The  $\clue$ measures the amount of information the restriction $X_U := \{X_v : v \in U\}$ of the full input $X_V \in \Omega^V$ to a subset $U\subseteq V$ has about the output random variable  $Z=f(X_V)$.
 
\bde[$L^2$-clue] \label{clue}
Let $f: (\Omega^V,\,\p) \lora \R$ and $U \subseteq V$. Then,
\beq
 \clue(f\,|\,U) :=  \frac{\Var(\E[f \,|\, \mathcal{F}_U])}{\Var (f)} =   \frac{\Cov(f, \E[f \,|\, \mathcal{F}_U])}{\Var (f)}  = \Corr^2(f, \E[f \,|\, \mathcal{F}_U])\,,
\eeq
where $\cF_U$ is the sigma-algebra generated by the variables $\{X_v : v \in U\}$. 
When $f$ is constant, so that $\Var(f)=0$, we say that $\clue(f\,|\,U)=0$ for any $U\subset V$.
\ede

Depending on the measure $\p$ and the function $f$, when is it possible that knowing $X_U$ for a small but carefully chosen subset $U$ (specified in advance, in a non-adaptive manner, independently of the values of the variables) will give enough information to estimate $f(X_V)$? The concept of sparse reconstruction formalizes this question:

\bde[Sparse reconstruction of functions]
Consider a sequence $f_n: (\Omega^{V_n},\,\p_n) \longrightarrow \R$. We say that there is \empha{sparse reconstruction} for $\{f_n\}$ w.r.t.~$\{\p_n  \}_{n \in \N}$ if there is a sequence of subsets  $U_n \subseteq V_n$  with $\lim_{n} \frac{\left|U_n\right|}{\left|V_n\right|}=0 $ such that
$$
 \liminf_{n}{ \clue(f_n\,|\,U_n)}>0\,.
$$
\ede

In the first part of this work, among other results, we proved the following. This statement is a slight generalization of \cite[Theorem 2.7]{GaPe}, but the proof is identical.

\bth[Clue of quasi-transitive functions for product measures \cite{GaPe}]\label{t.cluegen}
Let  $\Gamma \leq \mathrm{Sym}_V$ be a group acting on a finite set $V$, and consider a product probability measure space $\left(\Omega^V, \bigotimes_{v\in V} \pi_v \right)$ and a function $f \in L^2 \left(\Omega^V, \bigotimes_{v\in V} \pi_v \right)$ that are both invariant under the action of $\Gamma$. Then, for any $U \subseteq V$,
$$ 
\clue(f\,|\,U)\leq \ \frac{\left|U\right|}{\min_{1 \leq j \leq n }{|\Gamma \cdot j|}}.
$$ 
In particular, if the measure space is IID and $f$ is transitive, then
$$  \clue(f\,|\,U)\leq \frac{\left|U\right|}{n}.  $$  
\eth

Thus, sparse reconstruction for sequences of  transitive functions of IID variables is not possible. There is also a more general version for arbitrary functions, using random subsets $\cU$ with low revealment; see Theorem~\ref{t.cluerandom} below. 

In the present paper, we investigate what happens if we replace the product measure in Theorems~\ref{t.cluegen} and~\ref{t.cluerandom} with some different sequence of measures; our key examples will be Ising models on different transitive graphs. The general idea is that if, in a sequence of measures, there is not too much dependence, i.e., they are, in some sense, close to product measures (such as ``subcritical'', high temperature Ising models), then there is no sparse reconstruction; on the other hand, if the measures admit lots of dependencies, such as critical and supercritical (low temperature) Ising models, sparse reconstruction for some transitive function is possible.

Our spin systems, when defined on bounded degree graphs, will always converge locally to a spin system on an infinite transitive graph: via exhaustions or in the Benjamini-Schramm sense; sometimes they will also have a non-trivial scaling limit; in either case, our goal will typically be to connect properties of the limit system to sparse reconstruction questions in the finite systems. We will always use some property of the system that is related to independence (these notions will of course be defined in due course):
\bit{
\item ergodicity, tail-triviality \cite{LanRue};
\item being a finitary factor of IID measure, with a good control on the decay of the coding volume \cite{BS,Spinka.SSM}; 
\item being a divide-and-colour measure \cite{StT}, such as the FK random cluster representation of the Ising model, with a good control on the decay of the cluster sizes;
\item having the strong spatial mixing property \cite{MO1,DSS};
\item having a large entropy, very close to having independent bits, related to the ``propagation of chaos'' phenomenon \cite{chaos1}.
}

As explained in \cite[Section 1]{GaPe}, for IID spins, a small clue (no sparse reconstruction) theorem can be considered as a baby noise sensitivity result; see \cite{BKS,SS,GPS,OD,GS,TV} for background on noise sensitivity. Since no strategy for proving noise sensitivity has been extended yet to non-IID spins, it seems rather interesting to extend at least our Theorems~\ref{t.cluegen} and~\ref{t.cluerandom}. 

Another motivation is that the OSSS-technology \cite{OSSS} of low-revealment adaptive algorithms has been successfully applied to a large variety of spin systems \cite{DCRT, H:new, MV, DCGRS,AABT} to prove sharp threshold results, hence a natural question is the existence of sparse reconstruction: when do low-revealment {\it non-adaptive} algorithms exists?

A very recent development, related to both previous paragraphs, is that \cite{KLMM} shows how good mixing, i.e., a log-Sobolev inequality for the Glauber dynamics of a Markov random field, can imply KKL and Talagrand-type inequalities. Applications include sharp threshold results and large influence of certain small coalitions, not in the direction of our questions.

Finally, as we will see, small clue results are related to how entropy-factorization and strong spatial mixing may imply good bounds on the spectral gaps of Glauber block dynamics \cite{BCCPSV}, but optimized in a direction that we have not seen in the literature so far.

In the rest of this Introduction, we will first introduce in Subsection~\ref{ss.different} several natural versions of sparse reconstruction, and summarize their relationships to each other. In Subsection~\ref{ss.positive}, we give an overview of our positive results: when sparse reconstruction is possible. This is the easier direction, typically just proving that total magnetization can be approximated by the sum of spins over a sparse subset. In Subsection~\ref{ss.negative}, we list our main negative results: different conditions that imply that sparse subsets have small clue regarding \empha{any} function. All along, Ising models will be our key examples, introduced carefully in Subsection~\ref{ss.Isingintro}, but many of our results apply in greater generality.

\subsection{Different notions of sparse reconstruction}\label{ss.different}

We start by introducing a few different concepts of sparse reconstruction. For product measures these are all equivalent (all of them fails),  but in general, they are not.

\bde[Sparse Reconstruction] \label{d.SR}
Let   $\{(\Omega^{V_n},\,\p_n) : n \in \N  \}$  be a sequence of probability spaces where $\p_n$ is invariant under some $\Gamma_n\leq \mathrm{Sym}_{V_n}$ acting on $V_n$ transitively.

There is \empha{sparse reconstruction} (briefly: SR) for  $\{(\Omega^{V_n},\,\p_n) : n \in \N  \}$  if there exists a sequence of transitive functions $f_n: (\Omega^{V_n},\,\p_n) \longrightarrow \R$ such that $f_n$ is $\Gamma_n$-invariant and there is sparse reconstruction for $\{f_n  \}_{n \in \N}$.

There is \empha{non-degenerate sparse reconstruction} (briefly: non-deg SR) for the sequence $\p_n$ if there exists a non-degenerate sequence of $\Gamma_n$-invariant Boolean functions $\{f_n  \}_{n \in \N}$  (i.e., with non-vanishing variance) such that there is sparse reconstruction for $f_n$.
\ede

In Example~\ref{SR_nondSR} we give an example of a sequence of measures for which there is SR, but no non-deg SR. Indeed, it may be the case that while we have $\Corr^2(f_n, \E[f_n \,|\, \mathcal{F}_{U_n}]) \to 1$ for some   $\{f_n  \}_{n \in \N}$ , the correlation is dominated by a small probability event where $f_n$ is abnormally large, and with high probability we have no good guess on the value of $f_n$ (or any other function). This motivates the introduction of the concept of non-deg SR.

At the same time, we will show (see Corollary \ref{c.Wsr_sr})  that, whenever there is SR,  there is also a (possibly degenerate) sequence of Boolean functions for which there exists sparse reconstruction. That is, the difference between SR and non-deg SR lies in the fact that for the latter case we insist that the sequence of Boolean functions in question  be non-degenerate.
\medskip

There is an alternative version of sparse reconstruction, already featured in \cite[Theorem 1.2]{GaPe}, which does not require any symmetry of $\p_n$ and $f_n$:

\bde[Random Sparse Reconstruction] \label{RSR}
Let   $\{(\Omega^{V_n},\,\p_n): n \in \N  \}$  be a sequence of probability spaces. For every $n$, let $\mathcal{U}_n$ be a random subset of $V_n$, independent of the spin system $(\Omega^{V_n},\,\p_n)$, with the property that
$$
\delta_n:= \max_{j \in V_n}{\p[j \in \mathcal{U}_n]} \rightarrow 0.
$$
The quantity $\delta_n$ is called the revealment of $\mathcal{U}_n$.

There is \empha{random sparse reconstruction} (briefly, RSR) for   $\{(\Omega^{V_n},\,\p_n): n \in \N  \}$ if there is a sequence of Boolean functions $f_n: \{-1,1\}^{V_n} \longrightarrow \{-1,1\}$ and a  $\mathcal{U}_n$ as above, such that
$$
\E[\clue(f_n \; |\;\mathcal{U}_n)]> c
$$
for some $c>0$.
\ede

\bth[No RSR in product measures \cite{GaPe}]\label{t.cluerandom}
Let $f\in L^2(\Omega^{V},\pi^{\otimes V})$ be any function. Let $\mathcal{U}$ be a random subset of $V$, independent of the $\sigma$-algebra $\mathcal{F}_{V}$. Then 
$$
\E[\clue(f \; |\;\mathcal{U})] \leq  \delta(\mathcal{U}),
$$
where $\delta(\mathcal{U}) := \max_{j \in V}{\p[j \in \mathcal{U}]}$ is the revealment of $\mathcal{U}$.
\eth

For each type of sparse reconstruction, there may also be \empha{full reconstruction} (abbreviated as full SR, full non-deg SR and full RSR). By this we mean that
$$
 \lim_{n}{ \clue(f_n\,|\,U_n)}=1\,,
$$
or, in the case of RSR, $\lim_{n}{\E[\clue(f_n \; |\;\mathcal{U}_n)]} =1$.

It turns out that for each case of sparse reconstruction (SR, RSR, non-deg SR, non-deg RSR), the partial and full versions are equivalent, for any given spin system --- see Theorem~\ref{clueto1}. This is a highly non-trivial statement, especially for the non-deg cases. The reason for the non-trivality is that our natural idea to boost partial to full reconstruction is to take a sum of our translated or repeated guesses, but this destroys Booleanity, while the simplest Booleanization, going from the sum to majority, does not always work (see Example~\ref{SR_nondSR}, already mentioned above). Hence, in proving full non-deg reconstruction, one has to use crucially that partial reconstruction already held for a Boolean function. 
\medskip

There are several reasons why RSR, compared to SR, is of interest. First, it allows us to investigate sparse reconstruction for sequences of arbitrary (non-transitive) functions, and to investigate spin systems that do not exhibit symmetries. A natural example for the latter is the sequence of Ising models on the boxes $[0,n]^d$ with free boundary conditions. Another interesting example is the sequence $G_{n,d}$ of random $d$-regular graphs on $n$ vertices, which is known to converge to the $d$-regular infinite tree $\T_d$ in the \empha{Benjamini-Schramm (or local weak) sense}; see Section~\ref{basicconc} below, or \cite[Section 14.1]{PGG}, for instance. So, we may consider a Benjamini-Schramm approximation of a translation-invariant measure on $\T_d$ by approximating  measures in $G_{n,d}$; for instance, using finitary factor of IID measures (see Definition \ref{FFIIDdef}).

Another important point is that RSR can be directly applied to a measure $\p$ supported on an infinite spin system. Let $V$ be the set of vertices and let $\{V_n  \}_{n \in \N}$ be a monotone exhaustion of $V$. Let $\p_n$ be the projection of $\p$ onto $V_n$. We say that there is \empha{RSR via an exhaustion} for $(\Omega^{V},\,\p)$ if there exists an exhaustion $\{V_n\}$ such that there is RSR for the sequence  $\{(\Omega^{V_n},\,\p_n): n \in \N \}$. 

Equivalently, let  $\{f_n: \{-1,1\}^{V} \longrightarrow \{-1,1\} : n \in \N  \}$ be a sequence of finitely supported Boolean functions (that is, each function $f_n$ depends only on a finite set of coordinates). There is RSR for  $\{(f_n): n \in \N  \}$ on $(V, \p)$ if there exists a sequence of random subsets $\mathcal{U}_n \subseteq V$ with revealment $\delta_n \to 0$, such that for all $n \in \N$ we have $\E[\clue(f_n \; |\;\mathcal{U}_n)]> c$ for some $c>0$.

However, it is important to understand that it may be crucial if a spin system on an infinite graph is approximated in the Benjamini-Schramm sense or by an exhaustion. See for instance, Example \ref{e.devanSR}, or Question~\ref{q.free} for another possible example: the free Ising measure on $\T_d$.

Finally, yet another reason that makes RSR interesting is that it can be directly linked with the mixing properties of certain Glauber block dynamics (see Section~\ref{s.glauber}). 
\medskip

We now discuss the connection between SR and RSR. 

\bl[SR and RSR are equivalent for transitive functions]\label{SReqRSR}
For a sequence $\{(\Omega^{V_n},\,\p_n): n \in \N  \}$, where $\p_n$ is invariant under a group $\Gamma_n$ acting transitively on $V_n$, and $f_n$ is any sequence of $\Gamma_n$-invariant functions, there is SR if{f} there is RSR.

Furthermore, the same equivalence holds between full SR and full RSR.
\el

\bpf
First, assume that  SR holds for some transitive $f_n$. Then, from the deterministic set $U_n$ giving $\clue(f_n\,|\,U_n)>c$, we can define $\mathcal{U}_n$ as a uniformly random $\Gamma_n$-translate of $U_n$. Because of the invariance of $\p_n$ and $f_n$, the clue is invariant under automorphisms, while the revealment of $\cU_n$ is $\delta_n = \frac{\left|U_n\right|}{\left|V_n\right|}\to 0$. This also shows that full SR implies full RSR.

On the other hand, if there is RSR for a sequence $f_n$ of transitive functions, then, using that the clue can be at most $1$, from the expected clue we get that $\p[\clue(f_n \; |\;\mathcal{U}_n) >c]>c$ for some $c >0$.  At the same time, for any fixed $\eps>0$ and every large enough $n$, we have $\p\big[ |\mathcal{U}_n| > \eps |V_n| \big] <\eps$, using that $\E\, |\mathcal{U}_n| = o(|V_n|)$. So, 
$$
\p\big[\clue(f_n \; |\;\mathcal{U}_n) >c \text{ and }  |\mathcal{U}_n| \leq \eps |V_n|\big] 
\geq 1 - \p\big[\clue(f_n \; |\;\mathcal{U}_n) \leq c\big] - \p\big[ |\mathcal{U}_n| > \eps |V_n|\big]  \geq c-\eps,
$$
which is uniformly positive if $\eps$ is smaller than $c/2$, say. That is, there exist small sets with non-vanishing clue, showing SR. Moreover, if we have full RSR for $f_n$, then we can take $c$ arbitrarily close to 1, hence we also have full SR.
\epf

All of this means that there is possibly a difference between SR and RSR only for non-transitive functions. Note, however, that proving for a transitive function that there is RSR from some random $\cU_n$, hence obtaining the {\it existence} of a deterministic $U_n$, may not be completely satisfactory: without an explicit example of a $U_n$, this is just ``non-constructive sparse reconstruction''.  Hence it is of additional value to show an explicit sparse subset $U_n$ --- this will be the case, e.g., in our Theorem~\ref{critical_is} for the critical planar Ising model.

In the following example there is RSR for the value of a single bit (a highly non-transitive function), and also for a quasi-transitive function, but there is no SR.

\begin{ex}[A sequence of measures with RSR, but no SR] \label{SRvsRSR} 
  Let  $G_n = \Z_n\times \Z_2$; think of a ``top'' and a ``bottom'' circle $\Z_n$. Consider the following measure, invariant under the automorphism group of the graph $G_n$. First flip a fair coin; if it is heads, then the $n$ spins in the top circle are IID Bernoulli($\frac{3}{4}$) variables, while the bottom $n$ spins are IID Bernoulli($\frac{1}{4}$). If it is tails, the Bernoulli parameters are flipped. 
  
First, there is no SR in this system. Indeed, take a sparse sequence of subsets $U_n \subseteq V(G_n)$ and consider the $\sigma$-algebra $\mathcal{F}_n$ generated by $\mathcal{F}_{U_n}$ and the fair coin flip. For any sequence $f_n$ of transitive functions, we have $\clue(f_n \; |\;\mathcal{F}_{n}) \to 0$ by Theorem~\ref{t.cluegen}, since, conditioned on the fair coin, the spin system is distributed as a product measure, invariant under the $\Z_n$-rotation, acting with index 2. Then $ \clue(f_n \; |\;U_n) \leq \clue(f_n \; |\;\mathcal{F}_{n})$ implies that SR is not possible.

At the same time, a single bit can be (partially) reconstructed from a random set of vanishing density. Let $ 1/n \ll \delta_n \ll 1$ and let $\mathcal{U}_n$ be a random Bernoulli subset of $V(G_n)$ with density $\delta_n$. With probability tending to 1, the number of elements of $\mathcal{U}_n$, both in the top and in the bottom circle, tends to infinity. Then, by the weak law of large numbers, the probability that the average in the $\mathrm{Ber}(3/4)$ (or the $\mathrm{Ber}(1/4)$) side is more than (respectively less than) $1/2$ is tending to 1. Thus we can reconstruct the fair coin flip with probability tending to 1, and we can correctly guess the value of any given spin with probability at least $3/4-o(1)$. Since the original marginal of a single spin is just a fair coin flip, this is a guess with non-vanishing clue. Moreover, we have clue $1-o(1)$ about the index 2 quasi-transitive balanced Boolean function that describes whether there are more 1's in the top circle than in the bottom one.
\end{ex}

We do not know under what general conditions could SR and RSR be equivalent. See Question~\ref{q.SRvsRSR}.
\medskip

Next, the idea of trying to reconstruct some information in a spin system from sparse data might remind the Reader of the notion of {tail-triviality}. Namely, on an infinite graph $G(V,E)$, a probability measure $\mu$ on some $\Omega^V$ is called \empha{tail-trivial} if, for one or any exhaustion $V_n\nearrow V$ by finite subsets, the $\mu$-measure on any event that is measurable in the variables indexed by $V_n^c$ for every $n$, is either 0 or 1. For instance, \empha{Kolmogorov's 0-1 law} says (e.g., \cite[Theorem 9.20]{PGG}) that any product measure over $V$ is tail-trivial. The notion is easily seen to be equivalent to the following non-reconstruction property: for every $\eps>0$, if $A$ is an event depending on a finite number of spins $F\subset V$ with probability in $(\eps,1-\eps)$, then, if $n$ is large enough, and the distance of $U\subset V$ is at least $n$ from $F$, then the correlation between $A$ and any $U$-measurable event $B$ is at most $\eps$. Furthermore, a basic result connecting tail-triviality to ergodic theory is the \empha{Lanford-Ruelle lemma} \cite{LanRue}: a measure with the spatial Markov property is tail-trivial if{f} it is extremal in the set of all probability measures with the same local specifications (the conditional measures given the configuration on the complement of a finite set). Thus, for instance, if we have a unique spatial Markov measure for some local specifications, then it is tail-trivial, and hence a single spin (or any function depending on finitely many spins) cannot be reconstructed from far away.

Connecting this to our questions, if we are considering spin systems on finite graphs $G_n$, converging to a tail-trivial spin system on an infinite transitive graph $G$ either in the Benjamini-Schramm sense or via exhaustions, and there is a sequence of functions $f_n: \Omega^{V_n}\lora\R$ that are ``local'' in the sense that they converge to a measurable function $f: \Omega^V\lora \R$ in the obvious sense, then there is no reconstruction for these functions from far away, and in particular, there is no sparse reconstruction for them. It is also natural to consider sparse relaxations of tail-triviality on infinite graphs, as follows. These are small modifications of notions considered by Bal\'azs Szegedy \cite{Szy}.

\bde\label{d.stber}
Given $(\Omega^V,\p)$, a function $f: \Omega^V\lora \R$ is \empha{sparse Bernoulli reconstructable} if, for any $\eps>0$, letting $\cU_\eps$ be the random subset of $V$ that contains each element with probability $\eps$, independently of everything else, we have $\clue(f \,|\, \cU_\eps)=1$ almost surely.
\ede 

\bde\label{d.stgen}
Given $(\Omega^V,\p)$, a function $f: \Omega^V\lora \R$ is \empha{sparse reconstructable} (or \empha{sparse tail measurable}) if there exists a decreasing sequence of subsets $V\supset U_1 \supset U_2 \supset \dots$ such that the smallest distance between any two distinct elements of $U_n$ goes to infinity, but $\clue(f \,|\, \cU_n)=1$ for every $n$.
\ede 

Szegedy's motivation to consider sparse tail measurability was the following. While product measures are tail-trivial by Kolmogorov's 0-1 law, factor of IID measures, somewhat surprisingly, can be far from being tail trivial; e.g., the unique translation invariant measure on perfect matchings of the 3-regular tree $\T_3$ was shown to be an FIID measure in \cite{LyNaz}, while the tail sigma-algebra is easily seen to be the entire sigma-algebra. Nevertheless, he conjectured \cite{Szy} that any FIID measure is in fact sparse tail trivial: there exist no sparse (Bernoulli) reconstructable functions. Supporting his conjecture, he proved using entropy that this indeed holds for FIID measures on amenable transitive graphs. We will present a proof of this (so far unpublished) result in Subsection~\ref{FFIID.tt}, along with further non-reconstruction results for FIID measures.

The moral of this discussion is that, for many natural measures (for tail-trivial measures for sure, for FIID measures conjecturally), local functions will not yield sparse reconstruction, hence one should look at global functions such as magnetization. This brings us to the first half of our results: when sparse reconstruction is possible.

\subsection{Critical and super-critical (low temperature) models: reconstruction results}\label{ss.positive}

If a sequence of spin systems $\{ (G_n, \p_n) : n \in \N  \}$ is super-critical in the sense that spontaneous symmetry-breaking happens, the limit measure is non-ergodic, then it seems intuitively obvious that sparse reconstruction should be possible: a few spins already should reveal which limiting ergodic component (which pure phase) the finite system is close to. However, this is false in this generality, since phases can coexist on the same finite graph; see Example~\ref{nonerglim}. Nevertheless, if the finite graphs are bounded degree expanders (see Definition~\ref{d.expander}), then such phase co-existence is impossible, and we get the following result (see Section~\ref{basicconc}):

\bpr \label{SRnonerg0}
Let  $\{ (G_n, \p_n) : n \in \N  \}$ be a sequence of spin systems, where $\{ G_n  : n \in \N \}$ is an expander sequence with bounded degrees converging in the Benjamini-Schramm sense to a limiting spin system  $(G, \p)$, with $G$ transitive.

If the $\mathsf{Aut}(G)$-invariant measure $\p$ is not $\mathsf{Aut}(G)$-ergodic, then there is (Random) Sparse Reconstruction for $\{ (G_n, \p_n) : n \in \N  \}$.
\epr

Here, $\mathsf{Aut}(G)$ stands for the automorphism group of the graph $G$, and a measure $\p$ is called $\Gamma$-ergodic for some $\Gamma\leq \mathsf{Aut}(G)$ if every $\Gamma$-invariant event has $\p$-measure $0$ or $1$.

If the systems are super-critical in the sense of long-range order, or at least diverging susceptibility, we can also treat general graphs. Namely, assume that the finite graphs $G_n(V_n,E_n)$ are transitive, the measures $\p_n$ are defined on $\{-1,1\}^{V_n}$, and that the 
\empha{susceptibility} blows up:
\begin{equation}\label{e.suscfirst}
\Su(\sigma^n) := \sum_{v\in V_n} \Cov(\sigma^n_o,\sigma^n_v) \to \infty\, \textrm{ as }n\to\infty.
\end{equation}
Since the graphs are transitive, this is equivalent to the \empha{total average magnetization}
$$
M_n(\sigma^n) := \frac{1}{|V_n|} \sum_{ v\in V_n} \sigma^n_v
$$
having a standard deviation much larger than $1 / \sqrt{|V_n|}$. See Subsections~\ref{ss.Isingintro} and~\ref{ss.magn} for more details on susceptibility, especially for the Ising model on which we will now focus on. 

The following theorem is proved in Subsection~\ref{ss.supIsing}, using the general results of Section~\ref{General_SR}.

\bth[Critical and supercritical sparse reconstruction]\label{supercrit_is}
Let $\{ \sigma^n \}_{n \in \N}$ be the sequence of Ising models on a sequence of transitive graphs $\{G_n \}_{n \in \N}$, at inverse temperatures $\beta_n$ such that the susceptibility blows up: $\Su_{\beta_n}\to\infty$. 
\bit{
\item[{\bf (1)}] There is full sparse reconstruction for the total average magnetisation $M_n$ and sparse reconstruction for $\Maj_n:=\sign(M_n)$.
\item[{\bf (2)}] Assume, in addition, that the corresponding random cluster model FK$(p_n,q=2)$ on $G_n$ has a unique giant  cluster: there exists a constant $\lambda>0$ such that, for every $L>0$, the probability that the largest cluster has size at least $\lambda |V_n|$ and the second largest cluster has size at most $|V_n|/L$ goes to 1 as $n\to\infty$. Then there is full  sparse reconstruction for $\Maj_n$.
}
In both cases, the clue of an IID Bernoulli$(p_n)$ random subset tends to 1 with probability tending to 1, provided that $p_n \Su_{\beta_n} \to \infty$.
\bit{
\item[{\bf (3)}]  For the tori $G_n = \Z_n^d$, case (1) holds for $\beta_n \geq \beta_c(\Z^d)$, and (2) holds for $\beta_n\ge \beta > \beta_c(\Z^d)$.
\item[{\bf (4)}]  For the Curie-Weiss model on the complete graphs $G_n=K_n$,  for $\beta >\beta_c=1$, if $|U_n|\to\infty$, then $\clue(M_n\,|\,U_n)\to 1$ and $\clue(\Maj_n\,|\,U_n)\to 1$. For $\beta=\beta_c=1$, we have both full reconstruction results whenever $|U_n| \gg \sqrt{n}$.
}
\eth

For the lattices $\Z^d$, one shortcoming of these results is that the sparse subset from which we have reconstruction is not explicit. Another shortcoming is that we do not get full non-degenerate reconstruction at the critical temperature.  For the two-dimensional critical Ising model, using state-of-the-art scaling limit results \cite{CHI,CGN,CGN2}, we can fill these gaps in. The proof written up in Subsection~\ref{ss.planar} is essentially due to Christophe Garban.

\bth[Sparse reconstruction for critical planar Ising] \label{critical_is}
Consider the sequence of Ising measures with $\beta=\beta_c$ on the boxes $[0,n]^2$, with the free or the $+$ boundary condition. Then the total average magnetization $M_n$ can be reconstructed with high clue from a sublattice $H_n$ of mesh size $s_n$ as long as $s_n = o(n^{\frac{7}{8}})$. That is, 
$$
\lim_{n \to \infty}{\clue(M_n\,|\,H_n)} =1.
$$
Moreover,
$$
\lim_{n \to \infty}{\clue(\Maj_n\,|\,H_n)}=1.
$$
\eth

One may wonder if the above bound is optimal or close to optimal in terms of the size of the subset. In other words, how small a sequence of subsets has to be so that no sparse reconstruction is possible? As we shall see, Corollary~\ref{cor.FK}~(2) says that for $|H_n| \ll n^{1/8}$, so, for mesh size $s_n \gg n^{15/16}$, we have $\clue(\Maj_n\,|\,H_n) \to 0$.

It may appear from the above results that magnetization is the only global function that is worth looking at. But this does not seem to be the case. A simple-minded example is the critical Ising model on the tori $\Z^2_n$, conditioned on $M_n=0$, which we think does have RSR, but it is clear that it is not $M_n$ that can be reconstructed. See Question~\ref{q.robust} and the discussion after that.

A more natural example is the low temperature $+$ Ising measure, i.e., the Ising measure conditioned to have positive magnetization. In the statistical physics community this is usually regarded to be similar to the high temperature model. For instance, susceptibility is finite, and the free energy is an analytic function of $\beta>\beta_c(\Z^2)$ \cite{Ott}. However, there are some important differences. In contrast with high temperature models, the magnetization has a subexponential tail (because the Ising model is a Markov field, hence we only need to pay the price of the negative spins on the boundary of a large  cluster of minuses) and therefore it is not a finitary FIID measure \cite[Theorem 2.1]{BS}. Furthermore, it was proven in \cite{BM}, for not unrelated reasons, that the single site Glauber dynamics for the low temperature $+$ Ising measure on boxes of $\Z^2$ have no spectral gap. We expect that this model does have RSR; see Question~\ref{q.plus}. 

One more example where the function with RSR may be something quite different from total magnetization is the uniform random perfect matching on large girth graphs; see Question~\ref{q.pm}.

\subsection{Critical and sub-critical (high temperature) models: small clue results}\label{ss.negative}

A natural first idea to extend the IID small clue results of \cite{GaPe} is to consider \empha{factor of iid (FIID)} spin systems $(\sigma_v)_{v\in V(G)}$. This means that there is an underlying field of IID variables $(\omega_v)_{v\in V(G)}$, and then each spin $\sigma_v$ is given by applying the same measurable \empha{coding function} $\psi$ to the field $\omega$ viewed from the vertex $v$ as a root. See Subsection~\ref{FFIID.intro} for the precise definition and some examples. A usual way to strengthen this definition and make it even closer to IID measures is to insist that $\psi$ be not only measurable, but also depend on a bounded neighbourhood of the root (called a \empha{block factor}). In between are the \empha{finitary factors (FFIID)}, where $\psi$ is given by an algorithm that almost surely terminates after looking at some finitely many $\omega$-variables; i.e., for each $v\in V$, there is a finite coding radius $R_v$, measurable w.r.t.~the $\omega$-variables, such that $\sigma_v$ is determined by $\omega$ inside the ball $B_{R_v}(v)$ around $v$. See Definition~\ref{FFIIDdef}. 

A nice feature of FIID spin systems on an infinite graph $G$ is that they can naturally be approximated by spin systems on any sequence of finite graphs $G_n$ converging to $G$ in the Benjamini-Schramm sense, and then sparse reconstruction questions make perfect sense. One may think that if a FFIID system has \empha{finite expected coding volume}, then any sparse set of spins $(\sigma_v)_{v\in U}$ enquired to compute some function is determined by the IID variables $(\omega_v)_{v\in J}$ in some $J$, whose expected size is just a constant multiple of $|U|$, hence it is still sparse with high probability, and the earlier IID small clue results should apply. However, this sparse set $J$ now already depends on $\omega$, and we know that there are functions of IID bits that can be determined by \empha{adaptive} low revealment algorithms, hence we do not see how to complete this argument. We do nevertheless have some results in this direction. First, in Subsection~\ref{FFIID.intro}, we prove the following.

\bth[Almost no SR for exponentially decaying coding radius]\label{t.almost}
Let $\sigma$ be a finitary factor of IID spin system on $\Omega^{\Z^{d}}$ with the property that, for some $c>0$ and every $t>0$,
$$\p[R > t] <\exp(-ct),$$ 
where $R$ is the coding radius. Let $\sigma_n$ be any factor of IID approximation to $\sigma$ (in the sense of Definition~\ref{d.FIIDconv}) on the sequence $\Z_n^d$ of tori. Let $f_n: \Omega^{\Z^d_n} \longrightarrow \{-1,1\}$ be sequence of transitive, non-degenerate Boolean functions. Then, for any sequence of subsets $\{U_n \subseteq \Z^d_n\}_{n \in \N} $ with 
$$|U_n| \ll \left(\frac{n}{\log n}\right)^d,$$
we have $\clue(f_n \; | \;U_n) \to 0$.
\eth

This is sharp in the sense that in Example~\ref{e.devanSR} we give a sequence $(G_n,\p_n)$ converging to a FFIID measure on $\Z^2$ with exponentially decaying coding radius that does admit sparse reconstruction. However, in that example there exists a vertex in $G_n$ with probability tending to 1 for which the coding radius is larger than the diameter of $G_n$, resulting in a coding error, and it is basically this error that we can reconstruct. If we strengthen the convergence notion to exclude errors, then we do not have a counterexample; see Question~\ref{q.adam}. Supporting this conjecture, we have the following results, the two parts proved in Subsections~\ref{FFIID.tt} and~\ref{FFIID.magn}, respectively:

\bth[FFIID with ``robustly finite'' expected coding volume]\label{t.fvFFIID} 
Let $G(V,E)$ be an infinite transitive graph and $\p$ a finitary FIID measure on $\Omega^V$ with the coding radius $R(o)$ for $o\in V$ satisfying $\E |B_{2R(o)}|<\infty$.
\bit{
\item[{\bf (1)}] $\p$ is tail trivial.
\item[{\bf (2)}] Let  $\{(G_n, \p_n) \}$ be a sequence of FFIID spin systems converging to $(G, \p)$ in the FFIID sense of Definition~\ref{d.FIIDconv} or by exhaustions. Suppose there exists an $\alpha >0$ such that  $\Var(M_n) \geq \alpha/|V_n|$ for all $n \in \N$, where $M_n =\frac{1}{|V_n|}\sum_{j \in V_n}{\sigma_j}$ denotes average magnetization as usual. Then there is no (R)SR for magnetization.
}
\eth

Strong Spatial Mixing (see Definition~\ref{SSM} below; abbreviated as SSM) is a property for spin systems on $\Z^d$ that is stronger, closer to being IID-like, than being a FFIID measure with finite expected coding volume. (The fact that SSM implies the latter property is proved in \cite{Spinka.SSM}.) For such measures we will prove stronger results in Section~\ref{s.glauber}: the complete absence of sparse reconstruction.

\bth \label{t.norSR}
Let $\{ \mu_n \}$ be a sequence of measures on $S^{V_n}$ with finite state space $S$ satisfying the spatial Markov property  on the sequence of tori $V_n=\Z_n^d$ or on a sequence of rectangles $V_n=[0,n]^d$ such that the Strong Spatial Mixing property holds. Then there is no random sparse reconstruction on $\{ \mu_n \}$. 
\eth

This implies that there is no random sparse reconstruction for the Ising model on  $\Z_n^d$ for the entire high temperature regime. Indeed, it has been shown in \cite{DSS} that  the SSM condition holds for any $\beta <\beta_c$ for all dimensions. Note that, before \cite{DSS}, SSM was known to be equivalent to a uniform logarithmic Sobolev inequality for the Ising Glauber dynamics \cite{MOS}, satisfied for high enough temperatures for all $d\ge 2$ \cite{SZ,Marton}, and for all high temperature Ising measures for $d=2$ \cite{MOS}. Another example with SSM is the maximal entropy measure on proper $k$-colorings of $\Z^d$, with $k > 4d$; see \cite{PelSpi}.

For general graphs we have the following result for high temperature models:

\bth \label{t.norSRgen}
Let $\{ \mu_n \}_{n \in \N}$ be a sequence of Ising measures on a sequence $G_n(V_n, E_n)$ of finite graphs with maximum degree $d$ and with  $\beta < \beta_c(d)$, where $\beta_c(d)$ is the Tree Uniqueness Threshold (see \eqref{TUT}). Then there is no random sparse reconstruction for $\{ \mu_n \}_{n \in \N}$. 
\eth

The above two  results are obtained by understanding the mixing properties of some block Glauber dynamics, and the proofs  strongly rely on \cite{BCSV} and \cite{BCV}.
\medskip

The following proposition has a completely different proof, by extending the IID Fourier technology to Divide-and-Color (DaC) measures in Section~\ref{s.DaC}, given by partitioning the spins into clusters randomly, then flipping an independent fair coin for each cluster. 

\bpr\label{p.FK} Consider the Ising model $\sigma$ on any finite graph $G(V,E)$, at inverse temperature $\beta$, and total average magnetization $M(\sigma)$. In the associated random cluster model FK$(p,2)$, let $\cC_{v}$ be the cluster of vertex $v \in V$. Then, for any random set $\cU \subseteq V$ that is independent of $\sigma$, with revealment $\delta$, 
\beq\label{e.FKmagn}
\E\big[\clue(M\,|\,\cU)\big] \leq 
\delta\, \frac{\sum_{v \in V}{\E \left[|\mathcal{C}_{v}|^2\right]}}{\sum_{v \in V}{\E \left[|\mathcal{C}_{u}|\right]}}.
\eeq

More generally, if $f: \{-1,1\}^{V} \longrightarrow \R$ is odd, i.e., $f(-\omega) = -f(\omega)$, then
\beq\label{e.FKodd}
\E\big[\clue(f\,|\,\cU)\big]\leq 
\delta\, \E \left[\max_{v \in V}{|\mathcal{C}_{v}|}\right].
\eeq
\epr

The requirement that $f$ has to be odd is just a weakness of our method: we need to control the information that comes solely from knowing the cluster structure, which information is simply zero for odd functions. See the proof of Theorem \ref{DaCClueThm} for details.
  
Part~(1) of the next corollary is a direct consequence of the previous proposition. Part~(2) is a special case of the proposition, for the critical planar Ising measure, proved in Subsection~\ref{ss.planar}.

\bc\label{cor.FK}
Consider a sequence of Ising measures $\mu_n$ at some inverse temperatures $\beta_n$ on a sequence $G_n(V_n, E_n)$ of finite graphs.
\bit{
\item[{\bf (1)}] Assume that each $G_n$ is transitive, and the systems are subcritical in the sense that, in the random cluster FK representation, $\E_n \left[|\mathcal{C}_{u}|^2\right] / \E_n \left[|\mathcal{C}_{u}|\right]$ remains bounded. Then there is no random sparse reconstruction for the total average magnetization $M_n$. 
\item[{\bf (2)}]  Let $\{ \mu_n \}$ be the sequence of critical Ising measures (i.e., $\beta=\beta_c(\Z^2)$) on either tori $G_n=\Z_n^2$ or boxes $G_n=[0,n]^2$. Let $f_n: \{-1,1\}^{G_n} \longrightarrow \R$ be a sequence of odd functions. Then, for any sequence of random subsets $ \{ \mathcal{U}_n \subseteq  V_n\} _{ n \in \N}$ with 
$\delta(\mathcal{U}_n) \ll \frac{n^{\frac{1}{8}}}{n^2} = n^{-15/8}$, we have 
$$
\E[\clue(f_n\,|\,\cU_n)] \to 0.
$$
In particular this applies to total magnetisation.
}
\ec
 
So, for critical Ising on two-dimensional boxes, Theorem \ref{critical_is} says that we can find $|U_n| \gg n^{1/4}$ from which magnetisation can be reconstructed, while part~(2) of the above corollary says that for $|U_n| \ll n^{1/8}$ in tori, or for small revealment $\cU_n$ in boxes, the clue goes to $0$. See Question~\ref{q.IsingMplane}.

The condition in Corollary~\ref{cor.FK}~(1) is quite realistic. Hutchcroft \cite{H:new} proved that, for the random cluster FK$(p,2)$ model on any infinite transitive graph $G$, the volume has an exponential tail in the entire subcritical regime $p<p_c(2)$, hence $\E_{G,p,2} \left[|\mathcal{C}_{u}|^2\right] / \E_{G,p,2} \left[|\mathcal{C}_{u}|\right]$ is finite. One expects that if a sequence $G_n$ of transitive finite graphs converges to $G$ locally, then the random cluster models at the same $p$ value behave similarly. However, this locality question is very hard in general: even for Bernoulli percolation, $q=1$, it has been resolved only recently \cite{EaHu}. 
\medskip

For the Ising model on the complete graph $K_n$, called the \empha{Curie-Weiss model}, the notion of SSM makes no sense. The DaC approach does make sense, but it gives suboptimal results. However, with yet another completely different set of tools, using an information theoretic version of clue (Definition \ref{d.Infoclue}), we can give sharp results, matching the bound of Theorem~\ref{supercrit_is}~(4) for the critical case. The proof, given in Subsection~\ref{ss.CW}, uses Proposition~\ref{CWlargeent}, interesting in its own right, saying that the total $\log_2$-entropy of the system is $n-O_\beta(1)$ in the subcritical case, while $n-O(\sqrt{n})$ in the critical case. The elegant proof that we present for the proposition is due to Amir Dembo.

\bth \label{t.CWnoSR}
\bit{
\item[{\bf (1)}] 
There is no non-degenerate random sparse reconstruction for the high temperature Curie-Weiss model.
\item[{\bf (2)}] At $\beta = \beta_v$ for any non-degenerate function $f_n$ and any subset $|U_n| \ll \sqrt{n}$, we have $\clue(f_n \,|\, U_n) \to 0$.
 }
\eth

We will conclude the paper with a large number of open problems in Section~\ref{s.open}.

\subsection{The Ising and FK random cluster models}\label{ss.Isingintro}

As many of our results concerns the Ising model, we include here a brief introduction to the topic. The Ising model is one of the best understood and investigated model in statistical physics; for an introduction, see \cite[Section 13.1]{PGG} or \cite{FV}. Let $G(V,E)$ be a finite graph. The Ising measure on $G$ with inverse temperature $\beta>0$ is given by
$$
\p_{\beta, h}[\sigma ] : = \frac{1}{Z_{\beta}} e^{-\beta H(\sigma)},
$$
where $\sigma \in \{-1,1\}^{V}$ and the  so-called partition function $Z_{\beta}$ is a normalization factor that makes $\p_{\beta, h}$ a probability measure. The Hamiltonian $H:  \{-1,1\}^{V_n} \longrightarrow \R$  is defined as
$$
H(\sigma) = -J\sum_{(x,y) \in E}{\sigma_{x} \sigma_{y}} - h  \sum_{x \in V}{\sigma_{x}},
$$
for some $J >0$ and $h \in \R$. The parameter $J$ describes the strength of the interaction, while $h$ is the external magnetic field. Except for one of the definitions of susceptibility, formula~\eqref{issusc}, we will always assume $h=0$ in this paper.

One of the most remarkable features of the Ising model is that on many graph sequences it exhibits phase transition, that is, there are some critical inverse temperature values (denoted by $\beta_c$), where the behaviour of the model changes radically. 

The phase transition can be more clearly observed in an infinite model. The Ising measure on an infinite graph can be defined as a weak limit of a growing sequence of finite spanned subgraphs. The laws of finite models, however, may differ depending on what we do at the boundary of the finite subgraphs. One of the important features of the phase transition is, that typically, in particular on $\Z^d$, for $\beta > \beta_c$, one obtains different limiting measures if the outer boundary of the finite graphs is set to $-$ or $+$. On the other hand, for  $\beta \leq \beta_c$  there is a unique limiting Ising measure. The $\beta>\beta_c$ regime, due to the long range order, is usually called the supercritical phase, despite the temperature being low.

In the context of the Ising model, the average $M(\sigma) = \frac{1}{|V|}\sum_{x \in V}{\sigma_{x}}$ of spins is called the magnetization. Then we define
$$
m_{\beta,h} : = \E_{\beta,h}[M(\sigma)].
$$
Another important concept is susceptibility, defined as
\beq \label{issusc}
\Su_{\beta,h} := \frac{1}{\beta}\frac{\partial m_{\beta,h}}{\partial h}.
\eeq
The interpretation is that susceptibility expresses the reactivity of the magnetization to an external magnetic field. A short calculation gives that it can also be given in the following way: 
\beq\label{issusc2}
\begin{aligned}
	\Su_{\beta,h} &=  \frac{1}{|V|} \sum_{x,y \in V}{\Cov_{\beta,h}( \sigma_{x}, \sigma_{y})}\\
	&= \sum_{y \in V}{\Cov_{\beta,h}( \sigma_{x}, \sigma_{y})} \quad\text{for transitive graphs},\\
\end{aligned}
\eeq
in agreement with our earlier general definition \eqref{e.suscfirst}. 

It is an important fact that at critical temperature the susceptibility of the subgraph measures tends to infinity. This is also the case for supercritical Ising measures ($\beta > \beta_c$).  For high temperature models however  the susceptibility tends to a finite value. 

We shall also use the fact that the Ising model can be  described via a coupling with an edge percolation model. This is referred to as the FK--Ising model (FK stands for Fortuin and Kasteleyn), or more generally, the random cluster model \cite{GrimmFK}.  For a finite graph $G(V,E)$ the random cluster model with parameters $p$ and $q$ are given by the law 
$$
\phi_{p,q}(\omega) = \frac{1}{Z}\prod_{e \in E}{p^{\omega(e)}(1-p)^{1-\omega(e)}q^{k(\omega)}}
$$
where $q \geq 1$ $p \in (0,1)$ and $k(\omega)$ is the number of connectivity clusters of $\omega$.

It turns out that if one assigns an independent fair coin flip to each cluster of the resulting random graph, the law of the spin system will be that of an Ising model with inverse temperature $\beta$, given by $p = 1-e^{-\beta J}$. This connection often proves to be useful, and some concepts of the Ising model gain a new interpretation in the random cluster model. For instance, the susceptibility of the Ising model is just the expected cluster size of a uniformly chosen vertex.

The random cluster model can also be defined for infinite graphs via weak limits, and the above correspondence with the Ising model remains true, as well. The phase transition can also be interpreted in terms of the random cluster model: the critical inverse temperature $\beta_c$ is the critical point of the corresponding percolation where an infinite open cluster is created.

\subsection{Acknowledgments} 

We are indebted to Amir Dembo for his elegant proof of Proposition~\ref{CWlargeent}, to Christophe Garban for his proof of Theorem~\ref{critical_is} and for a coupling idea that inspired Lemma~\ref{l.eigenclue}, to Hugo Duminil-Copin for his suggestion to use  \cite{L} in the proof of Theorem~\ref{supercrit_is}~(1), to S\'ebastien Ott for telling us about \cite{DSS}, to Bal\'azs Szegedy for his Conjecture~\ref{c.SzyB} and Proposition~\ref{p.SzyB}, to Sasha Glazman and Ron Peled for info regarding Question~\ref{q.3col}, to \'Ad\'am Tim\'ar for his Definition~\ref{d.Adam}, and to a referee for many useful comments. 

A large part of this work was supported by the ERC Consolidator Grant No.~772466--NOISE, and was part of PG's PhD thesis at the Central European University, Budapest. Thanks to Christophe Garban and Bal\'azs Szegedy for reading the thesis and serving on the doctoral committee. Currently, GP is partially supported by the ERC Synergy Grant No.~810115--DYNASNET and by the Hungarian National Research, Development and Innovation Office, OTKA grant K143468.

\section{Sparse reconstruction in general spin systems: ergodicity, expansion, entropy} \label{basicconc}


As a warm-up, this section will discuss some relatively easy, general questions about sparse reconstruction in spin systems. The setup here is certainly not the most general possible, but it includes all the situations we are interested in and most of those we can think of.

We shall consider a sequence of probability measures on finite graphs, converging in the Benjamini-Schramm sense (specified below) to a measure on an infinite graph. The emphasis is on understanding how and when the properties of the limiting measure determine whether sparse reconstruction is possible in the respective sequence of measures.   

Let $G=(V, E)$ be a transitive, edge-labelled, infinite graph. That is, we assume that the group $\mathsf{Aut}(G)$ of label-preserving automorphisms acts on $V$ transitively. This will be our limiting object. The role of the edge labels is to be able to talk about spin systems on $G$ that are not invariant under the full unlabelled automorphism group, such as an Ising measure on $\Z^2$ with different weights on vertical and horizontal edges. In our setup this may be expressed by coloring vertical edges red, horizontal edges blue, respectively. 

We  consider a sequence $\{G_n(V_n, E_n) : n \in \N\}$ of edge-labelled graphs, which converges to $G$ in the \empha{Benjamini--Schramm sense}.  This means that, for every $R\in\N$ and every finite, rooted and edge-labelled  graph $B$, the probability that the $R$-ball around a uniformly randomly chosen vertex from $G_n$ is isomorphic to $B$  converges as $n \rightarrow \infty$. See \cite{BeSc} or \cite[Section 14]{PGG} for details. In general, Benjamini--Schramm convergence allows that the limit be a consistent distribution on the possible realisations of the $R$-balls, for all $R \in \N$. If $G$ and $G_n$ are deterministic transitive graphs, convergence  means that for every $R\in\N$  the $R$-neighbourhood in $G_n$ eventually becomes isomorphic to the  $R$-neighbourhood in $G$. But this setup allows for  sequences of non-transitive  or random graphs converging to $G$, as well. An interesting  example is the sequence of random $d$-regular graphs on $n$ vertices, converging to the $d$-regular tree $\mathbb{T}_{d}$.

For each $n$, we endow the configuration space  $ \{-1,1\}^{V_n}$ with a (usually $\mathsf{Aut}(G_n)$-invariant) measure $\p_n$.  Note that any particular  configuration  $\sigma \in \{-1,1\}^{V}$  on the vertex set can be interpreted as a vertex--labelled graph. Thus $G_n$, its vertex configuration endowed with the measure $\p_n$, can be identified with a vertex- and edge-labelled random graph, which we shall represent with the pair $( G_n, \p_n )$. Then, the sequence of finite vertex- and edge-labelled  random graphs $\{ (G_n, \p_n) : n \in \N  \}$  is required to converge to a transitive, $\mathsf{Aut}(G)$-invariant,  vertex- and edge-labelled  graph $(G, \p)$ in an annealed Benjamini--Schramm sense, in the joint randomness of the labels and the uniform random root.
%
%

One of the first natural questions is whether sparse reconstruction is an attribute of the limiting measure. That is, if $(G_n,\p_n)$ and $(G_n,\mathbb{Q}_n)$ are sequences with the same limit, then either both sequences admit SR, or both do not. The answer is, in general, negative.

It is possible to construct a sequence $(G_n,\p_n)$ converging to a product measure which admits sparse reconstruction:

\begin{ex}[Product limit, SR]\label{prodlim} 
Let $\p_n$ be the following measure on $ \{-1,1\}^{\Z_n}$. We choose a uniformly random $i \in \Z_n$, we flip a fair coin, and around $i$ in a neighborhood of size $\floor{ n^{\frac{2}{3}}}$ we make every spin in the interval $+1$ or $-1$ according to the coin flip. Outside this interval, the spins are IID coin flips. It is easy to check that this spin system converges to the product measure on $ \{-1,1\}^{\Z}$. At the same time, Majority admits sparse reconstruction. Indeed, one can choose $U$ simply to be the multipliers of $n^{\frac{1}{2}}$. With high probability, we can identify where the long $+$ or $-$ interval is, and again with high probability, whether it is $+$ or $-$ will tell us the result of Majority.
\end{ex}

Now let us investigate the same question under stricter conditions. It may seem intuitive that if the limiting spin system is not ergodic, then there should be sparse reconstruction for the sequence. This is the case for example for the supercritical Ising model on the tori $\Z_n^d$ (see Theorem~\ref{supercrit_is}~(3)). This is, however, not true in general. A simple example for a sequence even without RSR is the following.

\begin{ex}[Non-ergodic limit, no RSR]\label{nonerglim} 
Let $G_n$ be a path of $2n$ vertices, and put IID Bernoulli($\frac{3}{4}$) bits on the first $n$ vertices and IID Bernoulli($\frac{1}{4}$) bits to the remaining vertices. It is obvious that the limiting measure is non-ergodic; on the other hand, each $\p_n$ is a product measure, hence there is no RSR, by a straightforward generalization of \cite[Theorem 2.8]{GaPe}.
\end{ex}

Nevertheless, if $G_n$ is an expander sequence and the limiting measure is non-ergodic, we do always have RSR. For a review on expanders, see \cite{Lub}, but here are the relevant definitions.  

For a finite graph $G$, (one version of) the Cheeger constant is
$$
h(G) = \min \left\{ \frac{|\partial W|}{|W|} : W \subseteq V(G),\; |W| \leq \frac{1}{2} |V(G)| \right\},
$$
where $\partial W = \{ (x,y) \in E(G) :\; x \in W, \; y \notin W \}$. 

\bde[Expander sequence]\label{d.expander}
A sequence $\{ G_n(V_n, E_n) : n \in \N \}$ of bounded degree graphs with $|V_n| \to \infty$ is an $h$-expander sequence for some $h>0$ if 
$
h(G_n) > h
$
for every  $n \in \N$. 
\ede

\bpr \label{SRnonerg}
Let  $\{ (G_n, \p_n) : n \in \N  \}$ be a sequence of spin systems, where $\{ G_n(V_n,E_n) : n \in \N \}$ is an $h$-expander sequence with maximum degree $D$, and let $(G, \p)$ be the Benjamini-Schramm limiting spin system, with $G$ transitive.

If the $\mathsf{Aut}(G)$-invariant measure $\p$ is not $\mathsf{Aut}(G)$-ergodic, then there is (Random) Sparse Reconstruction for $\{ (G_n, \p_n) : n \in \N  \}$.
\epr

\bpf
By the non-ergodicity of  the limiting measure $\p$, there is an $\mathsf{Aut}(G)$-invariant event $A$ with non-trivial measure $\p[{A}]=c \in (0,1)$. As any measurable event can be approximated by an event depending on a finite subset of the coordinates (a cylinder event) with desired accuracy, for any  $\delta > 0$ one can find a cylinder event $A_{\delta}$ satisfying  
$$
\p[A \triangle A_{\delta} ] < \delta,
$$
where $A  \triangle B$ is the symmetric difference of the events $A$ and $B$. Now pick a root $o \in V$ and choose $R \in \N$ large enough so that $A_{\delta}$ is $B_R(o)$-measurable (that is, the coordinates on which $A_{\delta}$ depends are inside the ball $B_R(o)$). Since $G$ is transitive, all $R$-balls are isomorphic on $G$, so for any $x \in V$ we can define the  ${B_R(x)}$-measurable event  $A^{x}_{\delta}$, the respective translate of $A_{\delta}$. By the invariance of $A$, for every $x \in V$ the translated event $A^{x}_{\delta}$ is an equally good approximation of $A$,  and thus we have, for any $x,y \in V$,
$$
\p[A^{x}_{\delta} \triangle A^{y}_{\delta} ] < 2\delta.
$$
 
 The next step is to represent these events on the limiting sequence $\{ (G_n, \p_n) : n \in \N  \}$. If $x_n\in V_n$ is such that the $(R+1)$-ball around $x_n$  in $G_n$ is isomorphic to the $(R+1)$-balls in $G$ (which are all isomorphic to each other by transitivity), we can use this rooted isomorphism to define the event $A^{x_n}_{\delta}$ on the finite space $(G_n, \p_n)$. If $B_{R+1}(x_n)$ is not isomorphic to the  $(R+1)$-ball in $G$, then  $A^{x_n}_{\delta}$ is considered  unsatisfied. The point here is that if $(x_n,y_n) \in E_n$, and both $B_{R+1}(x_n)$ and $B_{R+1}(y_n)$ are isomorphic to the $(R+1)$-ball in $G$, then the joint distribution of $\big(\1_{A^{x_n}_{\delta}},  \1_{A^{y_n}_{\delta} }\big)$ agrees with that of $\big(\1_{A^{x}_{\delta}},  \1_{A^{y}_{\delta} }\big)$. Thus, for any $(x_n, y_n) \in E_n$ we have
\beq \label{eq.prdif}
\begin{aligned}
\p_n[\1_{A^{x_n}_{\delta}} \neq \1_{A^{y_n}_{\delta} } ]
&\leq   \p[A^{x}_{\delta} \triangle A^{y}_{\delta} ] + \1_{B_{R+1}(x_n) \not\cong B_{R+1}(o) \text{ or } B_{R+1}(y_n) \not\cong B_{R+1}(o)}\\
& \leq  2 \delta + \1_{B_{R+2}(x_n) \not\cong B_{R+2}(o)}.    
\end{aligned}
\eeq

Now, the sequence of events we are going to reconstruct is the majority of the events $A^{x_n}_{\delta}$. More precisely, for a fixed  $n$, let
$$
S^+_n := \sum_{v \in G_n} \1_{A^{v}_{\delta}}, \qquad  S^-_n := \sum_{v \in G_n} \1_{(A^{v}_{\delta})^c}, \qquad S_n := S^+_n-S^-_n,
$$
and the majority:
$$
\Maj_n^{A_{\delta}} = 
\left\{
	\begin{array}{ll}
		1 & \mbox{if } S_n >0 \\
		-1, & \mbox{otherwise.} 
	\end{array}
\right.
$$
The sketch of our argument is the following. We try to reconstruct the sequence $\{\Maj_n^{A_{\delta}} : n \in \N \}$. First, we are going to show that with high probability either $A^{v}_{\delta}$ holds for most $v \in V_n$, or  $A^{v}_{\delta}$ does not hold for most $v \in V_n$. Then we use this observation to show that $\{ \Maj_n^{A_{\delta}} : n \in \N  \}$ is a non-degenerate sequence of Boolean functions, and there is sparse reconstruction for this sequence from  the set of coordinates $B_R(v)$ where  $v \in V_n$ is  chosen uniformly at random. Indeed, as  $A^{v}_{\delta}$ is  $B_R(v)$-measurable, we can compute   $A^{v}_{\delta}$ for a random $v$. As $A^{v}_{\delta}$ agrees for most  $v \in V_n$, it is  clear that  $A^{v}_{\delta}$ for a random $v$ is a good guess for $\Maj_n^{A_{\delta}}$. 
 
Now we elaborate this argument. Fix an  $\eps>0$, and  for some fixed large $n$  denote by $B_{\eps}$ the event that both $\{ v \in V_n :  A^{v}_{\delta} \}$ and its complement have at least $ \eps|V_n| $ elements. On the event $B_{\eps}$, by the expansion property, there are at least $|V_n|h\eps \geq 2 |E_n|h\eps/D $ edges between $\{ v \in V_n :  A^{v}_{\delta} \}$ and its complement. (Recall that $D$ is the maximal degree of $G_n$ and thus $D|V_n| \geq 2|E_n| $.) So, if we define the random variable $J_n$ as the number of edges $(u,v)$ in $G_n$ such that $\{ \1_{A^{v}_{\delta} } \neq \1_{{A^{u}_{\delta} }}\}$, then $B_{\eps}$ implies 
$$J_n \geq 2 |E_n|h\eps/D.
$$ 
At the same time, by \eqref{eq.prdif}, for any $\eta>0$, if $n$ is large enough, we have  
$$
\E[J_n] <  2 \delta|E_n| + 2\eta |V_n| \leq (2\delta+4\eta) |E_n|,
$$ 
because of Benjamini--Schramm convergence: for $n$ large enough, the probability that $B_{R+2}(v_n)$ with $v_n \in V_n$ chosen uniformly at random is isomorphic to $B_{R+2}(o)$ in $G$ is at least $1-\eta$, or, in other words, the number of vertices $v_n \in V_n$ such that $ B_{R+2}(v_n) \not\cong B_{R+2}(o)$ is at most $\eta |V_n|$.

Thus Markov's inequality says that 
$$
 \p_n[ B_{\eps}] \leq \p_n\left[J_n \geq \frac{2h |E_n| \eps}{D}\right] \leq \frac{D (\delta+2\eta)}{h \eps}.
$$
So, we may choose  $\delta$ and $\eta$ small enough to have $\p_n[ B_{\eps}] \leq \eps$ for every large enough $n$. 

Let us write $A_\eps^+$ for the event that $|\{v \in V_n :  A^{v}_{\delta} \}|>(1-\eps)|V_n|$ and  $A_\eps^-$ for the event that $|\{v \in V_n :  (A^{v}_{\delta})^c \}|>(1-\eps)|V_n|$. On the event $A_\eps^+$ we have $S_n > (1-2\eps)|V_n|$, and on $A_\eps^-$ we have $S_n < -(1-2\eps)|V_n|$. Thus, for $\eps< 1/2$, we have that $B_\eps^c$ is the disjoint union of  $A_\eps^+$ and $A_\eps^-$, and these events imply $\Maj_n^{A_{\delta}}=+1$ and $\Maj_n^{A_{\delta}}=-1$, respectively.

We are now ready to show that  $\{ \Maj_n^{A_{\delta}} : n \in \N  \}$ is a non-degenerate sequence. Since $\p[A^{x}_{\delta} ] > c-\delta$, we have 
$$\E_n[S^+_n] > (1-\eta) |V_n| (c-\delta) > (c-\delta-\eta) |V_n|.$$
On the other hand, by the law of total expectation for $\E_n[S^+_n]$, we have
\begin{align*}
\E_n[S^+_n] &= \E_n[S^+_n |\; A_\eps^+ ] \,\p_n[A_\eps^+] +  \E_n[S^+_n |\; A_\eps^- ] \, \p_n[A_\eps^-] +  \E_n[S^+_n |\; B_\eps ] \, \p_n[B_\eps]\\
&\leq |V_n| \, \p_n[A_\eps^+] +   \eps |V_n| \, (1-\p_n[A_\eps^+]) + |V_n|\eps.
\end{align*}
Thus, 
$$
\frac{c-\delta-\eta-2\eps}{1-\eps} < \p_n[A_\eps^+] \leq \p_n[\Maj_n^{A_{\delta}}=+1] ,
$$
which is at least $c/2$ if $\delta$, $\eta$, $\eps$ are small enough.

We can apply the same argument to $A^-_{\eps}$ to get that
$$
\frac{1-c}{2} <  \p_n[A_\eps^-]  <   \p_n[\Maj_n^{A_{\delta}} = -1] .$$
Thus the sequence $\{ \Maj_n^{A_{\delta}} : n \in \N  \}$ is indeed non-degenerate. 

Now it is easy to see that there is RSR for  $\{ \Maj_n^{A_{\delta}} : n \in \N  \}$ from the random set $B_R(v)$ where  $v \in V_n$  is chosen uniformly (recall that $A^{v}_{\delta}$ is  $B_R(v)$-measurable). Since $G_n$ is of bounded degree, the revealment of this random set tends to $0$.

Let our guess for $ \Maj_n^{A_{\delta}}$ be  $g_{v}:=2\1_{A^{v}_{\delta}}-1$. It is clear that, conditioned on $B^{c}_{\eps}$, the probability of guessing  $\Maj_n^{A_{\delta}}$ correctly is at least $1-\eps$. Thus
$$
\p_n[\Maj_n^{A_{\delta}} \neq g_{v}] < 2\eps\,.
$$
Since we can take $\eps=\eps_n\to 0$ as $n\to\infty$, and two Boolean functions are highly correlated if{f} they are equal with high probability, we get RSR for $\{ \Maj_n^{A_{\delta}} : n \in \N  \}$.
%
%
\epf

Using a similar argument, it is easy to show that in case $(G, \p)$ is not $\mathsf{Aut}(G)$-ergodic, then there is Random Sparse Reconstruction for $(G, \p)$ from any exhaustion, as well. 
\medskip

We have seen that for some non-ergodic measure on some $(G, \p)$ it is possible that a sequence of spin systems converging to $(G, \p)$ admits SR, while another such sequence does not. In particular, this is the case for sequences of product measures. This raises the question whether there are graphs and measures with the property that no matter how we choose the  approximating sequence, we always have sparse reconstruction. See Question~\ref{q.robust}.

\medskip

Another line of questions concerns whether it is true that, if in some sense $\p_n$ contains less randomness than $\mathbb{Q}_n$, and $\mathbb{Q}_n$ admits SR, then $\p_n$ admits SR as well. Of course, the important point here is how we make the expression ``contains less randomness'' precise. A natural attempt is to express the degree of randomness in a sequence with asymptotic entropy.

We start by recalling the definition  of entropy. The entropy of a discrete random variable $Z$ is $\HH(Z):= - \sum_z \p[Z=z] \log \p[Z=z]$  (with $\log$ denoting $\log_2$).

\bde[Asymptotic entropy]\label{d.asympent}
Let $\{ \p_n \}$  be a sequence of measures. The asymptotic entropy of the sequence is
$$
\mathcal{H}(\{ \p_n \}) : = \lim_{n \rightarrow \infty}{\frac{\HH(\p_n)}{\sum_{v \in V_n} {\HH(\p_n|_{\sigma_v})}} }
$$
if it exists, where $\p_n|_{\sigma_v}$ is the distribution of an individual spin $v$.
\ede

It turns out that $\mathcal{H}(\{\mathbb{Q}_n \}) > \mathcal{H}(\{ \p_n \}) $  and $\{ \mathbb{Q}_n\}$ having SR does not imply that $\{ \p_n \}$ has SR. First, Example~\ref{prodlim} above of a spin system that weakly converges to a product measure and still admits SR testifies that it can happen that the asymptotic entropy is $1$ (as large as it can possibly be) but still there is Sparse Reconstruction.

There also exists a sequence $(G_n, \p_n)$ with $\mathcal{H}(\{ \p_n \})=0$ in such a way that $(G_n, \p_n)$ admits no SR. This is a version of Example~\ref{nonerglim}. Consider the cycle $\Z_n$, cut it into two equal halves at a uniform random place, then place IID $\mathrm{Ber}(\eps_n)$ and $\mathrm{Ber}(1-\eps_n)$ bits on the two halves, respectively, for some sequence $\eps_n \to 0$. Then the marginal $\p_n|_{\sigma_0}$ is $\mathrm{Ber}(1/2)$, with constant entropy, while, using the chain rule for entropy,
$$
\HH(\p_n) = \HH(\p_n \,|\, \mathsf{cut}) + \HH(\mathsf{cut}) = \HH(\mathrm{Ber}(\eps_n)) \, n + \log n = o(n),
$$
where $\mathsf{cut}$ is the uniformly random place where the cycle is cut. At the same time, there is no SR, since transitive functions of the spins do not care where $\mathsf{cut}$ is, and hence, from their point of view, this is the same as Example~\ref{nonerglim}. (Nevertheless, RSR does exist: in which half of the cycle $\mathsf{cut}$ is contained can easily be guessed.)

We do not know if such an example exists for ergodic measures. See Question~\ref{q.zero}
\medskip

Another way of expressing that $\p$ has no more randomness than $\mathbb{Q}$ is to say that $(G,\p)$ is a factor of $(H,\mathbb{Q})$. In fact, according to Ornstein-Weiss theory \cite{OW.isom}, for FIID measures on amenable groups this is equivalent to the entropy question above. Indeed, here is a possible counterexample in this factor setting: it might happen that there is SR for the sequence $\mathbb{Q}_n$, while the sequence $\p_n$, which is a factor of $\mathbb{Q}_n$, does not admit SR. For example, one can take $k_n$ independent copies of critical Ising model on the torus, for some $k_n \rightarrow \infty$. We prove in Theorem~\ref{critical_is} that there is sparse reconstruction for the the two-dimensional critical Ising model, and it is easy to see that if there is sparse reconstruction for a model, then there is also sparse reconstruction for the spin system with arbitrary number of independent copies of it. However, if the factor spin system is the product of the spins of all the $k_n$ copies at a given vertex, it is easy to see that, for $k_n$ large enough (say, doubly exponential in the size of the torus), the law of this spin system is very close to the uniform measure in total variation distance, hence no SR is possible. 

So, it appears that sparse reconstruction is not monotone with respect to the most common ways of measuring information. It would be interesting to find some invariant of sequences of spin systems which behaves well with respect to sparse reconstruction.

\section{Magnetization vs.~majority, full vs.~partial sparse reconstruction in general systems}
\label{General_SR}

\subsection{Reconstruction of the total magnetization}\label{ss.magn}

In this section, we consider the following general setup. We have an infinite sequence of  probability spaces $\{(\Omega_n, \p_n ) \}_{n \in \N}$, and for each  $n \in \N$ we have a system of random variables $X^n:=\{ X^n_v : v \in V_n \}$ on $\Omega_n$, where the cardinalities of the finite sets $V_n$ (the set of coordinates) is approaching infinity with $n$. We require each variable $X^n_v$ to have a finite second moment.

In case we talk about sparse reconstruction, we additionally require that for each $n$,  the joint distribution of  $X^n$ has to be invariant under some group action $\Gamma_n$ acting transitively on $V_n$. In particular, the one-dimensional marginals $X^n_v$ have the same distribution for any  $v \in V_n$ in this case.

For a system of random variables $X:=\{ X_v : v \in V \}$ we define the averaging operator $M$ as
\beq \label{magn}
M[X] : = \frac{1}{|V|}\sum_{u \in V}{X_{ u}}.
\eeq
Here $M$ refers to magnetization. In the statistical physics context the variables $X_{ u}$ are spins hence the term. At the same time we introduce a pair of somewhat different, more general magnetization operators associated with a random subset $\mathcal{U}$ of the coordinates: 

\beq \label{eq.pagn}
P^{\mathcal{U}}[X] : = \sum_{S \subseteq V} {\p[\mathcal{U}=S]\E[X \,|\, \mathcal{F}_{S}]},    
\eeq
and
\beq \label{eq.magn}
M^{\mathcal{U}}[X] : = \sum_{S \subseteq V} {\p[\mathcal{U}=S]\frac{\E[X \,|\, \mathcal{F}_{S}]}{\D(\E[X \,|\, \mathcal{F}_{S}])}},    
\eeq
where $\D(X): = \sqrt{\Var(X)}$ denotes the standard deviation of a random variable. It is worth noting that $P^{\mathcal{U}}$ is the direct  generalization of $M$, not $M^{\mathcal{U}}$ .  Still, in this Section  $M^{\mathcal{U}}$ will play a central role. 

More generally,  if we have a system of random variables $ \{ Y_S : S \subseteq V \}$, where $Y_S$ is $\mathcal{F}_{S}$-measurable, then we can define
\beq \label{eq.pagnY}
P^{\mathcal{U}}[Y_{\{ S \subseteq V \}}] : = \sum_{S \subseteq V} {\p[\mathcal{U}=S]Y_{S}},   
\eeq
and
\beq \label{eq.magnY}
M^{\mathcal{U}}[Y_{\{ S \subseteq V \}}] : = \sum_{S \subseteq V} {\p[\mathcal{U}=S]\frac{Y_{S}}{\D(Y_{S})}}.   
\eeq
Importantly, the variance of $P^{\mathcal{U}}[Y_{\{ S \subseteq V \}}]$ and $M^{\mathcal{U}}[Y_{\{ S \subseteq V \}}]$ has an interesting interpretation:
\beq \label{eq.avgcorr}
\begin{aligned}
\Var\big(M^{\mathcal{U}}[Y_{\{ S \subseteq V \}}]\big) &= \sum_{S,T \subseteq V} {\p[\mathcal{U}=S] \, \p[\mathcal{U}=T] \,\frac{\Cov(Y_{S}, Y_{T})}{\D(Y_{S})\D(Y_{T})}}\\
& =  \E\big[\Corr(Y_{\mathcal{U}_1}, Y_{\mathcal{U}_2} \,|\,\mathcal{U}_1, \mathcal{U}_2)\big],
\end{aligned}
\eeq
where $\mathcal{U}_1$ and $\mathcal{U}_2$ are two independent copies of $\mathcal{U}$. Note, in particular, that this average correlation is always non-negative.

In the same way one can verify that
\beq \label{eq.avgcov}
\Var\big(P^{\mathcal{U}}[Y_{\{ S \subseteq V \}}]\big) =  \E\big[\Cov(Y_{\mathcal{U}_1}, Y_{\mathcal{U}_2} \,|\,\mathcal{U}_1, \mathcal{U}_2)\big].
\eeq


The next lemma is going to imply that in case the average correlation defined above is sufficiently high for a system of random variables, then there is full random sparse reconstruction for the generalized magnetisation.

\bl \label{l.randombound}
Consider the system of random variables $X:=\{ X_v : v \in V \}$ on a probability space $(\Omega,\,\p)$. Let  $\mathcal{U} \subseteq V$ be a random subset of coordinates and suppose we have for all $S \subseteq V$ an $\mathcal{F}_{S}$-measurable random variable $Y_{S}$ with $\Var(Y_{S})=1$.
Let $\mathcal{U}^{\oplus k} : = \bigcup_{i=1}^{k}{\mathcal{U}_i}$, where $\mathcal{U}_1,\mathcal{U}_2 \dots, \mathcal{U}_k $ are independent copies of $\mathcal{U}$. Then 
$$
\E\big[\clue(M^{\mathcal{U}}[Y_{\{ S \subseteq V \}}]\,|\,\mathcal{U}^{\oplus k})\big]  \geq 1 - \frac{1}{1 + k \, \E \big[\Corr(Y_{\mathcal{U}_1}, Y_{\mathcal{U}_2} \,|\; \mathcal{U}_1, \mathcal{U}_2)\big] }.
$$ 
  \el

\bpf
We are going to estimate $\E[\clue(M^{\mathcal{U}}[Y_{\{ S \subseteq V \}}]\,|\,\mathcal{U}^{\oplus k})]$ by the average squared correlation between the magnetization of the system and the sample magnetization  $\sum_{i=1}^{k}{Y_{\mathcal{U}_i}}$. Introducing the shorthand $M^{\mathcal{U}}:=M^{\mathcal{U}}[Y_{\{ S \subseteq V \}}]$, we have 
$$
\E\big[\clue(M^{\mathcal{U}}\,|\,{\mathcal{U}^{\oplus k}})\big] 
=
\E \left[ \Corr^2 \big( M^{\mathcal{U}}, \E[M^{\mathcal{U}} \,|\, \mathcal{F}_{\mathcal{U}^{\oplus k}}]\;\bigm|\;\mathcal{U}^{\oplus k}\big)\right] 
\geq 
\E \left[ \Corr^2\Big( M^{\mathcal{U}}, \sum_{i=1}^{k}{Y_{\mathcal{U}_i}}\;\Bigm|\;\mathcal{U}^{\oplus k}\Big) \right],
$$
using that conditional expectation maximizes the absolute value of correlation among all $\cF_{\mathcal{U}^{\oplus k}}$--measurable random variables.

We give a lower bound on $\E \left[ \Corr^2\big( M^{\mathcal{U}}, \sum_{i=1}^{k}{Y_{\mathcal{U}_i}}\,\bigm|\,\mathcal{U}^{\oplus k}\big)\right]$ using that, by the Cauchy-Schwarz inequality for $\E[ (X/Y) \, Y]$, we have $\E[X^2/Y^2] \geq \E^2[X]/\E[Y^2]$ for any pair of non-negative random variables $X$ and $Y$. Thus,
\beq \label{jensenvot}
\begin{aligned}
\E \left[ \Corr^2\Big( M^{\mathcal{U}}, \sum_{i=1}^{k}{Y_{\mathcal{U}_i}}\;\Bigm|\;\mathcal{U}^{\oplus k}\Big) \right]
&= \frac{1}{\Var(M^{\mathcal{U}}) } 
\E \left[ \frac{\Cov^2\big( M^{\mathcal{U}}, \sum_{i=1}^{k}{Y_{\mathcal{U}_i}}\;\bigm|\; \mathcal{U}^{\oplus k} \big)}{\Var\big(\sum_{i=1}^{k}{Y_{\mathcal{U}_i}}\;\bigm| \; \mathcal{U}^{\oplus k}\big)}\right] \\ 
&\geq 
\frac{\E^2 \left[ \Cov\big( M^{\mathcal{U}}, \sum_{i=1}^{k}{Y_{\mathcal{U}_i}}\;\bigm|\; \mathcal{U}^{\oplus k} \big)\right]}{\Var(M^{\mathcal{U}}) \, \E \left[\Var\big(\sum_{i=1}^{k}{Y_{\mathcal{U}_i}}\;\bigm| \; \mathcal{U}^{\oplus k}\big)\right]}. 
\end{aligned}
\eeq
Observe that
\beq \label{eq.covest}
\begin{aligned}
\E \left[ \Cov\Big( M^{\mathcal{U}}, \sum_{i=1}^{k}{Y_{\mathcal{U}_i}}\;\Bigm|\;\mathcal{U}^{\oplus k}\Big) \right]
& =   
\sum_{S \subseteq V}  \sum_{i=1}^{k} {\p[\mathcal{U} =S] \, \E \left[\Cov\big(Y_{S}, Y_{\mathcal{U}_i} \bigm| \mathcal{U}_i \big) \right]}\\
 &=  \sum_{i=1}^{k} \sum_{S \subseteq V} \p[\mathcal{U} =S] \sum_{T \subseteq V} \p[\mathcal{U} =T] \, \Cov(Y_{S}, Y_{T}) \\
&  =  k\, \E \left[ \Corr\big(Y_{\mathcal{U}_1}, Y_{\mathcal{U}_2} \,\bigm|\, \mathcal{U}_1, \mathcal{U}_2\big) \right],
\end{aligned}
\eeq
using first that the random sets $\mathcal{U}_i$ are independent of everything else, and second that  $\D(Y_{S})=1$ for all $S \subseteq V$. Furthermore,
\beq \label{eq.evarest}
\begin{aligned}
 \E \left[\Var\Big(\sum_{i=1}^{k}{Y_{\mathcal{U}_i}}\;\Bigm| \; \mathcal{U}^{\oplus k}\Big)\right] 
&= 
\sum_{i=1}^{k} {\E \left[\Var(Y_{\mathcal{U}_i} \bigm| \mathcal{U}_i)\right]}  + \sum_{i=1}^{k} \sum_{j\neq i}{\E\left [\Cov\big(Y_{\mathcal{U}_i}, Y_{\mathcal{U}_j} \,\bigm|\, \mathcal{U}_i, \mathcal{U}_j\big)\right]}\\
 &= k + (k-1)k \, \E\left [\Corr\big(Y_{\mathcal{U}_1}, Y_{\mathcal{U}_2} \,\bigm|\, \mathcal{U}_1, \mathcal{U}_2\big)\right]. 
\end{aligned}
\eeq
Now we can give a lower bound for the average clue of $\cU^{\oplus k}$. Substituting back \eqref{eq.avgcorr}, \eqref{eq.covest} and \eqref{eq.evarest}   into \eqref{jensenvot}, we obtain that
\begin{align*}
\E \left[ \Corr^2\Big( M^{\mathcal{U}}, \sum_{i=1}^{k}{Y_{\mathcal{U}_i}}\;\Bigm|\;\mathcal{U}^{\oplus k}\Big) \right]
&\geq
\frac{k^2 \, \E\left [\Corr\big(Y_{\mathcal{U}_1}, Y_{\mathcal{U}_2} \,\bigm|\, \mathcal{U}_1, \mathcal{U}_2\big)\right]}{k + (k-1)k \,  \E\left [\Corr\big(Y_{\mathcal{U}_1}, Y_{\mathcal{U}_2} \,\bigm|\, \mathcal{U}_1, \mathcal{U}_2\big)\right] }\\ 
&\geq  \frac{k \,  \E\left [\Corr\big(Y_{\mathcal{U}_1}, Y_{\mathcal{U}_2} \,\bigm|\, \mathcal{U}_1, \mathcal{U}_2\big)\right]}{1 + k \,  \E\left [\Corr\big(Y_{\mathcal{U}_1}, Y_{\mathcal{U}_2} \,\bigm|\, \mathcal{U}_1, \mathcal{U}_2\big)\right] } \\ 
&=1 - \frac{1}{1 + k\,  \E\left [\Corr\big(Y_{\mathcal{U}_1}, Y_{\mathcal{U}_2} \,\bigm|\, \mathcal{U}_1, \mathcal{U}_2\big)\right]},
\end{align*}
finishing the proof of the lemma.
\epf

We can now use  Lemma~\ref{l.randombound} to establish useful sufficient conditions for RSR and SR.
 
\bpr[High average correlation implies full RSR]\label{nonconstr}
Suppose that for every $n \in \N$ we have a  system of random variables $X^n:=\{ X^n_v : v \in V_n \}$ with $\E[(X_v^n)^2]< \infty$ on a common probability space $(\Omega^{V_n},\,\p_n)$, with $\lim_{n \to \infty} |V_n| = \infty$.
 
Suppose that 
$$
\overline{\Corr}_{u, v \in V_n}(X_u^{n}, X^{n}_{v}): = \frac{1}{|V_n|^2}\sum_{u \in V_n}{\sum_{v \in V_n}{\Corr(X_u^{n}, X^{n}_{v})}} \gg \frac{1}{|V_n|}. 
$$
Then there is  a sequence of integers $k_n \ll |V_n|$  with a corresponding random subset $\mathcal{H}^{k_{n}}$ with revealment $\delta(\mathcal{H}^{k_{n}}) \leq \frac{k_n}{|V_n|} $  such that  for any $\eps>0$ we have
$$
\lim_{n \rightarrow \infty} \p_n \left[ \clue(M[X^n]\,|\,\mathcal{H}^{k_{n}}) > 1-\eps \right] =1, 
$$ 
where  $M[X^n] : = \frac{1}{|V_n|} \sum_{v \in V_n} {\frac{X^n_v}{\D(X^n_v)} }$. That is, there is full random sparse reconstruction for $\{X^n,\,\p_n\}_{n \in \N}$. If the joint distribution of $X^n$ is transitive, there is also full sparse reconstruction.
\epr

\bpf
Without loss of generality we may assume that $\D(X_v) = 1$ for every $v \in V_n$. By our assumptions, one can choose the sequence $k_n$ in such a way that
$$  
|V_n| \gg {k_n} \gg \frac{1}{\overline{\Corr}_{u, v \in V_n}(X_u^{n}, X^{n}_{v})}.
$$
Let us  denote by $\mathcal{V}_n$ the random set which assigns uniform probability to each singleton: $\p[ \mathcal{V}_n = \{v\} ]= \frac{1}{|V_n|}$ for every  $v \in V_n$. Now we can apply Lemma~\ref{l.randombound} with $\mathcal{U}=\mathcal{V}_n$ and $Y_{\{v\}} = X_v$ (for non-singleton $S \subseteq V_n$ we can set $Y_{S} = 0$) to get that 
$$
\E\big[\clue(M[X^n]\,|\,\mathcal{V}_n^{\oplus k_n})\big]  \geq  1 - \frac{1 }{k_n \overline{\Corr}_{u, v \in V_n}(X_u^{n}, X^{n}_{v}) +1}  \longrightarrow 1. 
$$

Clue can be at most 1, hence $\clue(M[X^n]\,|\,\mathcal{V}_n^{\oplus k_n})$ tends to $1$ in probability, as well. 

In case the joint distribution of $X^n$ is transitive, $M[X^n]$ is clearly also transitive (invariant under the same transitive action), hence Lemma~\ref{SReqRSR} applies, giving full SR. 
\epf

We would like to highlight the important special case when the joint distribution of $X$ is invariant under a transitive group action of $\Gamma$ on $V$. If $\gamma \in \Gamma$ and $v \in V$ we denote by $\gamma \cdot v$ the image of $v$ by $\gamma$. This action determines in turn an action of $\Gamma$  on $\Omega^{V}$ by $(\gamma \cdot X)_{v}: = X_{\gamma^{-1} \cdot v}$ for $X \in \Omega^{V}$.

We introduce the \emph{susceptibility} $\Su(X)$ of the random variable $X$ as 
\beq \label{susc}
\Su(X)  :=  \frac{1}{|\Gamma_{v}|} \sum_{\gamma \in \Gamma}{\Cov( X_v, X_{\gamma\cdot v})},
\eeq
 where $v \in V$ arbitrary and $\Gamma_{v}$ is the stabilizer subgroup of $v$. As long as the action of $\Gamma$  on $V$ is transitive, this definition does not depend on the choice of $v \in V$. In case the action of $\Gamma$ on $V$ is not free, that is, the stabilizer subgroup of a vertex is non-trivial, then for every $X$ we count every term of the form $\Cov( X_v, X_u)$ exactly $|\Gamma_{v}|$ many times. Indeed, as $|\Gamma| = |V||\Gamma_{v}|$ we have $|\Gamma_{v}|$ times ``too many'' terms in the susceptibility, hence the factor $1/|\Gamma_{v}|$. Warning: in case $X$ has additional symmetries it is possible that there are still repetitions in the sum of $\Su(X)$. 
 
The concept of susceptibility is borrowed from the Ising model. It can be shown that for the Ising model this quantity measures the change in the magnetic field of the system upon a small change in the external magnetic field, hence the name. In case absolute convergent, the susceptibility can be defined for countably infinite product spaces, as well.  

Also in general, susceptibility and magnetization are related. By $\Gamma$-invariance of the measure, for any $\alpha \in \Gamma$ we have  
$$
\sum_{\gamma \in \Gamma}{\Cov( X_v, X_{\gamma\cdot v})} =   \sum_{ \gamma \in \Gamma}{\Cov( X_{\alpha\cdot v}, X_{\gamma\cdot (\alpha \cdot v)})} =  \sum_{ \delta \in \Gamma} \Cov( X_{\alpha\cdot v}, X_{\delta \cdot v}),
$$
and therefore, using $|\Gamma|=|V| |\Gamma_v|$ and canceling the $|\Gamma_v|$ repetitions for each $v\in V$
\begin{align*}
\Var (M[X] ) &= \frac{1}{|V|^2} \sum_{v,w} \Cov(X_v,X_w) \\
&= \frac{1}{|\Gamma|^2}\sum_{\gamma \in \Gamma}{ \sum _{\alpha \in \Gamma}{ \Cov( X_{\alpha\cdot v}, X_{\gamma\cdot v})}}  = \frac{1}{|\Gamma|}\sum _{\gamma \in \Gamma} {\Cov( X_v, X_{\gamma\cdot v })}.
\end{align*}
Thus, we have the following relationship between (the non-normalized) $M[X]$ and $\Su(X)$:
\beq  \displaystyle \label{trans1}
\Var(M[X] ) = \frac{|\Gamma_{v}|\Su(X)}{|\Gamma|} =  \frac{\Su(X)}{|V|},
\eeq
using again that  $|\Gamma| = |V||\Gamma_{v}|$.

In the sequel, to  simplify notation, we are going to avoid the factor $1/|\Gamma_{v}|$, by assuming that the action of the group on the coordinate set is free. We note that all the results extends trivially for general transitive actions.

\bc \label{nonconstsusc}
Suppose that $\{\sigma^n: n \in \N\}$ is a sequence of  $\{-1,1\}^{V_n}$-valued random variables with  distribution invariant under the transitive free group action of $\Gamma_n$ on $V_n$, and that $\Var(\sigma_v) = s^2$ for every $v \in V_n$.

If $\Su(\sigma^n) \rightarrow \infty$, or equivalently $\Var(M[\sigma^{n}])\gg \frac{1}{|V_n|}$, then there is  a sequence of numbers $k_n=o(|V_n|)$  such that for a uniform random subset $\mathcal{H}^{k_{n}}$ of size $k_n$ and for any $\eps>0$ 
$$
\lim_{n \rightarrow \infty} \p_n\left[ \clue\big(M[\sigma^{n}]\;\big|\;\mathcal{H}^{k_{n}}\big) > 1-\eps\right] =1. 
$$ 
\ec
\bpf
It is straightforward to check that  $\Su(\sigma^n) \rightarrow \infty$ is equivalent to $\overline{\Corr}_{u, v \in V_n}(\sigma_{u}^n, \sigma^n_{ v}) \gg \frac{1}{|V_n|}$ when $\Var(\sigma_j)$ is constant, hence Proposition~\ref{nonconstr} applies.
\epf

\subsection{Reconstruction of majority} \label{ss.maj}

It is important to note that, in general, sparse reconstruction of the total magnetization does not imply sparse reconstruction for the Booleanized version, the sequence $\Maj_n$ of majorities (or any other non-degenerate sequence of Boolean functions), by the following simple example.

\begin{ex}\label{SR_nondSR}
Let us take the convex combination of a uniform IID measure on $\{-1,1\}^n$ and another measure on the same space in which all the spins are $+1$ or all the spin are $-1$ with probability $\frac{1}{2}$, respectively. Namely, let $\sigma_n$ be sampled from the uniform IID measure with probability $1-\frac{1}{\sqrt{n}}$ and from the second measure with probability $\frac{1}{\sqrt{n}}$. Let $M_n =\frac{1}{n} \sum_{j \in [n]} {\sigma_j}$. It is clear that in this mixed system $\Var(M_n) \gg 1/n$, and consequently, by Proposition~\ref{nonconstr}, the magnetization can be reconstructed. On the other hand, since the measure is close to an IID in total variation, there is no sparse reconstruction for Majority, or any other non-degenerate Boolean function.
\end{ex}

This example also shows that SR does not imply non-deg SR. Indeed, for this sequence of measures there is sparse reconstruction, but no non-degenerate sparse reconstruction.

Nevertheless, in this subsection we will give some useful conditions under which non-deg RSR, by reconstructing majority, is in fact possible.
   
We first state a lemma that we are going to use to infer Random Sparse Reconstruction  for the majority, whenever the magnetization has high fluctuations with a positive probability.
\bl \label{l.guess2ndSR}
Let $f: \{-1,1\}^{V}  \longrightarrow \{-1,1\}$ and  $\cU \subseteq V$ a random subset independent of $\p$ of $\{-1,1\}^{V}$. Let  $g_{\mathcal{U}}: \{-1,1\}^\cU \lora \{-1,1\}$ be a random function such that for fixed $U \subseteq V$, $g_{U}$ is $\mathcal{F}_U$-measurable.
Let $\mathcal{F}_{\mathcal{U}}$ denote the smallest $\sigma$-algebra such that $\mathcal{U}$ and all the spins in $\mathcal{U}$ are $\mathcal{F}_{\mathcal{U}}$-measurable and suppose that there is a $\mathcal{F}_{\mathcal{U}}$-measurable event $A$ on $\{-1,1\}^{V}$ with  $\p[A] =p$ such that
$$
 \p[f \neq g_{\mathcal{U}} \,|\,A ]< \eps.
$$
Then, with $|\E[f]| = \mu$,
\begin{equation}\label{eq.2casebound}
 \E\left[\clue(f \,|\,\mathcal{U}) \right]\geq
\begin{cases}
         \frac{p - \mu^2}{1 - \mu^2}-5\eps & \text{if } \mu \leq p\,,\\
         \frac{p}{1-p}\frac{1-\mu}{1+\mu}-5\eps & \text{if }  \mu > p\,.
   \end{cases}
\end{equation}
Moreover, if $\mu \leq p$, then
$$
 \p\left[\clue(f \,|\,\mathcal{U})\geq \frac{p - \mu^2}{1 - \mu^2}-5\sqrt{\eps} \right] > 1-\sqrt{\eps}, 
$$
and a corresponding bound applies if $\mu > p$.
\el

\bpf Recall that $$\E\left[\clue(f \,|\,\mathcal{U}) \right] \geq \frac{\Cov (f, g_{\mathcal{U}})}{\Var(f)},$$ 
hence we want to get a lower bound on the covariance.
 
Observe that $\Var(h) = 1 - \E^2[h]$, for any $h$ mapping to $ \{-1,1\}$. Consequently,
$$
\big|\Var(f \,|\,A) - \Var(g_{\mathcal{U}}  \,|\,A)\big| \leq \big|\E[f\,|\,A] -\E[ g_{\mathcal{U}} \,|\,A]\big| \cdot \big|\E[f\,|\,A] +\E[ g_{\mathcal{U}} \,|\,A]\big| \leq  4\eps,
$$ 
using that $\p[f \neq g_{\mathcal{U}} \,|\,A ]< \eps$. Therefore, we have 
 \beq \label{eq.thencorr1}
  \begin{aligned}
 \Cov (f, g_{\mathcal{U}}\,|\,A) 
 &=\frac{1}{2}\big(\Var(f \,|\,A) + \Var(g_{\mathcal{U}}  \,|\,A) -   \Var(f - g_{\mathcal{U}}  \,|\,A) \big) \\
 &\geq \Var(f \,|\,A) - 4 \eps.
  \end{aligned}
 \eeq 
By the law of total covariance, noting that $\E [f\,|\,\mathcal{U}] = \E [f]$, we have 
$$ 
 \Cov (f, g_{\mathcal{U}}) = \E [\Cov (f, g_{\mathcal{U}}\,|\,\mathcal{U})] + \Cov(\E [f \,|\, \cU ], \E [g_{\mathcal{U}}\,|\,\mathcal{U}]) = \E [\Cov (f, g_{\mathcal{U}}\,|\,\mathcal{U})].
$$
Using the law of total covariance again, this time for the indicator random variable $\1_{A}$:
$$
 \Cov (f, g_{\mathcal{U}}) = p\Cov (f, g_{\mathcal{U}}\,|\,A) +  (1-p)\Cov (f, g_{\mathcal{U}}\,|\,A^{c}) + \Cov (\E [f\,|\,\1_{A}], \E [g_{\mathcal{U}}\,|\,\1_{A}]).
$$
Without loss of generality, we may assume that  $ \Cov (f, g_{\cU}\,|\,A^{c})\geq 0$. Indeed, if this was not the case, we could define a better guess $ \tilde{g}$ by
$$
 \tilde{g}_{\mathcal{U}}(\sigma) :=  g_{\mathcal{U}}(\sigma)
$$ 
whenever $(\sigma,\mathcal{U}) \in A$, and 
$$
 \tilde{g}_{\mathcal{U}}(\sigma)  \sim  2 \mathsf{Ber}(\p [f =1 \,|\,A^{c}]) -1
$$
 if $(\sigma,\mathcal{U}) \notin A$, independently of both $\sigma$ and $\mathcal{U}$. Observe that this is the point in the argument where we use that $A$ is $\mathcal{F}_{\mathcal{U}}$-measurable.  From the above it is clear that we may (and will) also assume that 
\beq\label{e.opti}
\E [g_{\mathcal{U}}\,|\,A^{c}] = \E [f \,|\,A^{c}].
\eeq
Now we go on with an estimate similar to that of \eqref{eq.thencorr1}. First, observe that
\begin{align*}
\Var(\E [g_{\mathcal{U}}\,|\,\1_A]) &= p(1-p)\left(\E [g_{\mathcal{U}}\,|\,A] -\E [g_{\mathcal{U}}\,|\,A^{c}] \right)^2 \\
&\geq p(1-p)\left(|\E [f\,|\,A] -\E [f \,|\,A^{c}]| - 2\eps \right)^2  \\
&\geq \Var(\E [f\,|\,\1_A]) -2\eps,
 \end{align*}
where we first used (\ref{e.opti}) and  $\left|\E [g_{\mathcal{U}}\,|\,A] - E [f\,|\,A] \right| \leq 2\eps$, and then, in getting the term $-2\eps$ after squaring, that $|\E [f\,|\,A] -\E [f \,|\,A^{c}]|\leq 2$ and that $p(1-p)\leq 1/4$.
Thus,
\beq \label{eq.covexpclose}
\begin{aligned}
\Cov (\E [f\,|\,\1_{A}], \E [g_{\mathcal{U}}\,|\,\1_{A}]) 
&\geq  \frac{1}{2}\left(\Var(\E [f\,|\,\1_{A}]) + \Var(\E [g_{\mathcal{U}}\,|\,\1_{A}])  -4\eps^2\right)\\
&\geq \Var(\E [f\,|\,\1_{A}])-\eps-2\eps^2,
 \end{aligned}
\eeq
where we used that $\Var(\E [f\,|\,\1_{A}] - \E [g_{\mathcal{U}}\,|\,\1_{A}] )<4\eps^2$.

We thus get the following estimate for the covariance:
\beq \label{eq.covatleast}
\Cov(f,g_\cU) \geq p\left(\Var(f \,|\,A) - 2\eps\right) + 0 + \Var\left( \E[f  \,|\, \1_{A}]\right) - \eps-2\eps^2.
\eeq
Now let $x:= \E [f\,|\,A]$, and let us assume, without loss of generality, that $\mu = \E[f] \ge 0$. Clearly, $\Var(f \,|\,A) = 1 - x^2$ and $\E [f\,|\,A^{c}] = \frac{\mu - p x}{1-p}$. It follows from elementary calculations that 
$$
 \Var(\E[f  \,|\, \1_{A}]) =  p(1-p)\left(\E [f\,|\,A] -\E [f \,|\,A^{c}]\right)^2 = \frac{p}{1-p} (x-\mu)^2. 
$$
Now we express, for given $p$ and $\mu$, the RHS of~(\ref{eq.covatleast}) without the $\eps$-dependent terms, as a function of $x$:
$$
C_{p, \mu}(x) = p(1-x^2) + \frac{p}{1-p} (x-\mu)^2.
$$
It is clear that $-1\leq x \leq 1$. As the quadratic term of $C_{p, \mu}$ is positive for any $p \in (0,1)$, the expression admits a minimum. Indeed, at the minimum,
$$
C'_{p, \mu}(x) = -2px + \frac{p}{1-p}(2x-2\mu) = \frac{2p}{1-p}(p x-\mu) = 0,
$$ 
which yields $x = \mu/p$. In case $\mu \leq p$, that is, if the global minimum lies in the range $x\in[-1,1]$, then
$$
C_{p, \mu}(x) \geq C_{p, \mu}\left(\frac{\mu}{p}\right) = p - \mu^2.
$$
If, on the other hand, $\mu > p$, then the minimum over $x \in [-1, 1]$ is at $x=1$, with 
$$
C_{p, \mu}(x) \geq C_{p, \mu}(1) =  \frac{p}{1-p} (1-\mu)^2.
$$
After writing back the bounds into \eqref{eq.covatleast} and dividing by the variance, we obtain \eqref{eq.2casebound}.

%
As for the second statement, define the $\mathcal{U}$-measurable event:
$$
B =\big\{ U \, : \,  \p[f \neq g_{U} \,|\,A, \{\cU=U\} ]< \sqrt{\eps}  \big\}.
$$
Note that the condition $\p[f \neq g_{\mathcal{U}} \,|\,A ]< \eps$ implies that
$$
\p[B]\geq 1-\sqrt{\eps}.
$$
Now we can repeat the proof of \eqref{eq.2casebound} conditioned on $B$ for clue instead of expected clue and using $\sqrt{\eps}$ instead of $\eps$.
\epf

\bpr \label{SR_maj} 
Let $\{\sigma^n \}_{n \in \N}$ be a sequence of  spin systems $(\{-1,1\}^{V_n},\p_n)$, as above. Suppose that there are some $c> 0$ and a sequence of positive numbers $a_n$ such that $\frac{a_n}{\sqrt{|V_n|}} \rightarrow \infty $, and for every large $n$ it holds that 
\beq \label{sygmunpaley}
\p_n\left[\Big|\sum_{j \in V_n}{\sigma^n_j}\Big|\geq a_n\right]> c.
\eeq
Furthermore, suppose that $\Maj_n$ is a non-degenerate sequence of functions with respect to  $\{\sigma^n \}_{n \in \N}$. In particular, let $\mu = \limsup_{n} {|\E[\Maj_n]|} < 1$.

Then there is a sequence $p_n \rightarrow 0$ and an $ \gamma(\mu,c)>0$ such that, for every $\eps>0$ it holds for every large enough $n$ that 
$$
 \p\left[\clue(\Maj_n\,|\,\mathcal{B}^{p_n}) \geq \gamma(\mu,c)  \right] \geq 1-\eps,
$$
where $\mathcal{B}^{p_n}$ is the random set in which every element is chosen independently with probability $p_n$, and thus there is non-deg RSR.
Futhermore, for fixed $\mu<1$, we have
$$
\lim_{c\to 1} \gamma(\mu,c) = 1.
$$
\epr 

\bpf
Let $S_n := \sum_{v \in V_n}{\sigma_v^n} $ and  given a  spin configuration $\sigma^n$ consider the observed magnetization $\Sigma_n = \sum_{v \in\mathcal{B}^{p_n}}{\sigma_v^n}$, where $\mathcal{B}^{p_n}$ is a Bernoulli subset with parameter $p_n$.
We are guessing $\Maj_n$ from  $\mathcal{B}^{p_n}$ as
$$
\Maj^{p_n} :=\sign (\Sigma_n). 
$$
Now let $A : = \{|\Sigma_n| \geq p_n\frac{a_n}{2} \} $. We are going to show first that, for some sequence $\{p_n \}  \to 0$,
\beq \label{eq.melso}
\p[A] \geq \p\left[\left|S_n \right|\geq a_n\right] -\eps > c-\eps,
\eeq
and second that
\beq \label{eq.masod}
\p\left[\Maj^{p_n}_n \neq \Maj_n   \;\vert\ A \right] \leq \eps.
\eeq 
As $A$ is clearly $\mathcal{F}_{\mathcal{B}^{p_n}}$-measurable, Lemma  \ref{l.guess2ndSR} can be applied to conclude our statement.

Let $P_n$ and $N_n$ be the number of $1$ and  $-1$ spins in  a Bernoulli sample $\mathcal{B}^{p_n}$, respectively, given the spin configuration $\sigma^n$. It is clear that  $\Sigma_n = P_n - N_n$, and, conditioned on  $S_n$, the random variables $P_n$ and $N_n$  both follow binomial distributions, with parameters $\left(\frac{|V_n|}{2} + \frac{S_n}{2}, p_n \right)$ and $\left(\frac{|V_n|}{2} - \frac{S_n}{2}, p_n \right)$, respectively. Therefore, $\D(P_n \,|\, S_n), \D(N_n \,|\, S_n) \leq \sqrt{|V_n| p_n}$.

For any $\alpha_n >0$, we have
\begin{align*}
\p\left[|\Sigma_n - p_n  S_n |\geq \alpha_n \;\big\vert \; S_n  \right] 
&=\p\left[|P_n - \E[P_n] - (N_n - \E[N_n])|\geq \alpha_n  \;\big\vert \; S_n \right] \\
&\leq \p\left[|P_n - \E[P_n]| \geq \alpha_n/2 \;\big\vert \; S_n \right] +\p \left[|N_n - \E[N_n]| \geq \alpha_n/2  \;\big\vert \; S_n  \right].
\end{align*}
Our interest is in the case when $\alpha_n = \beta a_n p_n$, for an arbitrary fixed $\beta>0$. Using the condition that $p_n \gg |V_n| / a_n^2$, for any $C>0$ one can choose $n$ large enough so that 
$$
\beta a_n p_n > C \sqrt{|V_n| p_n}. 
$$
Thus, by Chebyshev's inequality, for any $\eps>0$ one can choose large enough $C$ such that
\begin{align*}
  \p\left[|P_n - \E[P_n]| \geq \beta a_n p_n \;\big|\; S_n \right] \leq \p\left[|P_n-\E[P_n]| \geq  C\D(P_n \;\big|\; S_n )\;\big| \; S_n \right] 
<\eps/2.
\end{align*}
Obviously, the same argument shows that $\p \left[|N_n - \E[N_n]| \geq \beta a_n p_n \;\big\vert \; S_n \right]<\eps/2$, as well. We got that, for any $\beta>0$ and  $\eps>0$, if $n$ is large enough, then
\beq \label{eq.binombound}
\p\left[|\Sigma_n - p_n  S_n |\geq \beta a_n p_n  \;\big\vert \; S_n \right] < \eps.
\eeq
This observation will allow us to prove both \eqref{eq.melso} and \eqref{eq.masod}.

We start by \eqref{eq.melso}. First note that the conditioning in $\p\left[|\Sigma_n| < p_n\frac{a_n}{2} \;\big\vert \; |S_n| > a_n\right]$ can happen in two different ways: either $S_n > a_n$ or $S_n < -a_n$. The two cases lead to essentially the same estimation, so we only cover the first one. 
\beq \label{eq.nobigdiff}
\begin{aligned}
\p\left[|\Sigma_n| < p_n\frac{a_n}{2} \;\Big| \; S_n > a_n\right] &\leq \p\left[\Sigma_n < p_n\frac{a_n}{2} \;\Big| \; S_n > a_n\right]\\
&\hskip -4 cm \leq\p\left[\Sigma_n < p_n\frac{a_n}{2} \;\Big| \; S_n = a_n\right] \leq\p\left[|\Sigma_n - p_n  S_n |\geq  p_n\frac{a_n}{2}  \;\Big| \; S_n = a_n \right] < \eps, 
\end{aligned}
\eeq 
using \eqref{eq.binombound} with $\beta=1/2$. Together with a similar estimate for $\p\left[|\Sigma_n| < p_n\frac{a_n}{2} \;\big\vert \; S_n < -a_n\right]$, this  proves \eqref{eq.melso}.

Recall that  $A = \{|\Sigma_n| \geq p_n\frac{a_n}{2} \} $. To prove \eqref{eq.masod}, we split the probability according to the value of $\Maj_n$, that is, the sign of $S_n$:
\beq \label{eq.badguessest}
\begin{aligned}
\p\left[\Maj^{p_n}_n \neq \Maj_n   \;\big\vert\ A \right] 
&= \p\left[\Sigma_n <-p_n\frac{a_n}{2} \;\Big\vert\; A\right] \p\left[ S_n \geq 0 \;\Big\vert\; \Sigma_n <-p_n\frac{a_n}{2}\right]\\
&\qquad + \p\left[\Sigma_n > p_n\frac{a_n}{2} \;\Big\vert\; A\right] \p\left[ S_n < 0 \;\Big\vert\; \Sigma_n > p_n\frac{a_n}{2}\right].
\end{aligned}
\eeq 
We may assume that  $ \p[\Sigma_n <-p_n\frac{a_n}{2} \,\vert\, A]$ and $\p[\Sigma_n > p_n\frac{a_n}{2} \,\vert\, A] = 1-\p[\Sigma_n <-p_n\frac{a_n}{2} \,\vert \, A]$ are both bounded away from $0$. Otherwise, it is enough to consider one of the two terms. 

Using that $\p[B|A] = \p[A|B]\frac{\p[B]}{\p[A]}$, we can estimate the first term in  \eqref{eq.badguessest} as
\begin{align*}
\p\left[ S_n \geq 0 \;\Big|\; \Sigma_n <-p_n\frac{a_n}{2}\right]
= \p\left[\Sigma_n <-p_n\frac{a_n}{2} \;\Big|\; S_n \geq 0 \right]  \frac{\p[S_n \geq 0]}{\p[ \Sigma_n <-p_n\frac{a_n}{2}]}< \eps,
\end{align*}
as $\p[ \Sigma_n <-p_n\frac{a_n}{2}]^{-1}$ is bounded, since we assumed the probability to be bounded away from $0$, while $\p\left[\Sigma_n <-p_n\frac{a_n}{2} \,|\,S_n \geq 0 \right]$ converges to 0 because of \eqref{eq.binombound}, just as in \eqref{eq.nobigdiff}. The second term in~\eqref{eq.badguessest} can be bounded in the same way.
\epf

\subsection{From partial to full sparse reconstruction} 
The following statement, although it follows from some elementary facts by a straightforward calculation, has a few important consequences. 
\bl \label{corr_lemma} Let $X_V = \{ X_v  : v \in V \}$ be a system of random variables on a common probability space. Let $\{ Y_{S}  :  S \subseteq V \} $, where $Y_S$ is $X_S$-measurable  and $Z$ be $X_V$-measurable random variables with finite second moment. Let $\cU$ be a random subset of $V$  independent of $X_V$. Then, with the notation of (\ref{eq.magnY}),
$$
\Corr\big(Z, M^{\mathcal{U}}[Y_{\{ S \subseteq V \}}]\big) \, \E\left[\Corr\big(M^{\mathcal{U}}[Y_{S \subseteq V}], Y_{\mathcal{U}}\;\big|\; \mathcal{U}\big)\right] = \E\big[\Corr(Z,Y_{\mathcal{U}} \,|\,\mathcal{U})\big].
$$
\el

\bpf
 On the one hand, we have
$$
\Cov(Z, M^{\mathcal{U}}[Y_{\{ S \subseteq V \}}]) = \sum_{S \subseteq V} {\p[\mathcal{U}=S]\frac{\Cov(Z,  Y_{S})}{\D( Y_{S})}} = \D(Z) \, \E[\Corr(Z,Y_{\mathcal{U}} \,|\,\mathcal{U})],
$$
and therefore
$$
\Corr(Z, M^{\mathcal{U}}[Y_{\{ S \subseteq V \}}]) = \frac{\E[\Corr(Z,Y_{\mathcal{U}} \,|\,\mathcal{U})]}{\D(M^{\mathcal{U}}[Y_{\{ S \subseteq V \}}])}.
$$
On the other hand, using \eqref{eq.avgcorr},
\begin{align*}
 \E\left[\Corr\big(M^{\mathcal{U}}[Y_{S \subseteq V}], Y_{\mathcal{U}}\;\big|\;\mathcal{U}\big)\right] 
& = \frac{1}{\D(M^{\mathcal{U}}[Y_{\{ S \subseteq V \}}])}\sum_{S \subseteq V} {\p[\mathcal{U}=S] \frac{1}{\D( Y_{S})}\E\left[\frac{\Cov( Y_{S}, Y_{\mathcal{U}} \,|\,\mathcal{U})}{\D(Y_{\mathcal{U}} \,|\,\mathcal{U})]} \right]}\\
 &= \frac{1}{\D(M^{\mathcal{U}}[Y_{\{ S \subseteq V \}}])} \sum_{S,T \subseteq V} {\p[\mathcal{U}=S] \,\p[\mathcal{U}=T] \, \frac{\Cov( Y_{S}, Y_{T})}{\D( Y_{S})\D(Y_{T})}} \\
 &=  \frac{\Var(M^{\mathcal{U}}[Y_{\{ S \subseteq V \}}])}{\D(M^{\mathcal{U}}[Y_{\{ S \subseteq V \}}])} = \D(M^{\mathcal{U}}[Y_{\{ S \subseteq V \}}]),
 \end{align*}
which proves the statement.
\epf

Note that putting together the last calculation with \eqref{eq.avgcorr} also gives rise to the following, somewhat bizarre, identity:
\beq \label{strange}
\E\big[\Corr(Y_{\mathcal{U}_1}, Y_{\mathcal{U}_2} \,|\,\mathcal{U}_1, \mathcal{U}_2)\big] =  \E^2\left[\Corr\big(M^{\mathcal{U}}[Y_{S \subseteq V}], Y_{\mathcal{U}}\;\big|\;\mathcal{U}\big)\right] , 
\eeq
where $\mathcal{U}_1$ and $\mathcal{U}_2$ are two independent copies of $\mathcal{U}$. 

Also note that Lemma~\ref{corr_lemma} implies that
\beq\label{e.maxcorr}
 \E\left[\Corr\big(M^{\mathcal{U}}[Y_{S \subseteq V}], Y_{\mathcal{U}}\;\big|\;\mathcal{U}\big)\right]  = \max_Z  \E\big[\Corr(Z,Y_{\mathcal{U}} \,|\,\mathcal{U})\big]\,,
\eeq
where $\max_Z$ is taken over all $X_V$-measurable random variables $Z$. This shows that $M^{\mathcal{U}}[Y_{S \subseteq V}]$ is a very natural construction for finding functions with high clue.

The following simple lemma quantifies that a bounded random variable with large variance needs to be non-degenerate. 

\bl \label{l.nondegest}
Let $X$ be a random variable with $\D(X)\geq\sigma$, $\E[X] =0$ and $|X| \leq C$. Then
\beq\label{e.option1}
\p\left[ X > \frac{\sigma}{\sqrt{2}} \right]\geq \frac{\sigma^2}{4C^2} \quad\text{and}\quad \p[X \leq 0] > \min\left(1/2, \frac{\sigma^3}{8\sqrt{2}C^3} \right),
\eeq
or
\beq\label{e.option2}
\p\left[ X < -\frac{\sigma}{\sqrt{2}} \right] \geq \frac{\sigma^2}{4C^2} \quad\text{and}\quad \p[X \geq 0] > \min\left(1/2, \frac{\sigma^3}{8\sqrt{2}C^3} \right).
\eeq
\el

\bpf 
From the law of total expectation for $X^2$  we obtain the following estimate :
\begin{align*}
\p\left[ |X| \geq \frac{\sigma}{\sqrt{2}} \right] 
& = \frac{\Var(X) - \E[X^2\,|\,X^2\leq  \frac{\sigma^2}{2}]}{\E[X^2\,|\,X^2 > \frac{\sigma^2}{2}] - \E[X^2\,|\,X^2\leq  \frac{\sigma^2}{2}]}\\
& \geq \frac{\sigma^2 -  \frac{\sigma^2}{2} }{C^2} = \frac{\sigma^2}{2C^2}.
\end{align*}
Thus we have $\p\left[ X \geq \frac{\sigma}{\sqrt{2}} \right]\geq \frac{\sigma^2}{4C^2}$ or $\p\left[ X \leq -\frac{\sigma}{\sqrt{2}} \right]\geq \frac{\sigma^2}{4C^2}$. Suppose that the first option is satisfied. Suppose further that $\p\left[X \leq 0 \right] \leq 1/2$. Since $\E[X\,|\,X > 0] \geq \frac{\sigma^2}{4C^2} \frac{\sigma}{\sqrt{2}}$ and $-\E[X \,|\,X \leq 0 ] \leq C$, we have the following estimate from the law of total expectation:
$$ 
\p\left[X \leq 0 \right] = \p\left[ X > 0 \right] \frac{\E[X\,|\,X > 0]}{-\E[X\,|\,X \leq 0]} \geq \frac{\sigma^3}{8\sqrt{2}C^3},
$$
and hence \eqref{e.option1} holds. 

If we assume the second option $\p\left[ X \leq -\frac{\sigma}{\sqrt{2}} \right]\geq \frac{\sigma^2}{4C^2}$, then, by a symmetric argument, we get that \eqref{e.option2} holds.
\epf

We are now ready to state and prove the main result of this subsection.

\bth[From partial to full reconstruction]\label{clueto1}
Suppose that for every $n \in \N$ we have a  system of random variables $X^n:=\{ X^n_v : v \in V_n \}$ with $\E[(X_v^n)^2]< \infty$ on a common probability space $(\Omega^{V_n},\,\p_n)$  with $\lim_{n \to \infty} |V_n| = \infty$. 

In the statement that follows, in the case of SR, we further assume for every $n \in \N$ that the joint distribution of  $X^n$ is $\Gamma_n$-invariant, where $\Gamma_n$ acts transitively on $V_n$. 

If there is  RSR / non-degenerate RSR / SR / non-degenerate SR for $\{ X^n \}_{n \in \N}$, then there is also  full RSR / full non-degenerate RSR / full SR / full non-degenerate SR for  $X^n$.
\eth

\bpf 
We start by showing that RSR implies full RSR.
By assumption, there exist a sequence of random subsets $\mathcal{U}_n \subseteq V_n$ with revealment  $\delta(\mathcal{U}_n) \ll o(1)$, independent of $(\Omega^{V_n},\,\p_n)$, and a sequence of functions of $g_n: \Omega^{V_n} \longrightarrow \R$ satisfying for all $n \in \N$
$$
 \E[\clue_{X^n }(g_n\,|\,\mathcal{U}_n)]= \E\left[\Corr^2\big(Z^n, \E[Z^n \,|\, \mathcal{F}_{\mathcal{U}_n}] \;\big|\;\mathcal{U}_n\big)\right] >  c
$$
for some $c > 0$, where $Z^n := g_n(X^n)$.
Using Lemma \ref{corr_lemma} for $Y_S = \E[Z^n \,|\, \mathcal{F}_{S}]$ and $Z = Z^n$, together with the  fact that $0\leq \Corr(Z^n, Y_S) \leq 1$ for any $S \subseteq V_n$ (because $Y_S$ is a projection of $Z^n$, hence $\Cov(Y_S, Z^n)=\Var(Y_S)\ge 0$), we can write
\beq \label{magcorr} 
c <  \E\left[\Corr\big(Z^n, \E[Z^n \,|\, \mathcal{F}_{\mathcal{U}_n}] \;\big|\;\mathcal{U}_n\big)\right] \leq \E\left[\Corr\big(M^{\mathcal{U}_n}[Z^n], \E[Z^n \,|\, \mathcal{F}_{\mathcal{U}_n}] \;\big|\;\mathcal{U}_n\big)\right].
\eeq
Thus, by \eqref{strange}, we have
\beq \label{highcorr}
c^2 < \E\left[\Corr\big( \E[Z^n \,|\, \mathcal{F}_{\mathcal{U}_n}], \E[Z^n \,|\, \mathcal{F}_{\tilde{\mathcal{U}}_n}] \;\big|\;\mathcal{U}_n, \tilde{\mathcal{U}}_n \big)\right], 
\eeq
where  $\tilde{\mathcal{U}}_n$ is independent of  $\mathcal{U}_n$ with the same distribution.
 
We now apply Lemma \ref{l.randombound} for $Y_{S} = \E[Z^n \,|\, \mathcal{F}_{S}]$ and $\mathcal{U}=\mathcal{U}_n$ with a suitable integer $k_n$:
\beq \label{eq.avgcorr2clue}
\E\left[\clue\big(M^{\mathcal{U}_n}[Y]\;\big|\;\mathcal{U}_n^{\oplus k_n}\big)\right]  \geq 1 - \frac{1}{1 + k_n \, \E \left[\Corr\big(Y_{\mathcal{U}_n}, Y_{\tilde{\mathcal{U}}_n} \;\big|\; \mathcal{U}_n, \tilde{\mathcal{U}}_n\big)\right] } \geq  1 - \frac{1 }{k_n c^2 +1} .
\eeq
where we used \eqref{highcorr} for the second inequality. Recall here that $\mathcal{U}_n^{\oplus k_n}$ denotes  the union of  $k_n$ independent copies of $\mathcal{U}_n$.

Now choose the sequence of integers $k_n$ in such a way  that $k_n \to \infty$, but $\delta(\mathcal{U}_n^{\oplus k_n}) \leq \delta(\mathcal{U}_n) k_n \ll  1$. This is possible because $\delta(\mathcal{U}_n) \ll  1$. This choice of $k_n$ ensures that the sequence of random subsets $\{\mathcal{U}_n^{\oplus k_n} \}_{n \in \N}$ is sparse, while it is  also clear that, whenever $k_n \to \infty$, from~\eqref{eq.avgcorr2clue} we have
$$
\E\left[\clue\big(M^{\mathcal{U}_n}[Y]\;\big|\;\mathcal{U}_n^{\oplus k_n}\big)\right] \to 1.
$$ 
Thus there is full RSR for $M^{\mathcal{U}_n}[Y]$.
\medskip

We now turn to the Boolean case. Suppose that $\big\{g_n : \Omega^{V_n} \longrightarrow \{-1,1 \}\big\}_{n \in \N}$ is a non-degenerate sequence of Boolean functions with non-vanishing average clue with respect to the sequence of random sets $\{\mathcal{S}_n \subseteq V_n \}_{n \in \N}$.  Let $\sigma^n := g_n(X^n)$. Using the above argument we can find a sparse sequence of random sets $\{\mathcal{U}_n \subseteq V_n : n \in \N\}$ such that, for any $\eps>0$, for every large enough $n \in \N$,
\beq\label{e.nagyklu1}
\E\left[\clue\big(M^{\mathcal{U}_n}[\sigma^n]\;\big|\;\mathcal{U}_n\big)\right] > 1 -\eps.
\eeq

Before going further, let us collect a few basic properties of this $M^{\mathcal{U}_n}[\sigma^n]$, a measurable function of the configuration $X^n$, which will be useful in the sequel. First, we can repeat \eqref{magcorr} for $Z=M^{\mathcal{U}_n}[\sigma^n]$ to obtain from the high clue~\eqref{e.nagyklu1} that
\beq\label{e.nagyklu2}
\E\left[\Corr\big(M^{\mathcal{U}_n}[\sigma^n], \E[M^{\mathcal{U}_n}[\sigma^n] \,|\,\mathcal{U}_n] \;\big|\;\mathcal{U}_n\big)\right] > 1-\eta\,.
\eeq
Then, it is non-degenerate in the sense that it is uniformly bounded and has a uniformly positive variance. Indeed, by the definition~\eqref{eq.magn} and Jensen's inequality,
\beq\label{e.maxMU}
\begin{aligned}
 \max_{X^n} \big|M^{\mathcal{U}_n}[\sigma^n]\big| 
& \leq \E\left[\D^{-1}\big(\E[\sigma^n\,|\,\mathcal{F}_{\mathcal{U}_n}]\;\big|\;\mathcal{U}_n \big)\right] \\
& \leq \frac{1}{\E\left[\D\big(\E[\sigma^n\,|\,\mathcal{F}_{\mathcal{U}_n}]\;\big|\;\mathcal{U}_n \big)\right]}
= \frac{1}{\D(\sigma^n) \, \E^{\frac{1}{2}}\left[\clue(\sigma^n\,|\, \mathcal{U}_n)\right]}\,,
\end{aligned}
\eeq
showing the uniform boundedness: by assumption, $\sigma^n$ is non-degenerate, and $\E\left[\clue(\sigma^n\,|\, \mathcal{U}_n)\right]$ is bounded away from $0$, as there is RSR for $\sigma^n$ from $\cU_n$.  Regarding the variance, from~\eqref{eq.avgcorr},~\eqref{strange}, and~\eqref{e.maxcorr}, we have
\beq\label{e.VarMU}
\begin{aligned}
\Var[M^{\cU_n}[\sigma^n]] &= \E^2\left[\Corr\big(M^{\cU_n}[\sigma^n], \E[\sigma^n \,|\, \cF_{\cU_n}
] \;\big|\;\mathcal{U}_n\big)\right] \\
&\geq \E^2\left[\Corr\big(\sigma^n,  \E[\sigma^n \,|\, \cF_{\cU_n} ] \;\big|\;\mathcal{U}_n)\right] > c > 0\,,
\end{aligned}
\eeq
which is uniformly positive because $\clue(\sigma^n \,|\, \cU_n) \geq \clue(\sigma^n \,|\, \mathcal{S}_n)$ is  uniformly positive, already.
 
Now, for any $S \subseteq V_n$, define  
\beq\label{e.sigmaS}
\xi^n_S: = 
\left\{
	\begin{array}{ll}
		1 & \mbox{if }\E[M^{\mathcal{U}_n}[\sigma^n]\,|\,\mathcal{F}_{S}] > \alpha \\
		-1 & \mbox{otherwise};
	\end{array}
\right.
\eeq
we shall discuss in the sequel how to choose the value of $ \alpha$. Then, we introduce the shorthand notation:
$$
P^n:=P^{\mathcal{U}_n}[\xi^n_{\{S\subseteq V_n\}}]- \E[P^{\mathcal{U}_n}[\xi^n_{\{S\subseteq V_n\}}]].
$$
This is simply a centered version of the $P^{\mathcal{U}}$ operator as introduced in  \eqref{eq.pagnY} applied to the collection of variables $\xi^n_S$.
 Note that $\{\xi^n_S\}_{S \subseteq V_n}$, hence $P^n$ are $X^n$-measurable.

It is clear that $|P^n|\leq 2$ almost surely. It will be important below that, on the other hand, it has a uniformly positive variance.

\bl[$P^n$ is non-degenerate]\label{l.Pnondeg}
If the constant $\alpha>0$ in~\eqref{e.sigmaS} is chosen appropriately, then there is some $s>0$ such that $\Var(P^n) >s$  for all $n \in \N$.
\el

\begin{proof}
In order to give a lower bound on 
\beq\label{e.VarP}
\Var(P^n) = \sum_{S,T \subseteq V_n} \p[\cU_n=S]\, \p[\cU_n=T] \, \Cov(\xi^n_S,\xi^n_T)\,,
\eeq
our strategy is the following. From~\eqref{e.nagyklu2} we can see that, for typical subsets $S$ and $T$ sampled independently from $\cU_n$, the projections $\E[M^{\mathcal{U}_n}[\sigma^n]\,|\,\mathcal{F}_{S}]$ and $\E[M^{\mathcal{U}_n}[\sigma^n]\,|\,\mathcal{F}_{T}]$ are highly correlated with $M^{\mathcal{U}_n}[\sigma^n]$, hence they are also close to each other in the $L^2$ sense. We will use this fact to show that the level sets $\{ \E[M^{\mathcal{U}_n}[\sigma^n]\,|\,\mathcal{F}_{S}] > \alpha  \}$ defining the spins  $\xi^n_S$, for a suitable $\alpha$, coincide with probability close to $1$. Using that $\xi^n_S$ is non-degenerate, this implies that $\Cov(\xi^n_S, \;\xi^n_T)$ is bounded away from $0$. The atypical pairs $S,T$ contribute little to~\eqref{e.VarP}, completing the proof.

We will now work out the details. In light of~\eqref{e.nagyklu2}, suppose that $S\subseteq V_n$ satisfies
\beq \label{eq.corrcond}
\Corr\big(M^{\mathcal{U}_n}[\sigma^n], \E[M^{\mathcal{U}_n}[\sigma^n]\,|\,\mathcal{F}_{S}]\big) >1-\delta\,.
\eeq 
Then
$$
\frac{\D(\E[M^{\mathcal{U}_n}[\sigma^n]\,|\,\mathcal{F}_{S}])}{\D(M^{\mathcal{U}_n}[\sigma^n])} = \sqrt{\clue(M^{\mathcal{U}_n}[\sigma^n]\,|\, S)} = \Corr(M^{\mathcal{U}_n}[\sigma^n], \E[M^{\mathcal{U}_n}[\sigma^n]\,|\,\mathcal{F}_{S}])>1-\delta\,.
$$
Now note that, for any random variables $X, Y$ with finite variance on the same probability space,
$$
\Var(X -Y ) = \D(X)\D(Y) \left(\frac{\D(X)}{\D(Y)} + \frac{\D(Y)}{\D(X)} -2 \Corr(X, Y) \right).
$$
Applying this  with $X=M^{\mathcal{U}_n}[\sigma^n]$ and $Y=\E[M^{\mathcal{U}_n}[\sigma^n]\,|\,\mathcal{F}_{S}]$, using that $\D(M^{\mathcal{U}_n}[\sigma^n])\leq 1$, and assuming that $\delta$ is small enough:
\beq \label{eq.magcs1}
\Var(M^{\mathcal{U}_n}[\sigma^n] -\E[M^{\mathcal{U}_n}[\sigma^n]\,|\,\mathcal{F}_{S}] ) \leq 1 + \frac{1}{1 -\delta} -2(1-\delta)  \leq 4 \delta\,.
\eeq
This, by the triangle inequality, implies that
\beq \label{eq.projclose}
\Var\left(\E[M^{\mathcal{U}_n}[\sigma^n]\,|\,\mathcal{F}_{S}] -\E[M^{\mathcal{U}_n}[\sigma^n]\,|\,\mathcal{F}_{T}] \right) \leq 8\delta\,,
\eeq
whenever  \eqref{eq.corrcond} holds for both $S$ and $T$ (in place of $S$).

Let us choose a small enough interval $(a,b) \subseteq \big( \E[\sigma^n], \E[\sigma^n]+\D(M^{\mathcal{U}_n}[\sigma^n])/\sqrt{2}\big)$ such that
\beq \label{eq.pigeon}
\p[M^{\mathcal{U}_n}[\sigma^n] \in (a,b)] < \eps.
\eeq
Choose  $\delta$  in such a way that $4 \sqrt[4]{\delta} < b-a$ and  $4 \sqrt{\delta}<\eps$. Let $c := a + \sqrt[4]{\delta}$ and $d := b - \sqrt[4]{\delta}$. This will be our constant $c$ in the definition~(\ref{e.sigmaS}) of $\xi^n_S$.

Observe that, under condition \eqref{eq.corrcond}, we have 
\beq\label{e.Scd}
\p\big[\E[M^{\mathcal{U}_n}[\sigma^n]\,|\,\mathcal{F}_{S}] \in (c,d)\big] < \eps + 4\sqrt{\delta}<2 \eps\,.
\eeq
Indeed, if this was not the case, then, together with \eqref{eq.pigeon}, we would have that $M^{\mathcal{U}_n}[\sigma^n]-\E[M^{\mathcal{U}_n}[\sigma^n]\,|\,\mathcal{F}_{S}] \geq \sqrt[4]{\delta}$ holds with probability at least $4\sqrt{\delta}$, which would imply $\Var(M^{\mathcal{U}_n}[\sigma^n] -\E[M^{\mathcal{U}_n}[\sigma^n]\,|\,\mathcal{F}_{S}] )> 4 \sqrt{\delta} (\sqrt[4]{\delta})^2=4\delta$, contradicting~\eqref{eq.magcs1}.

Assuming still that \eqref{eq.corrcond} is satisfied for both $S$ and $T$, we have
\beq \label{eq.}
\begin{aligned}
\p[\xi^n_S=1, \;\;\xi^n_T=-1] 
& \leq  \p\left[\E[M^{\mathcal{U}_n}[\sigma^n]\,|\,\mathcal{F}_{S}]>d,\;\; \E[M^{\mathcal{U}_n}[\sigma^n]\,|\,\mathcal{F}_{T}] \leq \alpha \right]\\
&\hskip 0.5 cm + \p\left[ \alpha <\E[M^{\mathcal{U}_n}[\sigma^n]\,|\,\mathcal{F}_{S}]\leq d \right] \\
& < 2 \sqrt{\delta} +  2\eps \leq 3 \eps\,,   
\end{aligned}
\eeq  
using that $d-\alpha > 2 \sqrt[4]{\delta}$ and thus the first term can be at most $ 2 \sqrt{\delta}$ in light of \eqref{eq.projclose}, while using~\eqref{e.Scd} for the second term. Thus, by symmetry, $\p[\xi^n_S\neq \xi^n_T]<6\eps$. Now in order to conclude from this that $\Cov(\xi^n_S, \;\xi^n_T)$ cannot be very small, we need the fact that  $\xi^n_S$ and $\xi^n_T$ are non-degenerate. For this, we are going to use Lemma \ref{l.nondegest}.
 
First, obviously,
$$
\big|\E[M^{\mathcal{U}_n}[\sigma^n]\,|\,\mathcal{F}_{S}]\big| \leq \max_{X^n} \big|M^{\mathcal{U}_n}[\sigma^n]\big|, 
$$
hence~\eqref{e.maxMU} shows that $\E[M^{\mathcal{U}_n}[\sigma^n] \,|\, \mathcal{F}_{S}]$ is uniformly bounded. Next, 
$$
\Var\left(\E[M^{\mathcal{U}_n}[\sigma^n]\,|\,\mathcal{F}_{S}] \right) = \Var(M^{\cU_n}[\sigma^n]) \, \Corr^2\big(M^{\mathcal{U}_n}[\sigma^n], \E[M^{\mathcal{U}_n}[\sigma^n]\,|\,\mathcal{F}_{S}]\big) 
$$
is bounded away from 0, because of~\eqref{eq.corrcond} and~\eqref{e.VarMU}.
 
Thus,  Lemma \ref{l.nondegest} applied to $\E[M^{\mathcal{U}_n}[\sigma^n]\,|\,\mathcal{F}_{S}]$ says that, under \eqref{eq.corrcond}, $\xi^n_S$ is non-degenerate; in particular, $\E[\xi^n_S]$ is bounded away from $\pm 1$. Therefore, we have 
$$
\Cov(\xi^n_S, \;\xi^n_T) = \p[\xi^n_S=\xi^n_T]-\p[\xi^n_S\neq \xi^n_T]-\E[\xi^n_S]\E[\xi^n_T] \geq 1-6\eps-\E[\xi^n_S]\E[\xi^n_T]>s\,,
$$ 
for some uniform constant $s>0$. 

Now, for any $\delta>0$, if $\eta$ in~\eqref{e.nagyklu2} is small enough, then the event
$$
\left \{ \Corr\left(M^{\mathcal{U}_n}[\sigma^n], \E[M^{\mathcal{U}_n}[\sigma^n]\,|\,\mathcal{F}_{\mathcal{U}_n}]\,|\, \mathcal{U}_n \right) >1-\delta \right\}
$$
happens with probability at least $1-\delta$. Thus $\Cov\big(\xi^n_{\cU_n},\xi^n_{\tilde\cU_n}\big) > s$ holds with probability at least $(1-\delta)^2$, where $\tilde\cU_n$ is an independent copy of $\cU_n$. On the other hand, $|\Cov\big(\xi^n_{\cU_n},\xi^n_{\tilde\cU_n}\big)|\leq 2$ always holds. Altogether, \eqref{e.VarP} is uniformly positive, finishing the proof of Lemma~\ref{l.Pnondeg}.
\end{proof}

Continuing the proof of Theorem~\ref{clueto1}: using $P^n$, we are now going to show that some sequence of Boolean functions can be reconstructed from the random set  $\mathcal{U}_n^{\oplus k_n}$, where the sequence $\{ k_n \}_{n \in \N}$ will be appropriately chosen. We will also choose a sequence $\{ d_n \}_{n \in \N}$, which together satisfy
\begin{equation}\label{eq.gap}
1 \gg  \frac{d_n}{|V_n|} \gg \frac{1}{\sqrt{k_n}}\,. 
\end{equation}
Let us fix an $n \in \N$, and define the events
$$
A^n_j : =  \left\{\frac{(2j-1) \, d_n}{|V_n|}   \leq P^n\leq \frac{(2j+1) \, d_n}{|V_n|} \right\}
$$  
and
$$
B^n_j : = \left\{ \frac{2 |j| d_n}{|V_n|}\leq  |P^n| \right\},
$$ 
for 
\beq\label{e.j}
j = -\left\lceil{\frac{q|V_n|}{2d_n}}\right\rceil , \dots,-1,0, 1, \dots, \left\lceil{\frac{q|V_n|}{2d_n}}\right\rceil,
\eeq
with $q=\D( P^n)/\sqrt{2}$. This will be a reasonable choice because for $j\ge 0$ we have $\{ P^n \geq q \} \subseteq B^n_j$ and $\{  P^n < 0 \} \subseteq (B^n_j)^{c}$, while for $j\leq 0$ we have  $\{ P^n \leq -q \} \subseteq B^n_j$ and $\{  P^n > 0 \} \subseteq (B^n_j)^{c}$, and therefore $|P^n|<2$ and Lemma~\ref{l.Pnondeg} and Lemma~\ref{l.nondegest} together tell us that, for each $n$ we can choose either any $j \ge 0$ or any $j \leq 0$ from the set~(\ref{e.j}) to ensure that the sequence of Boolean functions  $\{ f_n := \1_{B^n_j} \; \}_{n}$ is non-degenerate. To simplify notation in what follows, let us assume that $n$ is such that any $0 \leq j \leq \frac{q|V_n|}{2d_n}$ works. The other case is very similar.
 
In the spirit of the proof of Proposition \ref{SR_maj}, we now show that conditioning on the complement of $A^n_j $, the event $B^n_j$ can be reconstructed from  the random set  $\mathcal{U}_n^{\oplus k_n}$. For each  $i=1,2, \dots , k_n$ we have an independent sample $\mathcal{U}^{i}_n$ of $\mathcal{U}_n$, and a corresponding random variable $\zeta_i:=\xi^n_{\mathcal{U}^i_n}$. 

Let $p=p(X^n) := \p[ \xi^n_{\mathcal{U}_n} = 1\;| \;X^n ]$. Then $\sum_{i=1}^{k_n}{\frac{1+\zeta_i}{2}}$ is binomially distributed with parameters $(k_n,p)$, and our natural guess for $P^n$ is $ Q^n:=\frac{1}{k_n}\sum_{i=1}^{k_n}{\zeta_i}-\E[\xi^n_{\cU_n}]$. Clearly, 
$$
\E[Q^n\;| \;X^n]=P^n \quad\text{and}\quad \D(Q^n\;| \;X^n)= 2\sqrt{\frac{p(1-p)}{k_n}} \leq \frac{1}{\sqrt{k_n}}.
$$
Now, this implies that
$$
\p\left[ \frac{2j d_n}{|V_n|} \leq Q^n \;\middle| \; \frac{(2j-1) d_n}{|V_n|} > P^n \right]  \leq \p\left[ \big| Q^n - \E[Q^n ] \big| \geq \frac{d_n}{|V_n|}  \;\middle| \;X^n \right].
$$
We can establish the same bound for $\p\left[ \frac{2j d_n}{|V_n|} > Q^n \;\middle| \; \frac{(2j+1) d_n}{|V_n|} < P^n \right]$.  We can conclude that, in light of \eqref{eq.gap}, by Chebyshev's inequality for any $n$ large enough, we have
$$
\p\left[\1_{B^n_j} \neq \1_{\left\{  \frac{2j d_n}{|V_n|}\leq  Q^n \right\}} \; \middle| \; (A^n_j)^c \right] < \frac{\eps}{2},
$$
where $\eps>0$ arbitrary.

Since $d_n \ll |V_n| $ and the events $A^n_j$ are mutually exclusive, for large enough $n$ there exists some $j$ such that $\p[A^n_j] < \frac{\eps}{2}$. For every such $j$:
$$
\p\left[\1_{B^n_j} \neq \1_{\left\{  \frac{2j d_n}{|V_n|}\leq  Q^n \right\}} \right] 
\leq
\p\left[\1_{B^n_j} \neq \1_{\left\{  \frac{2j d_n}{|V_n|}\leq  Q^n \right\}} \; \middle| \; (A^n_j)^c \right]  + \p[A^n_j] < \eps.
$$
Thus we conclude that for arbitrary $\eps >0$  one can find large enough $n$ such that for some $j$ the Boolean function $\1_{B^n_j} $ can be guessed correctly with probability $1-\eps$ from the random subset $\mathcal{U}_n^{\oplus k_n}$ with revealement  $ \delta(\mathcal{U}_n^{\oplus k_n}) \leq k_n \delta(\mathcal{U}_n)$. Since, for Boolean functions, high correlation and high probability of agreement are equivalent, we obtain full non-deg RSR, assuming that $\mathcal{U}_n^{\oplus k_n}$ is indeed sparse.

We are left to show that the sequences $\{ k_n \}_{n \in \N}$ and $\{ d_n \}_{n \in \N}$ can be realized to satisfy all the necessary conditions. Indeed, let $k_n ={\delta}^{-\frac{3}{4}}(\mathcal{U}_n)$ and $d_n = |V_n|\delta^{\frac{1}{4}}(\mathcal{U}_n)$. One can easily verify that, in this case,
$$
1 \gg \frac{d_n}{|V_n|} \gg \frac{1}{\sqrt{k_n}} \qquad \text{and} \qquad k_n \delta(\mathcal{U}_n)\ll 1
$$
are all satisfied. The last inequality ensures that $\mathcal{U}_n^{\oplus k_n}$ is sparse, finishing the proof for RSR.

\medskip

Finally, suppose that our sequence admits (non-deg) SR; that is, there is a sequence of transitive (non-deg Boolean) functions $\{g_n : \Omega^{V_n} \longrightarrow \{-1,1 \}\; \}_{n \in \N}$  and a corresponding sequence of  subsets ${U}_n \subseteq V_n$ with   $|U_n| = o(V_n)$.  In this case, the functions $M^{\cU_n}[Y]$ and $P^n$ and $\1_{B^n_j}$ that we have constructed above are also transitive, and therefore Lemma~\ref{SReqRSR} translates the full (non-deg) RSR shown above to full (non-deg) SR.
\epf

Let us spell out a quantitative version of Theorem~\ref{clueto1}.

\bc \label{c.quantnoSR}
Suppose for each $n \in \N$ we have a  system of random variables $X^n:=\{ X^n_v : v \in V_n \}$ with $\E[(X_v^n)^2]< \infty$ on a common probability space $(\Omega^{V_n},\,\p_n)$  with $\lim_{n \to \infty} |V_n| = \infty$. 

In the case of SR, for every $n \in \N$ we further assume that the joint distribution of  $X^n$ is $\Gamma_n$-invariant, where $\Gamma_n$ acts transitively on $V_n$ .

If there is  no RSR for $\{ X^n \}_{n \in \N}$, then there is a $C>0$ such that for every sequence $\{ f_n \}$ of functions, 
$$
\E[\clue(f_n\,|\,\mathcal{U}_n)]\leq C \sqrt{\delta(\mathcal{U}_n)}.
$$

If there is  no SR for $\{ X^n \}_{n \in \N}$, then there is a $C>0$ such that for every sequence $\{ f_n \}$ of transitive functions, 
$$
\clue(f_n\,|\,U_n)\leq C \frac{\left|U_n\right|}{|V_n|}.
$$
\ec

\bpf
Aiming at a contradiction, suppose that  for any $K > 0$ there exists an index $n=n_K \in \N$ with an $f_n: \{-1,1\}^{V_{n}} \longrightarrow \R$ and an $\mathcal{U}_{n} \subseteq V_{n}$ such that
\beq \label{eq.clue2macs0}
\E[\clue_{X^{n} }(f_n\,|\,\mathcal{U}_{n})] > \sqrt{K \delta(\mathcal{U}_{n})}.
\eeq
In the proof of Theorem \ref{clueto1} we showed (see \eqref{eq.avgcorr2clue}) that, if for a spin system $X$ there exists a  function $f: \{-1,1\}^{V} \longrightarrow \R$  and a random set $\mathcal{U} \subseteq V$ with $ \E[\clue_{ X}(f\,|\,U)] > c$, then for every $\ell \in \N$ there exists a  $g: \{-1,1\}^{V} \longrightarrow \R$ and a random set $\Tilde{\mathcal{U}}$ satisfying  $\delta(\Tilde{\mathcal{U}})\leq \ell \delta(\mathcal{U})$ such that 
\beq \label{eq.clue2macs}
\E[\clue_{X_{V}}(g\,|\,\Tilde{\mathcal{U}})] >  1 - \frac{1 }{\ell c^2 +1}.
\eeq
Let $\ell=\ell_{K} : = \lceil\frac{1}{K \delta(\mathcal{U}_{n})}\rceil$, with $n=n_K$, and consider \eqref{eq.clue2macs} with $f = f_n$ and $\mathcal{U} = \mathcal{U}_{n}$ from~\eqref{eq.clue2macs0}, and $\ell =\ell_n$.  First note that the resulting sequence of subsets $\Tilde{\mathcal{U}}_{n}\subseteq V_{n}$ is sparse by the choice of $\ell_{K}$:
$$
\delta(\Tilde{\mathcal{U}}_{n}) \leq \ell_K  \delta(\mathcal{U}_{n}) \ll 1, \quad \text{as }K\to\infty. 
$$
Second, by \eqref{eq.clue2macs0} and \eqref{eq.clue2macs} we have 
$$
\E[\clue_{X^{n} }(g_n\,|\,\Tilde{\mathcal{U}}_{n})] > 1 - \frac{1 }{\ell_K K\delta(\mathcal{U}_{n}) +1}\geq \frac{1}{2},
$$
 which contradicts the assumption that there is no RSR in the system.

For the corresponding statement for SR, we use the same  argument with  transitive functions $f_n: \{-1,1\}^{V_{n}} \longrightarrow \R$ and corresponding subsets ${U}_{n}$. Then  $\mathcal{U}_{n}$, as usual, can be defined  as a uniformly chosen $\Gamma_n$-translate of ${U}_{n}$. Recall that in such a way we have $\delta(\mathcal{U}_{n}) = \frac{\left|U_n\right|}{|V_n|}$. 

Again, we start with the indirect assumption that  for any $K>0$ there exist large enough index $n=n_K \in \N$  with a transitive $f_n: \{-1,1\}^{V_{n}} \longrightarrow \R$ and an $U_{n} \subseteq V_{n}$ such that
\beq \label{eq.clue2macsSR}
\clue_{X^{n} }(f_n\,|\,U_{n}) > K \frac{\left|U_{n}\right|}{|V_{n}|}.
\eeq
Now, in contrast to the general case,  $\clue_{X^n }(g_n \,|\, \mathcal{U}_n)$ is constant (does not depend on the random translation), and therefore, using the notation $f_{n}(X^{n}) = Z^{n}$,   
$$
 \E[ \clue_{X^{n} }(f_{n}\,|\,\mathcal{U}_{n})]= \Corr^2(Z^{n}, \E[Z^{n} \,|\, \mathcal{F}_{U_{n}}] ) = \E^2[\Corr(Z^n, \E[Z^n \,|\, \mathcal{F}_{\mathcal{U}_{n}}] \,|\,\mathcal{U}_{n})] >  c\,.
$$
Thus, instead of \eqref{highcorr}, the identity~\eqref{strange} now gives the stronger bound
$$
c< \E[\Corr(\E[Z^n \,|\, \mathcal{F}_{\mathcal{U}_n}], \E[Z^n \,|\, \mathcal{F}_{\tilde{\mathcal{U}}_n}] \,|\,\mathcal{U}_n, \tilde{\mathcal{U}}_n)], 
$$
which, via~\eqref{eq.avgcorr2clue}, results in an improved version of \eqref{eq.clue2macs}:
$$
\clue_{X_{V}}(g\,|\,\Tilde{U}) >  1 - \frac{1 }{\ell c +1}.
$$
This allows us to get a contradiction with \eqref{eq.clue2macsSR} the same way as in the general case.
\epf

It is worth noting the different bounds we obtained for SR and RSR. We know that, for product measures, we have  the upper bound $\delta(\mathcal{U}_{n})$ instead of $\sqrt{\delta(\mathcal{U}_n)}$, so here SR and RSR share the same sharp upper bound. Similarly, for the high temperature Curie-Weiss model, for any non-degenerate sequence of Boolean functions, Lemma~\ref{largeentSR} and the comparison between $L^2$-clue and I-clue in \cite[Proposition 3.3]{GaPe} give the bound $O(\delta(\cU_n))$. On the other hand, when proving no RSR for spin systems on $\Z^d$ with Strong Spatial Mixing, we will obtain a quantitative bound \eqref{e.badcube} that is worse than the general $\sqrt{\delta(\mathcal{U}_n)}$ we obtained above. We do not know if the bound $C\sqrt{\delta(\mathcal{U}_n)}$ is ever sharp. See Question~\ref{q.RSRsqrt}.


\section{Factor of IID measures}\label{FFIID}

\subsection{Introduction and a first result} \label{FFIID.intro}

In this section, we investigate sequences of spin systems that converge to finitary factor of IID systems. As this class of measures is one of the most straightforward extensions of product measures, it is an obvious choice trying to understand them. Moreover, some of the Ising measures can also be described in this framework. 

\bde[(Finitary) Factor of IID systems]\label{FFIIDdef}
Let $\cG_*$ be the set of locally finite, possibly edge-labelled, graphs $G$ with a distinguished root $o\in V(G)$. Let $$
\cG_{*,[0,1]}=\left\{(G,o,\omega) : (G,o)\in \cG_*,\ \omega\in [0,1]^{V(G)} \right\}
$$ 
be the set of $[0,1]$ vertex-labelled rooted graphs, equipped with the ``local product topology'': two labelled rooted graphs are close to each other if they are isomorphic in a large neighbourhood of the root, including the edge-labels, and the $[0,1]$ vertex-labels are close in a large sub-neighbourhood of the root. Then let $\Omega$ be a probability measure space, and
 $\psi: \cG_{*,[0,1]} \lora \Omega$ be a Borel-measurable function that is invariant under all the automorphisms of every $G$ fixing the root.
 
For any locally finite graph $G(V,E)$, the \empha{factor of IID (FIID)} random configuration $\sigma\in\Omega^V$ generated by the \empha{coding map} $\psi$ is defined as  
$$
\sigma_v:=\psi(G,v,\omega), \quad v \in V(G),
$$
where $\omega \sim \mathrm{Unif}[0,1]^{V}$ is an underlying IID source field.

%

The coding map $\psi$ is called \empha{finitary}, and the resulting field $\sigma$ is a \empha{finitary factor of IID (FFIID)} if, additionally, there almost surely exists a random coding radius $R=R(G,v,\omega)<\infty$, for which it holds that $\psi(G,v,\omega)$ is determined by $\left\{\omega_u : u \in B_{R}(v) \right\}$, including the value of $R$.

Most often, $\psi$ is defined only for a fixed infinite transitive graph $G$, and one just says that $\sigma$ is given by an equivariant measurable function $\Psi$ of $\omega$: i.e., $\Psi$ commutes with the automorphism group of $G$. The connection between the two descriptions is that $\Psi(\omega)(v) = \psi(G,v,\omega)$. 
\ede

FIID processes, under the name Bernoulli factors, have a rich history in ergodic theory. Most importantly, \cite{OW.isom} proved that every discrete-valued Bernoulli factor over an amenable group is in fact measurably isomorphic to an IID system, and a complete invariant for isomorphism is the Kolmorov-Sinai entropy calculated along F{\o}lner exhaustions. Moreover, every free ergodic process of positive entropy factors onto an IID process with the same entropy; this is true on any countable group \cite{Sew}. Furthermore, it is easy to see that every invariant process on an amenable group can be weakly approximated by FIID processes. On the other hand, on non-amenable groups, there are nice processes that are far from being FIID: a not fast enough two-point correlation decay \cite{LyNaz,BaVi}, or the violation of certain entropy inequalities \cite{BGH,CsHaVi} exclude FIID-ness. 

For discrete-valued process, being a finitary FIID is  equivalent to almost sure continuity of the factor map. Another motivation is computational: a FFIID measure can be generated from an IID measure via a distributed local algorithm that terminates in finite time almost surely \cite{GR}. Of course, from a practical point of view, it is desirable to have some reasonable  bound on the size of $B_{R}(v)$, the number of IID variables that need to be revealed to learn $\sigma_v$. For instance, the condition of \empha{finite expected coding volume}, $\E[\Vol_v] := \E[|B_{R}(v)|] < \infty$, or slightly more strongly, $ \E[|B_{2R}(v)|] < \infty$, will come up for us several times.

Some examples of FIID and FFIID processes are in order. It was proved in \cite{BS}, using positive associations and coupling-from-the-past \cite{PW}, that the plus (and the minus) Ising measures are FIID on any graph; however, on amenable transitive graphs, when they are not equal (i.e., at supercritical inverse temperatures), they are not finitary FIIDs, due to the sub-exponential large deviations. On the other hand, at the critical temperature on $\Z^d$, the unique Gibbs measure is FFIID, but it is not a finitary FIID with finite expected coding volume. 

On non-amenable transitive graphs, the above surface-order large deviations phenomenon for low temperature Ising models does not exist anymore, and in fact, they are FFIID measures, for $\beta$ high enough \cite{RaySpin}. The free Ising measure on the $d$-regular tree $\T_d$ is not a FIID for large enough $\beta$, where the plus and minus measures are the only extremal ones; on the other hand, it is a FIID in part of the intermediate regime where it is extremal, but different from the plus and minus measures; this was proved in \cite{NSZ} by constructing a strong solution to an infinite-dimensional stochastic differential equation. Whether the free Ising measure is ever FFIID in this regime, that is completely open. 

Beyond the Ising model: finitely dependent processes on amenable groups, and Markov random fields on $\Z^d$ with Strong Spatial Mixing (see Definition~\ref{SSM}), or on general graphs with a high noise condition (such as Definition~\ref{ASSM}) are FFIID; see \cite{Spinka.dep} and \cite{Spinka.SSM}, respectively. Stochastic domination, with either positive or negative associations, can also be used on amenable graphs to prove FFIID-ness \cite{Adam.stochdom}; this applies, e.g., to the Uniform Spanning Forest. On non-amenable groups, there are fewer techniques and results: apart from the results mentioned above, some notable examples are the uniform perfect matching measure \cite{LyNaz}, the Wired Uniform Spanning Forest \cite{ARS}, random interlacements \cite{BRR}, and an exotic process with $\mathrm{Unif}[0,1]$ marginals on the vertices, and infinite clusters spanned by equal labels \cite{Peti}. See also \cite{Lyons} for a not very recent survey.
\medskip

The question we investigate here is what conditions on the limiting measure imply that there is or there is no sparse reconstruction (or some of its variant) for some measures $\p_n$ converging to a FFIID measure $\p$. One of the key questions is what notion of convergence is suitable in this context. In light of the examples in Section~\ref{basicconc}, weak convergence is too weak in this setting. Also, we would like a sort of convergence that keeps the factor of IID structure.

One may consider two natural ways to link a sequence of measures to a finitary factor of IID measure on an infinite graph. The first option, already mentioned in the introduction, is to take an exhaustion $V_n$ of $V(G)$ via finite sets, for example by balls, and project the measure $\p$ onto $V_n$ to get $\p_n$. While being a general convergence concept, it still suits the FFIID setting well, since  $\p_n$ can be obtained by applying the translated factor maps to the IID sourse on the infinite graph. However, as we lose the symmetry of the original graph, in this setting we can only talk about RSR. 

Another approach is to introduce a stronger concept of convergence, tailored specifically for FIID measures. We  only consider sequences of finite spin systems which themselves are factors of IID and, in particular, are generated by (a possibly truncated version of) the same local algorithm that we see in the limit.

%

\bde[(F)FIID approximation]\label{d.FIIDconv}
Consider a sequence $G_n$ of finite graphs converging to $G$ in the Benjamini-Schramm sense. We take a sequence $\psi_\eps$ of functions depending only on $\{\omega(u) : u\in B_{r(\eps)}\}$ for some finite non-random $r(\eps)<\infty$, in such a way that $\psi_\eps$ converges to $\psi$ in probability, as $\eps\to 0$. Then let $\eps_n$ converge to 0 so slowly that the proportion of vertices in $G_n$  whose $r(\eps_n)$-neighbourhood is isomorphic to that in $G$ is at least $1-\eps_n$, and let 
$$\sigma_n(v):=\psi_{\eps_n}(G_n,v,\omega_n), \quad\textrm{for } \omega_n \sim \mathrm{Unif}[0,1]^{V(G_n)}\,.
$$
For finitary $\psi$, an even more natural approximation $\psi_n$ is the following. For $v\in V(G_n)$, let $\rho_n(v)$ be the largest radius such that the $\rho_n(v)$-ball in $G_n$ around $v$ is isomorphic to the ball in $G$. If $\omega_n$ on $G_n$ is such that the coding radius $R$ in this ball does not exceed $\rho_n(v)$, then we can just set $\sigma_n(v):=\psi(G_n,v,\omega_n)$; otherwise, we let $\sigma_n(v)$ be an arbitrary value measurable w.r.t.~the $\omega_n$-configuration inside this ball.
\ede

Simple heuristics suggest that if  $\{ (G_n, \p_n): n \in \N \}$ converges to a FFIID measure with finite expected coding volume, there should be no sparse reconstruction. Indeed, for each $v \in U_n$ one can learn the value of $\sigma_v$ by asking  $\E[\Vol_v]$  IID labels on average, which is a finite number. Thus, one uses around $\E[\Vol_v] |U_n|$ IID variables in total, which is much less than the total number of IID variables in the system, so by Theorem~\ref{t.cluegen}, there should be no sparse reconstruction.

The problem is that here the coordinates are queried iteratively, according to an adaptive algorithm, and in this case a transitive function can actually be sometimes determined from a low revealment set of coordinates (see Theorem~\ref{reveal}). 

Nevertheless, in our case the algorithm is localised,  and with high probability it will reveal a deterministic and not much larger subset, the union of balls around the vertices in $U_n$. This idea is enough to bound the $\clue$, under some additional conditions on the behavior of the coding radius $R$.

%

\bpf[Proof of Theorem~\ref{t.almost}]
First observe that if $\sigma$ is an FFIID system on $\Z^{d}$ with exponentially decaying coding radius, then so is the approximating $\sigma_n$ on $\Z^d_n$ from Definition~\ref{d.FIIDconv}, with the same constant $c$.

Let us first take an $n \in \N$ and a positive integer $L_n$ to be determined later, and define the event that all the vertices in $U_n$ have coding radius less than $L_n$:
$$
A_n : = \left\{\forall \, u\in U_n : \; R_u < L_{n}  \right\}\,.
$$
By the union bound, we have the following estimate:
\beq \label{eq.A_nbound}
\p[A_n ] = 1 - \p[ \exists\, u\in U_n : \; R_u \geq L_n] 
\geq 1-|U_n| \exp(-c L_n)\,.
\eeq
By definition, whenever $A_n$ happens, the spins in $U_n$ can be calculated from at most $|U_n| L_{n}^d$ independent, $\mathrm{Unif}[0,1]$ distributed labels. Denote the set of the vertices by $J_n :=  \bigcup_{u\in U_n} B_{L_n}(u)$ and let $\mathcal{J}_n := \sigma \{ \omega(v)  : v \in J_n \}$ be the $\sigma$-algebra generated by the respective labels.

Now we can consider the function $g_n:= f_n \circ \Psi_n$,  where $\Psi_n(\omega)(v)=\sigma_n(v)$ is the equivariant FFIID approximation.  Obviously, whenever $f_n$ is transitive, so is  $g_n$, and since the variables $\omega(u)$ are i.i.d.,  Theorem~\ref{t.cluegen} tells us that 
$$
\clue_{\omega_n}\left(g_n \; | \;J_n\right)\leq \frac{|J_n|}{n^{d}}.
$$
Now we argue that if $\lim_n{\p[A_n ]}=1$, then $\clue_{\omega_n}\left(g_n \, | \, J_n\right)$ is close to $\clue_{\sigma_n}(f_n \, | \, U_n)$, so in case ${|J_n|/n^{d}} \to 0$, we have $\clue(f_n \,| \,U_n ) \to 0$, as well. 

Indeed, let $Z_n : = g_n(\omega_n) = f_n(\sigma_n)$, and introduce the sigma-algebra $\mathcal{G}_n: = \mathcal{J}_n \vee \cF_n(U_n)$, where $\cF_n(U_n)$ is the sigma-algebra generated by the variables $\sigma_n(U_n)$. It is clear that  $\E[Z_n |\mathcal{G}_n] = \E[Z_n|\mathcal{F}_n(U_n)]$.  Using that $\mathcal{J}_n \subseteq \mathcal{G}_n$, by Pythagoras's theorem we have 
$$
\big\| \E[Z_n|\mathcal{G}_n] \big\|^2=\big\| \E[Z_n| \mathcal{J}_n ] \big\|^2 + \big\|  \E[Z_n | \mathcal{G}_n] - \E[Z_n |\mathcal{J}_n ]  \big\|^2\,,
$$
where $\|\cdot\|^2$ denotes the standard $L^2$-norm with respect to $\E$.
Subtracting the common squared expectation $\E^2[Z_n]$, and using that whenever $A_n$ holds (which is a $\mathcal{J}_n$-measurable event) $\E[Z_n | \mathcal{G}_n] - \E[Z_n |\mathcal{J}_n ] = 0$,  we get 
$$
\big| \Var\left(\E[Z_n |\mathcal{G}_n]\right)  - \Var\left(\E[Z_n | \mathcal{J}_n ]\right)  \big| \leq 2 (1 - \p[A_n]).
$$
Using  that $f_n$ is non-degenerate, that is, $\Var(Z_n)$ is bounded away from $0$, we get as promised, that $\clue_{\omega_n}(g_n | J_n) \to 0$ is equivalent to $\clue_{\sigma_n}(f_n | U_n)$ in case  $\lim_n{\p[A_n ]}=1$. 

What is left to be done is to ensure that $\lim_n{\p[A_n ]}=1$. According to \eqref{eq.A_nbound}, this is satisfied whenever
$$
 \displaystyle
  \lim_{n} |U_n|\exp(-c L_n) = 0.
$$
Second, we require that $J_n$ is a sparse subset of coordinates:
$$
|J_n| \leq  |U_n| L_n^d \ll |\Z^d_n| = n^{d}.
$$
It is a straightforward calculation to verify that in case we choose $L_n = K \log n$ with sufficiently large $K$, and $|U_n|\ll (n/ \log n)^{d}$, both of the above hold. 
\epf

It is a classical result from \cite{BS} that the Ising measure in the high temperature regime satisfies the conditions of Theorem~\ref{t.almost}, giving an almost sharp non-reconstruction result. Theorem~\ref{t.norSR} will give the sharp result. 

%

Somewhat surprisingly, the bound of Theorem \ref{t.almost} turns out to be sharp in general, at least in the case when we use the setup based on Benjamini-Schramm convergence (as opposed to an exhaustion):

\begin{ex}\label{e.devanSR}
 We first construct a measure that is a finitary factor of fair IID  $\{-1,1 \}$ bits on ${\Z}$, and a sequence of tori converging to it in the FFIID sense, such that one can reconstruct a sequence of functions with tiny but positive probability from a set of three coordinates. We will then use this to define an FFIID measure on $\Z^2$ and an FFIID approximation that shows that Theorem~\ref{t.almost} is close to being sharp. 
 
 The local algorithm $\psi$ on $\{-1,1 \}^\Z$ decides the spin value $\sigma_{0}$ as follows. Read the bits starting from $0$ going to the positive direction, until one finds two consecutive bits with the same value. This value will be the spin $\sigma_{0}$. 

It is clear that the coding radius of this simple FFIID system decays exponentially. It is also straightforward to check that  almost surely, the strings $- + -$ and  $+ - +$ are not contained in the configuration $\sigma$.  

Now we define a sequence  $\sigma_n$ of  FFIID spin systems on $\Z_{2n}$  that converges to $\sigma$ in the sense of Definition~\ref{d.FIIDconv}. The algorithm $\psi_n$ on $\Z_{2n}$ wants to follow the same rule as $\psi$. However, it happens with probability $1/2^{2n-1}$ that $\psi$ does not stop in $2n$ steps, when the configuration is alternating (that is,  $\sigma_n(i)\sigma_n(i+1) = -1$ for any $i \in \Z_{2n}$), and then we set $\sigma_n(0) = \omega(0)$. In this case, the entire $\sigma_n$ is alternating, and this event can be reconstructed from any $3$ consecutive spins: we can see  $- + -$ or  $+ - +$ if and only if the entire sequence is alternating. 

Of course, this is a highly degenerate sequence of events, as the respective  probabilities decay exponentially. Nevertheless, one can modify this example to obtain a sequence of measures in which a non-degenerate sequence  can be reconstructed from a sparse subset of coordinates.  We define the same algorithm $\psi$ on $\{-1,1 \}^{\Z^2}$, so the algorithm only queries bits to the positive direction along the $x$-axis. Now we define the sequence of FFIID spin systems on the graphs $G_n = \Z_{2n} \times \Z_{2^{2n-1}}$, and we generate the spins with $\psi_n$ on every copy of $\Z_{2n}$.
Clearly, the number of lines  on  which $\psi_n$ does not stop has a distribution close to Poisson$(1)$, and therefore the sequence of events $\{\exists k \in [2^{2n-1}] :\; \sigma_n(\;,k)\; \text{ is alternating} \}$ is non-degenerate, and it can be reconstructed from $3 \cdot 2^{2n-1}$ spins. The density of this subset is not much lower than the upper bound $n^2/\log^2 n$ in Theorem~\ref{t.almost}.
\end{ex}
 
Still, this example feels unsatisfactory. The original idea was to link sparse reconstruction in a sequence of FFIID spin systems  to the properties of the limiting FFIID measure. In this case, however, it is the error of the approximation that we can reconstruct. In order to get around this, \'Ad\'am Tim\'ar suggested the following simple strengthening of the concept of convergence for FFIID spin systems. 

\bde[Regular FFIID approximation]\label{d.Adam}
For a sequence of  FFIID spin systems $\{ (G_n, \p_n) : n \in \N\}$ let  $ \rho_n$ be the largest radius such that for every $v \in V_n$ the ball $B_v(\rho_n)$ is isomorphic to the corresponding  $\rho_n$-ball on $G$.  We say that a sequence of FFIID spin systems $(G_n, \p_n)$ converges to the FFIID $(G, \p)$ \empha{regularly} if the sequence of approximations $\p_n$ satisfies 
$$
\lim_{n \rightarrow \infty} \p_n[ \exists\, v \in V_n :\; R_v > \rho_n] = 0.
$$
That is, we require that the approximation is good enough to ensure that the probability of any error to occur tends to $0$. 
\ede

Using this condition, we restate the main question on FFIID measures in Question~\ref{q.adam}.

\subsection{Versions of tail-triviality}\label{FFIID.tt}

As we explained around Definitions~\ref{d.stber} and~\ref{d.stgen}, giving slightly different notions of sparse tail-triviality, FIID measures can be far from being tail-trivial. So, one can either consider these sparse relaxations, as in Szegedy's Conjecture~\ref{c.SzyB} and Proposition~\ref{p.SzyB}, or a strengthening of the factor of IID property, as in Proposition~\ref{p.fvFFIIDtt} below.

\bpr[Amenable FIID sparse tail-triviality \cite{Szy}]\label{p.SzyB}
Consider a right Cayley graph $G$ of an amenable group $\Gamma$. Let $\sigma$ be a $\pm 1$-valued factor of IID system on $\Gamma$. Then any sparse Bernoulli reconstructable function $f:\{-1,1\}^\Gamma\lora\R$ is constant.
\epr

Note that it is important here to talk about discrete-valued spin systems only, because we will use entropy in the proof.

\bpf
Assume that there is a non-constant sparse Bernoulli reconstructable function. Then there is also such a non-constant Boolean function $f:\{-1,1\}^\Gamma\lora\{0,1\}$, since any superlevel set of the original function is again sparse Bernoulli reconstructable. Now consider the spin system $\zeta_x:=f\big(x^{-1}\sigma\big)$, where $x^{-1}\sigma_{y}:=\sigma_{xy}$ is the $x$-left-translated configuration. Since the Bernoulli subset $\cU_\eps$ from Definition~\ref{d.stber} has the same distribution as its $x$-translate, we also have $\clue( \zeta_x \,|\, \cU_\eps)=1$ almost surely, for every $x$. This means that for almost every realization of $\cU_\eps$, we can obtain the field $\zeta$ as a measurable function of $\sigma |_{\cU_\eps}$; i.e., $\zeta$ is a factor of the FIID spin system $\eta^\eps := \left( \sigma_x \mathbf{1}_{\{x\in \cU_\eps\}}\right)_{x\in V}$, equivariant w.r.t.~$\Gamma$.

On the other hand, $\zeta$ is an FIID spin system with non-trivial marginals, hence it has a positive entropy per site: for any F{\o}lner exhaustion $F_n$ of $\Gamma$,
\beq\label{e.h}
\lim_{n\to\infty} \frac{1}{|F_n|} \, \HH\big(\zeta |_{F_n}\big) = h > 0
\eeq
exists, is positive, and is independent of the F{\o}lner exhaustion, where $\HH(Z)$ is the $\log_2$-entropy of the random variable $Z$. (See \cite{OW.SMB,OW.isom} for the entropy theory of amenable group actions.) However, the entropy per site of the process  $\eta^\eps$ is at most the entropy per site of the random subset $\cU_\eps$, which is almost surely $\HH(\mathsf{Ber}(\eps))=-\eps \log_2 \eps - (1-\eps) \log_2(1-\eps)$, which is smaller than the $h$ of~(\ref{e.h}) if $\eps>0$ is small enough. This contradicts the fact that $\zeta$ is a factor of $\eta^\eps$.
\epf

For FIID processes on non-amenable groups, the situation is completely unclear. It may appear that, instead of the limiting entropy per site along a F{\o}lner sequence, one could easily use sofic entropy along a Benjamini-Schramm approximation, assuming that the graph $G$ admits one (see \cite{Bowen} for definitions and background on sofic entropy). However, when taking the approximation of the FIID process along the sofic approximation $G_n$ of the infinite graph $G$, it may happen that the small errors in the approximation kill the reconstruction algorithm --- hence we do not see how to derive a contradiction. On the other hand, we do not have any counterexamples to Conjecture~\ref{c.SzyB}, either. We do not know if the example of the invariant perfect matching of $\T_3$ by \cite{LyNaz} has a trivial sparse tail, although we expect so. We also expect that there exist sofic approximations $\{G_n\}$ to $\T_3$ such that any approximation $\p_n$ of the invariant perfect matching measure admits random sparse reconstruction for some {\it non-local} function, which of course would not contradict Conjecture~\ref{c.SzyB}.
\medskip

Under a quantitative condition on the factor map for a FIID measure, tail triviality actually holds. This is part~(1) of Theorem~\ref{t.fvFFIID}, which is a good preparation for part~(2), proved in the next subsection. It is quite possible that we are not the first to observe this tail-triviality; similar arguments can be found in \cite{BS, dHSt1, dHSt2}.

\bpr[Tail-triviality for FFIID with ``robustly finite'' expected coding volume]\label{p.fvFFIIDtt} If $G(V,E)$ is an infinite transitive graph and $\p$ is a finitary FIID measure on $\Omega^V$ with the coding radius $R(o)$ for $o\in V$ satisfying $\E |B_{2R(o)}|<\infty$, then $\p$ is tail trivial.
\epr

\bpf Consider a tail-measurable function $f$. For every $\eps>0$, it has an approximation $f_\eps$ that depends only on the $\sigma$-variables in some finite ball $B=B_s(o)$, with $s=s_\eps$. The set of IID variables that $\sigma |_B$ depends on is contained in the finite random set $W:=\bigcup_{x \in B} B_{R(x)}(x)$, where $R(x)$ is the coding radius for $x$. This set $W$ is contained in some ball $B_\rho(o)$, where, by a union bound,
\beq\label{e.rhotail}
\p[ \rho > r] \leq |B| \, \p\big[R(x) > r-s\big]\,. 
\eeq

Consider now the variables $\sigma_y$ along a sphere, $y \in S_r(o):=\{x : d(o,x)=r\}$. Each of them is a function of the IID variables in $B_{R(y)}(y)$. The expected number of vertices $y\in S_r(o)$ for which $W \cap B_{R(y)}(y) \not= \emptyset$, hence $f_\eps(\sigma)$ may not be independent of $\sigma_y$, is at most $|S_r| \, \p\big[ R(y) > r - \rho \big]$.
Hence, the expected number of spins $\sigma_y$ whose coding set intersects $W$ is at most $|B_{s}(o)|$ plus
\begin{align*}
\sum_{r > s} |S_r| \, \p\big[ R(y) > r - \rho \big] 
&\leq \sum_{r > s} |S_r| \, \Big\{ \p\big[ R(y) > r/2-s/2 \big]  + \p\big[ \rho > r/2+s/2 \big] \Big\} \\
&\leq \sum_{r > s} |S_r| \, \Big\{ \p\big[ R(y) > r/2-s/2 \big]  + |B| \, \p\big[ R(x) > r/2-s/2 \big] \Big\},\textrm{ by }\eqref{e.rhotail},\\
&\leq 2 |B| \sum_{r > s} |S_r| \, \p\big[ R(o) > r/2-s/2 \big] \\
&\leq C |B| \sum_{t > 0} |S_{2t+s}| \, \p\big[ R(o) > t \big], \textrm{ with }C\textrm{ depending on }G,\\
&\leq C |B|^2 \sum_{t > 0} |S_{2t}| \, \p\big[ R(o) > t \big] \\
&\leq C |B|^2 \sum_{t>0} |B_{2t}|  \, \p\big[ R(o) = t \big], \textrm{ by summation by parts},\\
&\leq C |B|^2 \, \E |B_{2R(o)}| < \infty, \textrm{ by assumption.}
\end{align*}
Therefore, the set of variables $\sigma_y$ whose coding set intersects the coding set of $\sigma_o$ is finite almost surely. That is, for a large enough (random) radius $r$, there is no such $y$ outside of $B_r(o)$. This clearly implies that $\p$ is tail trivial. 
\epf

%

\subsection{No RSR for magnetisation}\label{FFIID.magn}

In this section we are going to show that the magnetization cannot be reconstructed from a finitary factor of IID sequence converging --- in either of the two senses specified in the beginning of Section~\ref{FFIID.intro}  --- to a finitary factor of IID sequence measure, under two conditions: the first is that of a robust version of finite expected coding volume (the same as in Proposition~\ref{p.fvFFIIDtt}), the second is that magnetization is 'non-degenerate', that is, its variance is not much lower than in the i.i.d. case. 

The following lemma will be useful to control the variance of certain functions. Moreover, it also gives an upper bound for  clue, using the fact that $\Var (P^{\mathcal{U}}[Y_{U\subseteq V}])=\E[\Cov(Y_{\mathcal{U}_1}, Y_{\mathcal{U}_2})]$, with the definition of $P^\cU$ in~(\ref{eq.pagnY}). 

\bl  \label{l.randomperes} 
Let $X$ be a FFIID process with finite coding volume on a transitive graph $G$ with its coding radius $R(o)$ for $o\in V(G)$ satisfying $\E |B_{2R(o)}|<\infty$. Let $V \subseteq V(G)$ be finite and $\mathcal{U}_1, \; \mathcal{U}_2$  independent, equidistributed random subsets of $V$ with revealment $\delta$. Consider a system of random variables $ \{ Y_S : S \subseteq V \}$, such that $Y_S$ is $\mathcal{F}_{S}$-measurable, with $\left|Y_S\right|\leq K $ almost surely. Then
$$
\E\left[|\Cov(Y_{\mathcal{U}_1}, Y_{\mathcal{U}_2} \bigm| \mathcal{U}_1, \mathcal{U}_2  )|\right] \leq 2 \delta K^2 \, \E|\mathcal{U}_1| \, \E |B_{2R(o)}|
$$
\el

\bpf
\begin{align*}
\E\big[\big|\Cov(Y_{\mathcal{U}_1}, Y_{\mathcal{U}_2} \,|\, \mathcal{U}_1, \mathcal{U}_2  )\big|\big] &= \p\big[ \exists\, u \in \mathcal{U}_1,\, v \in \mathcal{U}_2 \,:\,R_{u} + R_{v}\geq d(u, v) \big] \times \\
&\hskip 1 cm \times \E\Big[\left|\Cov\left(Y_{\mathcal{U}_1},Y_{\mathcal{U}_2}\bigm|  \exists\, u \in \mathcal{U}_1,  \; v \in \mathcal{U}_2 \,:\,R_{u} + R_{v}\geq d(u, v) )\right) \right| \Big]\\
& \leq  K^2 \p\big[ \exists\,  u \in \mathcal{U}_1,  \, v \in \mathcal{U}_2 \,:\,R_{u} + R_{v}\geq d(u, v) \big]
\end{align*}
Using the union bound and that $\mathcal{U}_1$ and  $\mathcal{U}_2$ are independent of both the spin system and  one another, we get:
\begin{multline*}
\p\big[ \exists\,  u \in \mathcal{U}_1,  \, v \in \mathcal{U}_2 \,:\,R_{u} + R_{v}\geq d(u, v) \big] \\
\leq \sum_{u \in V} \left\{ \p\left[ u \in \mathcal{U}_1 \right] \sum_{k \geq 0} \p\big[ \exists\,    v \in \mathcal{U}_2 \,:\, d(u, v) = k \big] \, \p\left[ R_{u} + R_{v} \geq k\right] \right\}.
\end{multline*}
As $\{R_{u} + R_{v} \geq k\}$ implies that  $\{R_{u} \geq k/2\}$ or $\{R_{v} \geq k/2\}$ most hold, we can write 
$$
 \p\left[ R_{u} + R_{v} \geq k\right]\leq \p\left[ R_{u} \geq k/2\right] + \p\left[ R_{v} \geq k/2\right] = 2\p\left[ R_{v} \geq k/2\right].
$$
Thus, using summation by parts we can rewrite the inner sum as
\begin{align*}
\sum_{k=0}^{\infty}\p\big[ \exists \,   v \in \mathcal{U}_2 \,:\, d(u, v) =k\big] \, \p\left[ R_{u} + R_{v} \geq k\right] & \\
 & \hskip -5cm \leq 2\sum_{j=0}^{\infty}{\Big(\p[ \exists\,   v \in \mathcal{U}_2 \,:\, d(u, v) =2j] + \p[ \exists\;    v \in \mathcal{U}_2 \,:\, d(u, v) =2j-1]\Big) \, \p\left[ R_{v} \geq j\right]} \\
& \hskip -5cm  = 2 \sum_{j=0}^{\infty}\p\left[ \exists  \,  v \in \mathcal{U}_2 \,:\, d(u, v) \leq 2j\right] \, \p\left[R_{v} =  j\right].
\end{align*}
At the same time, using that the revealment of $\mathcal{U}_2$ is $\delta$, we have the following estimate:
$$
\p\left[ \exists\,    v \in \mathcal{U}_2 \,:\, d(u, v) \leq k\right]\leq |B_k(u)|\,\delta.
$$
Putting all this together we end up with the following bound:
\begin{align*}
\E\left[|\Cov(Y_{\mathcal{U}_1}, Y_{\mathcal{U}_2} \bigm| \mathcal{U}_1, \mathcal{U}_2  )|\right]
&\leq  2  K^2 \delta \sum_{u \in V}\p\left[ u \in \mathcal{U}_1 \right]  \sum_{k=0}^{\infty}|B_{2k}(u)|\p\left[R_{v} =  k \right]\\
& =  2  K^2  \delta \sum_{u \in V}\p\left[ u \in \mathcal{U}_1 \right] \E[|B_{2R_{v}}(v)|] \\
& =  2  K^2 \delta \, \E[|\mathcal{U}_1|] \, \E[|B_{2R_{v}}(v)|],
\end{align*}
proving the lemma.\epf

\begin{rem}\label{r.Peres}
Lemma \ref{l.randomperes} can be thought of as a randomized and generalized version of Theorem 4.3 in \cite{BS}. That theorem states that if $\sigma$  is a  finitary factor of IID measure on $\Z^{d}$ with finite expected coding volume and non-negative correlations, then
$$
\chi(\sigma) = \sum_{j \in \Z^{d}} {\Cov(\sigma_{0}, \sigma_{j}} )\leq C \, \E |B_{R(o)}| < \infty\,.
$$
This is in fact an easy corollary of Lemma \ref{l.randomperes}. Indeed, let $\mathcal{U}$  be a uniformly random singleton from  a finite $V \subset \Z^{d}$. Then, by  Lemma~\ref{l.randomperes},
$$
\frac{1}{|V|}\sum_{j \in V} {\Cov(\sigma_{0}, \sigma_{j}}) \leq 2 \frac{1}{|V|}  \E |B_{2R(o)}|.
$$
Since $V$ is arbitrary, taking the limit by exhaustions shows that 
$$
\chi(\sigma) \leq 2 \E |B_{2R(o)}| \leq 2^{d+1} \E|B_{R(o)}|.
$$
\end{rem}

If we have a FFIID system $\sigma$, then for a positive integer $L$ one can consider an $L$ block factor of IID $\sigma^{L}$ on the same graph that approximates (the distribution of) $\sigma$ as follows. We set $\sigma_v^{L} = \sigma_v$ whenever the local algorithm $\psi$ generating $\sigma$ stops before going outside the ball $B_v(L)$. Otherwise, we sample $\sigma_v^{L}$  according to the distribution of $\sigma_{v}$ conditioned on $B_v(L)$ independent of everything outside of $B_v(L)$. The approximating block factor can also be defined  on a finite graph $G_n$ approximating $G$, as in Definition~\ref{d.FIIDconv}.

The following lemma is an easy extension of Lemma~\ref{l.randomperes}, that we will use to show that, under appropriate conditions, $f(\sigma)$ and $f(\sigma^{L})$ are close to each other.

\bl  \label{l.rperes_error} 
Let $\sigma$ be a FFIID process with finite coding volume on a transitive graph $G$ with its coding radius $R(o)$ for $o\in V(G)$ satisfying $\E |B_{2R(o)}|<\infty$. For some fixed $L \in \N$, we denote by $\sigma^{L}$ the corresponding block factor of IID system.  Let $\mathcal{U}_1$ and $\mathcal{U}_2$  be independent, equidistributed random subsets of $V(G)$ with revealment $\delta$.  Consider a system of random variables $ \{ Y_S : S \subseteq V(G) \}$, such that $Y_S$ is $\mathcal{F}_{S}$-measurable, with $\left|Y_S\right|\leq K $ almost surely. Then
$$
\E\left[\Cov\left(Y_{\mathcal{U}_1}(\sigma)- Y_{\mathcal{U}_1}(\sigma^{L}), Y_{\mathcal{U}_2}(\sigma)-Y_{\mathcal{U}_2}(\sigma^{L})\bigm| \mathcal{U}_1, \mathcal{U}_2 \right)\right] \leq 2 K^2 \delta_{\mathcal{U}_1} \E|\mathcal{U}| \, \E\left[|B_{2R(o)}|\1_{R_v>L}\right].
$$
\el

\bpf
First observe that $Y_{S}(\sigma)-Y_{S}(\sigma^{L})=0$ whenever $R_{u} \leq L$ for all $u \in S$. Thus, we have the following, slightly modified version of the proof of Lemma~\ref{l.randomperes}:
\begin{multline*}
\E\left[\Cov\left(Y_{\mathcal{U}_1}(\sigma)- Y_{\mathcal{U}_1}(\sigma^{L}), Y_{\mathcal{U}_2}(\sigma)-Y_{\mathcal{U}_2}(\sigma^{L}\bigm| \mathcal{U}_1, \mathcal{U}_2 \right)\right]  \\
\leq K^2 \p\big[ \exists\,  u \in \mathcal{U}_1,  \, v \in \mathcal{U}_2 \,:\, R_{u} + R_{v}\geq d(u, v) \text{ and }   R_{u}, R_{v} > L  \big].
\end{multline*}
Continuing along the argument of Lemma~\ref{l.randomperes}, using that $\mathcal{U}_1$ and  $\mathcal{U}_2$ are independent of both the spin system and of each other, we get:
\begin{multline*}
\p\big[ \exists\,  u \in \mathcal{U}_1,  \, v \in \mathcal{U}_2 \,:\,R_{u} + R_{v}\geq d(u, v)  \text{ and }   R_{u}, R_{v} > L  \big] \\
\leq \sum_{u \in V}\p\left[ u \in \mathcal{U}_1 \right] \sum_{k=0} \p\left[ \exists\,    v \in \mathcal{U}_2 \,:\, d(u, v) =k \right] \, \p\left[ R_{u} + R_{v} \geq k, \     R_{u}> L, \ R_{v} > L \right],
\end{multline*}
and 
$$
 \p\left[ R_{u} + R_{v} \geq k,\   R_{u}> L, \ R_{v} > L \right] \leq  2\p\left[ R_{v} \geq \max( k/2, L+1)\right],
$$
thus, using summation by parts gives us
\begin{multline*}
\sum_{k=0}^{\infty}\p[ \exists \,   v \in \mathcal{U}_2 \,:\, d(u, v) =k] \, \p\left[ R_{u} + R_{v} \geq k,\   R_{u}> L, \ R_{v} > L \right]\\
 \leq 2 \sum_{j=0}^{\infty}\Big(\p[ \exists\,    v \in \mathcal{U}_2 \,:\, d(u, v) = 2j-1] + \p[ \exists\,    v \in \mathcal{U}_2 \,:\, d(u, v) = 2j] \Big) \, \p\left[ R_{v} \geq \max( j, L+1)\right] \\
 = 2 \sum_{j=L+1}^{\infty}\p\left[ \exists  \,  v \in \mathcal{U}_2 \,:\, d(u, v) \leq 2j\right] \, \p\left[R_{v} =  j \right].
\end{multline*}
From here we finish the proof as in Lemma \ref{l.randomperes}.
\epf
We are going to introduce a few concepts that we are going to use in the sequel.

We introduce the natural inner product $(f,g) = \E[fg]$ on the space of real functions on the hypercube with respect to the  the uniform measure $\p_{1/2}:=(\frac{1}{2}\delta_{-1} + \frac{1}{2}\delta_{1} )^{\otimes V }$.
 
\begin{definition} [Fourier-Walsh expansion]
For any $f \in L^2( \{-1,1\}^{V}, \p_{1/2})$ and $\omega \in\{-1,1\}^V$,
\begin{equation}
f(\omega)=\sum_{S\subset V} \widehat{f}(S) \chi_S(\omega), \quad \quad \chi_S(\omega):=\prod_{i\in S}\omega_i\quad (\text{and } \chi_S(\emptyset):= 1).
\end{equation}
\end {definition}

The functions $\chi_S$ form an orthonormal basis with respect to the inner product, so Parseval's formula applies: 
$$
 \sum_{S \subseteq V}{\widehat{f}(S)^2} =\left\|f\right\|^2 .
$$
Noting that $\widehat{f}(\emptyset) = \E[f]$, we also have 
\beq \label{VarFourier}
 \Var(f) = \sum_{\emptyset \neq S  \subseteq V}{\widehat{f}(S)^2}.  
\eeq 

It often turns out to be useful to interpret the squared Fourier coefficients $\widehat{f}(S)^2$ as a measure on all the subsets of $V$. This was, for instance, instrumental in studying noise sensitivity~\cite{GPS}.

\bde[Spectral sample]\label{d.Spec}
Let $f \in L^2( \{-1,1\}^{V}, \p_{1/2})$. The \empha{spectral sample} $\Spec_f$ of $f$ is a random subset of $V$ chosen according to the distribution 
$$
\mathbb{P}[\Spec_{f} = S] = \frac{\widehat{f}(S)^2}{\left\|f\right\|^2}, \quad \text{for any} \; S \subseteq V.
$$
\ede

The following fact makes the Fourier-Walsh expansion a powerful tool to quantify $\clue$. For any  function $f: \{-1,1\}^{n} \longrightarrow \R \;$  we have
\beq \label{condexpFour}
\E[f \,|\, \mathcal{F}_T]= \sum_{S\subseteq T} {\widehat{f}(S)\chi_S}\,.
\eeq
In particular, from \eqref{VarFourier} we get:
\beq \label{VarFourierCond}
\Var(\E[f \,|\, \mathcal{F}_T])= \sum_{\emptyset \neq S\subseteq T} {\widehat{f}(S)^2}.
\eeq 
Hence, the notion of clue translates  to the spectral sample language neatly, as follows: 
\beq \label{clueFourier}
\clue(f\,|\,U) =\p [\Spec_{f} \subseteq U \,|\,\Spec_{f} \neq  \emptyset],
\eeq
where we used \eqref{VarFourier} and  \eqref{VarFourierCond}. One of the proofs of Theorem~\ref{t.cluegen} in~\cite{GaPe} started with this interpretation of clue.

The following result can be considered in our context as a natural generalization of the Fourier-Walsh transform on the hypercube for general product measures; see \cite{ES} or \cite[Section 8.3]{OD}. Again, this was already used in~\cite{GaPe}.

\bth [Efron-Stein decomposition]\label{EfronStein}
For any $f\in L^2(\Omega^n,\pi^{\otimes n})$, there is a unique decomposition
$$f=\sum_{S\subseteq [n]} f^{S}\,,$$
where $f^{S}$ is a function that depends only on the coordinates in $S$, and if $S \not\subseteq T$ then $(f^{S},g)=0$ for any $\mathcal{F}_T$--measurable function $g$.
\eth

We are also going to use some results on randomized algorithms. We start by introducing the necessary concepts.

A randomized algorithm $\mathcal{A}$ for a function $\Omega^{V} \longrightarrow \R$ queries the coordinates $\omega_j$ for $j \in V$ one by one, in such a way that the decision of the coordinate to ask next might be made based on the outcome of the values already learned and on external randomness as well.

The revealment of a randomized algorithm for a Boolean function $f$ is the maximum probability that a particular bit is queried during the algorithm. 

\begin{definition} [Revealment]
Let $J_{\mathcal{A}} \subseteq V$ denote the random set of coordinates queried by the algorithm $\mathcal{A}$ until it learns the value of $f$.The revealment of $f$ with respect to $\mathcal{A}$ is 
\begin{equation} 
\delta_f(\mathcal{A}) =\max_{j \in V}{\mathbb{P}[j \in J_{\mathcal{A}}]} 
\end{equation}
\end {definition}

We are going to use the following version from \cite[Theorem 1.8]{SS}: 

\bth [Revealment and spectrum] 
\label{reveal}
Let $\mathcal{A}$ be a randomised algorithm and $f :\Omega^{V} \times \Theta \longrightarrow  \R $ a random function ($\Theta$ denotes the randomness governing $\mathcal{A}$) such that the value of $f$ can be computed from the coordinates revealed by $\mathcal{A}$. Then 
\begin{equation}
\frac{\sum_{\left|S\right| =k}{\left\|{f}^S\right\|_{2}^2}}{\left\|f\right\|_{2}^2}\leq \delta_f(\mathcal{A}) k 
\end{equation}
where $f = \sum_{S \in V} {{f}^S}$ denotes the Efron-Stein decomposition of $f$.
\eth

This is a generalization of the corresponding result from  \cite{SS} in two ways. First, we allow any kind of product space, while the original result is stated only for the uniform hypercube. Second, we allow $f$ to be a random function, while in the original setting $f$ needs to be measurable with respect to the sigma-field generated by $\Omega^{V}$. The steps of the original proof, however, carry over in this more general setting.

We are also going to need the following fact from \cite[Lemma 2.9]{GaPe}:

\bl \label{proj}
Let $f, g \in  L^2(\Omega^V,\p)$ with 
$$
\Corr(f,g)\geq 1 -\eps. $$
Let $U \subseteq V$. If  
$$
\clue(f\,|\,U) \geq c, 
$$
then
$$
\clue(g\,|\,U)\geq  c  - 5\sqrt{\eps}.
$$ 

\el

%
%
%

We are now ready to prove our theorem on having no (R)SR for magnetization in FFIID spin systems with robustly finite expected coding volume, under the technical conditions of having a regular FFIID approximation (Definition~\ref{d.Adam}) and non-degenerate variance.

%

\bpf[Proof of part~(2) of Theorem~\ref{t.fvFFIID}]
First let us fix an $n \in \N$, and set $F = V_n$. 
Let $\mathcal{U} \subseteq F$ be a random subset of the spins with revealment $\delta$. Let $f : \{-1,1\}^{F} \longrightarrow  \R $ be a function of the spins which is a `degree $D$ polynomial'; that is, there exists a  positive integer $D$,  such that whenever $|S|>D$ then $\widehat{f}(S)=0$. This is a slight generalization from the statement of the theorem, where we had $f=M$, the average magnetization, hence $D=1$. Although we will not be able to handle the case $D>1$ completely, we would like to highlight the role of $D$ in the argument. 

Consider an $L$-block factor $\sigma^L$ approximating $\sigma$, as defined before Lemma~\ref{l.rperes_error}, with $L$ to be determined later. Let $\tilde{F}:=\bigcup_{v\in F} B_L(v)$ denote the finite set of IID labels one possibly needs to know to compute the set of spins $\sigma^L_{F}$. Note that the coding volume of $\sigma^L$, by construction, is obviously upper bounded by the coding volume $\Vol$ of the original field $\sigma$.

Let $g$ be a random function that is computed from the spins revealed by $\mathcal{U}$. This will be interpreted as our guess for a $\sigma^L_{F}$-measurable target random variable based on the random subset $\mathcal{U}$. Thus, $g$ is measurable with respect the sigma-algebra jointly generated by the uniform IID labels $\omega_{\tilde{F}}$ and the randomness that samples $\mathcal{U}$. 
%
Observe that $g$ can be computed by a randomized algorithm from the  IID labels $\omega_{\tilde{F}}$ by first sampling $\mathcal{U}$ and after using the local algorithm to compute $\sigma_{u}$ for each $u \in \mathcal{U}$. It is clear that the revealment of this algorithm is at most $\delta \E[\Vol]$. Therefore, Theorem~\ref{reveal} tells us that, for any $S \subseteq \tilde{F}$, and positive integer $k$,
\beq \label{eq.SS}
\sum_{\left|S\right| =k}{\left\| g^S\right\|_{2}^2}\leq k \,\delta \, \E[\Vol]\left\| g\right\|_{2}^2,
\eeq
where $g^S$ denotes the corresponding term from the Efron-Stein decomposition of $g$. 

We will use this to give an upper bound on $\clue (Y\,|\,\mathcal{U})$, where  $Y := f(\sigma^{L})$.
We first note that in the Efron-Stein decomposition of $Y$ with respect to $\tilde{F}$ we have $Y^{S} = 0$ (as a function) whenever $|S|>D |B_{o}(L)|$. Indeed, this follows directly from the condition on $f$ and that $Y$ uses at most $|B_{o}(L)|$ IID labels to compute a `spin'.

It will be useful to assume in the sequel that both $Y$ and $g$ are centered, and thus $\Var(Y) = \left\| Y \right\|^{2}_{2}$ and $\Var(g) = \left\| g \right\|^{2}_{2}$ hold. It is clear that both the clue and the behavior of a randomised algorithm are unaffected by translation with a constant.
 
By the law of total covariance, using that $Y$ is independent of $\mathcal{U}$, we have
$$
  \Cov(Y, g) = \E\left[\Cov(Y, g \,|\,\mathcal{U} )\right] . 
$$
Now let us express $\Cov(Y, g)$ in terms of  the Efron-Stein decomposition of $f$ and $g$ on the coordinate set $\tilde{F} \cup \{\mathcal{U}\}$. Recall that $Y$ does not depend on  $\mathcal{U}$ and its  Efron-Stein decomposition is supported on sets not greater than $D |B_{o}(L)|$. Therefore, 
\beq  \label{eq.SS2}
\begin{aligned}
\Cov(Y, g) = & \sum_{\emptyset \neq S \subseteq \tilde{F}  \;: \; |S|\leq D |B_{o}(L)| }{\E[Y^S g^S]} \\
\leq &\left\| Y\right\|_{2} \sqrt{\sum_{0<\left|S\right|\leq D |B_{o}(L)|}{\left\| g^S\right\|_{2}^2}}\\ \leq & \left\| Y\right\|_{2}\sqrt{\sum_ {j=1}^{D |B_{o}(L)|}j \delta \E[\Vol]}  \left\| g \right\|_{2}\\
 \leq  & \left\| Y\right\|_{2}\left\| g \right\|_{2} D |B_{o}(L)|\sqrt{\delta \E[\Vol]}\,,
\end{aligned}
\eeq
where we used the Cauchy-Schwarz inequality and \eqref{eq.SS}. 

Now we will use the estimate \eqref{eq.SS2} to bound the average clue. For $U \subseteq F$, we can assume that $g$, given $\cU=U$, is just $\E[Y \,|\,\mathcal{F}_{U} ]$,  since this is clearly the best possible guess in the $L^2$ sense. 

We will need the fact that $\left\| g \right\|_{2}\leq \left\| Y \right\|_{2} $. This follows from the law of total variance for $g$:
$$
\Var(g) = \E[\Var(g\,|\, \mathcal{U})] + \Var(\E[g \,|\,\mathcal{U}]).
$$
On the one hand, $\Var(g\,|\, {U}) = \Var(\E[Y \,|\, \mathcal{F}_{U}])\leq \Var(Y)$ for any $U \subseteq F$, so the first term is at most $\Var(Y)$. On the other, $\E[g \,|\, \mathcal{U}] = \E\left[  \E[Y \,|\, \mathcal{F}_{\mathcal{U}}]  \right] = \E[Y]$ and therefore the second term vanishes. 

This leads to the following estimate for the average clue:
\beq \label{eq.noSRforF}
\E\left[\clue(Y \,|\,\mathcal{U})\right] = \frac{\E\left[\Cov(Y, \E[Y \,|\, \mathcal{F}_{\mathcal{U}}] \,|\,\mathcal{U})\right]}{\Var(Y)} =  \frac{\Cov(Y, g)}{\Var(Y)} \leq  D |B_{o}(L)|\sqrt{\delta \E[\Vol]}. 
\eeq
In case $\E[\Vol] < \infty$ and $\delta \to 0$, clearly $\E\left[\clue(Y \,|\,\mathcal{U})\right] \to 0$ as well.

The next step is to show that under certain conditions on the variance and/or the polynomial structure of the respective functions,  $Y=f(\sigma^L)$ and $ f(\sigma)$ are close to  
one another with respect to $L^2$ distance. This will be done using Lemma~\ref{l.rperes_error}. 

Indeed, any function $f : \{-1,1\}^{F} \longrightarrow  \R $ can be written in the form $f = P^{\mathcal{S}}[Y_{S \subseteq F}]$ for a suitable choice of random subset $\mathcal{S} \subseteq F$ and corresponding collection  $\{Y_{S} : S \subseteq F\}$ of random variables, where  $Y_{S}$ is $\mathcal{F}_{S}$ measurable. One  way to do this is to take the Fourier-Walsh decomposition of $f$ and set $Y_S : = \sign (\widehat{f}(S))\big(\sum_{T \subseteq F}{|\widehat{f}(T)|}\big) \chi_S$ and let $\p[S \in \mathcal{S} ] = |\widehat{f}(S)|/\sum_{T \subseteq F}{|\widehat{f}(T)|}$. When the spectrum of $f$ is supported on subsets of size no greater than $D$, this $\mathcal{S}$ also has this property.

Now observe that
\beq
\begin{aligned} \label{eq.explainlemma}
 \Var\left(f(\sigma^{L}) - f(\sigma)\right) &=  \Var(P^{\mathcal{S}}[Y_{S \subseteq F}(\sigma)]-P^{\mathcal{S}}[Y_{S \subseteq F}(\sigma^{L})])\\
& =\Var(P^{\mathcal{S}}[Y_{S \subseteq F}(\sigma) - Y_{S \subseteq F}(\sigma^{L}) ])\\
&= \E\big[\Cov\left(Y_{\mathcal{S}_1}(\sigma)-Y_{\mathcal{S}_1}(\sigma^L), Y_{\mathcal{S}_2}(\sigma)-Y_{\mathcal{S}_2}(\sigma^L)\,\big|\, \mathcal{S}_1,\mathcal{S}_2 \right)\big]. 
\end{aligned}
\eeq
For the last equality, we used \eqref{eq.avgcov}. Recalling that $ \E[\Vol]<\infty$, from~(\ref{eq.explainlemma}) we obtain via Lemma~\ref{l.rperes_error} that, for any given $\eps>0$, there exists large enough $L$ such that 
\beq\label{eq.goapro}
\begin{aligned}
\Var(f(\sigma) - f(\sigma^{L}) ) 
 &\leq  C K^2 \delta_{\mathcal{S}} \E[|\mathcal{S}|] \,\E[\Vol\1_{R_v>L}]\\
 &\leq \eps K^2  \E[|\mathcal{S}|] \, \delta_{\mathcal{S}},
\end{aligned}
\eeq
where $K = \max \{|Y_{\mathcal{S}}| \}$.

Let us consider a convergent sequence of FFIID measures on $\{-1,1\}^{V_n}$,  a  sequence of functions $\{f_n : \{-1,1\}^{V_n} \longrightarrow  \R \,:\, n \in \N \}$   with the  corresponding random sets $\mathcal{S}_n$ and random variables  $\{Y^n_{S} \;:\; S \subseteq V_n \}$ such that $Y^n_{S}$ is $\mathcal{F}_{S}$ measurable and   $f_n = P^{\mathcal{S}_n}[Y_{S \subseteq V_n}] $. 

The point is that whenever for a sequence of functions we have  
\beq  \label{eq.impbo}
\Var(f_n) \geq \beta K_n^2  \E[|\mathcal{S}_n|] \delta(\mathcal{S}_n)
\eeq 
for some uniform constant $\beta>0$, then~\eqref{eq.goapro}  translates to the bound  
\beq \label{eq.impb}
  \Var(f_n(\sigma) - f_n(\sigma^{L}) )\leq \eps\Var(f_n),
\eeq   
which promptly implies that $f_n(\sigma)$ and $f_n(\sigma^{L})$ are highly correlated.
Indeed, from the triangle inequality it follows that 
$\D(f_n(\sigma^{L})) \geq \D(f_n(\sigma))(1- \sqrt{\eps})$. Thus, using that  
$$
\D(f)\D(g)  2 (1 - \Corr(f, g ) \leq \Var(f- g)
$$ 
for any $f,g \in L^2$ whenever \eqref{eq.impb} holds,  for any small $\eps >0$ we get that
\begin{align*}
  2 \left( 1 - \Corr\left(f_n(\sigma), f_n(\sigma^{L})\right)\right) \leq  \frac{\Var(f_n(\sigma) - f_n(\sigma^{L})) }{ \D(f_n(\sigma))\D(f_n(\sigma^{L}))} \leq \frac{\eps\Var(f_n(\sigma))}{ \D(f_n(\sigma))\D(f_n(\sigma^{L}))} \leq \frac{\eps}{1 - \sqrt{\eps}}\leq 2 \eps,
\end{align*}
and thus
\beq  \label{eq.varER} 
\Corr\left(f_n(\sigma), f_n(\sigma^{L})\right) \geq 1-\eps.
\eeq
Our claim now follows from Lemma \ref{proj}, which states that if two random  variables are highly correlated, then the respective clues (with respect to the same subset of spins)  need to be close to each other, as well. Suppose $ f_n(\sigma_{V_n})$ admits an average clue higher than $c>0$ from $\mathcal{U}_n \subseteq V_n$. Since the $\clue$ is at most $1$, this would imply that the probability of $\clue(f_n(\sigma)\,|\,\mathcal{U}_n  ) > \tilde{c}$ is uniformly positive, say greater than $p>0$. Then  one can choose an appropriate $\eps>0$ and a corresponding $L \in \N$ in such  a way that
$$
\E [\clue(f_n(\sigma^{L})\,|\,\mathcal{U}_n )] \geq \tilde{c}p/2, 
$$
which is in contradiction with \eqref{eq.noSRforF}. 

Now that we outlined the general argument, we are left to verify  
that the lower bound in~\eqref{eq.impbo} for $\Var(f_n)$ is satisfied and thus Lemma~\ref{l.rperes_error} allows us to conclude that $f_n(\sigma)$ and $f_n(\sigma^{L})$ are highly correlated. 

Indeed, for the average magnetization $f_n=M_n$, we can take a uniformly random singleton  for $\mathcal{S}$; that is, $\p[\mathcal{S} = \{v\}] =1/|V_n|$  and $Y_{\{v\}} = \sigma_{v}$, for every $v \in V_n$. Clearly, $M_n = P^{\mathcal{S}}[Y_{S\subseteq V_n}]$, so in this setting we have $\E[|\mathcal{S}|] = 1$ and $K_n = 1$ and $ \delta_{\mathcal{S}} = 1/|V_n|$,  
hence~\eqref{eq.impbo} is satisfied whenever $\Var(M_n)\geq \beta/n$, for some uniform $\beta >0$.
\epf

\begin{rem}\label{r.positive}
A shortcoming of this result is that it uses the properties of the specific sequence to conclude that there is no RSR for magnetization, instead of relying only on the properties of the limiting factor of IID measure. We do not think that this is completely possible, but here is another condition on the approximating sequence that some might find more natural. Namely, assume for the sequence of FFIID measures converging to our measure either in the Benjamini-Schramm sense or by exhaustions that the correlations between pairs of spins are non-negative. Then, if the limiting FFIID measure has robustly finite expected coding volume, then Remark~\ref{r.Peres} says that the susceptibility is finite, and then the pointwise convergence and the non-negativity of correlations imply that the susceptibilities in the finite systems remain bounded, while the lower bound on the variance of the magnetization is also satisfied, and  Theorem~\ref{t.fvFFIID}~(2) applies. 
%
\end{rem}
 
\begin{rem}
The above argument implies no RSR for many more sequences of functions. In fact, any sequence will work with a uniformly bounded Fourier-Walsh expansion (that is, every function $f_n$ in the sequence satisfies $\widehat{f_n}(S) = 0$ whenever $|S|>D$) and --- this is the more problematic condition --- \eqref{eq.impbo} must hold. While this translates to a natural lower bound on the variance in the case of magnetization, it does not in general.

\end{rem}

\section{Glauber dynamics and strong spatial mixing} \label{s.glauber}

\subsection{Sparse reconstruction and block dynamics}
\label{ss.gap}

We  continue working with a rather general setup:  $V$ and $\Omega$ are finite sets, $\p$ is a probability measure on $\Omega^{V}$, and  $\mathcal{U}$ is a random subset of $V$, independent of $\p$. In case we want to be able to discuss SR, we assume further that there is a group $\Gamma$ acting on $V$ transitively and that $\p$ is $\Gamma$-invariant.

Then $\mathcal{U}$ induces a Markovian  dynamics $(X^{\mathcal{U}}_t)_{t\in \N}$ on $\Omega^{V}$ as follows. At each step, one samples a subset of $V$ according to the law of $\mathcal{U}$, and updates the configuration  $\Omega^{V\setminus \mathcal{U}}$ according to the measure $\p^{ \mathcal{U}}$ on $\Omega^{V\setminus \mathcal{U}}$ induced by the boundary condition  $\Omega^{ \mathcal{U}}$, independently of the previous configuration on $V\setminus \mathcal{U}$. It is easy to see that this chain, called the heat-bath block dynamics for $\p$ induced by $\cU$, with transition operator denoted by $P_\cU$, is reversible (see \cite[Section 6.1]{PGG} or \cite[Section 1.6]{LPW} for the definition), and in case the revealment is $\delta_{\mathcal{U}} < 1$, the unique stationary measure of $X^{\mathcal{U}}_t$ is $\p$. 

A simple, but crucial observation is that the transition operator of the block dynamics coincides with the operator $P_\cU$, as defined in \eqref{eq.pagn}. Indeed, it is clear from the construction of $(X^{\mathcal{U}}_t)$ that, for any $\cF_{V}$-measurable function $f$ and any $\sigma \in \Omega^{V}$  with $X^{\mathcal{U}}_0 = \sigma$, we have
$$
\E_{\sigma}[f(X^{\mathcal{U}}_1)] =\sum_{S \subseteq V} {\p[\mathcal{U} = S] \,
\E[f(X) \,|\, \mathcal{F}_S](\sigma)},
$$
where $\E_{\sigma}$ denotes the expectation given  $X^{\mathcal{U}}_0 = \sigma$.

We introduce a few notions used in the theory of discrete time Markov chains. First, if $\mu$ and $\nu$ are two probability measures on $\Omega^{V}$, their total variation distance is defined as
$$
\|\mu - \nu \|_{TV} = \frac{1}{2}\sum_{x \in \Omega^{V}}{|\mu(x) - \nu(x) |}.
$$

It is a classical result that every irreducible Markov chain on a finite set admits a unique stationary distribution $\mu$, and if we run the chain long enough, it will converge (in any reasonable distance between measures, at least in the Ces\`aro sense) to this $\mu$. In particular, for an $\eps>0$, the mixing time $T_{\mathrm{mix}}(\eps)$ of an irreducible Markov chain with transition matrix $P$ is the minimal $t \in \N$ such that, for every $x \in \Omega^{V}$, we have
$$
\|\mu - P^t(x,\cdot) \|_{TV} \leq \eps,
$$
where $P^t(x,\cdot)$ is the distribution after $t$ steps of the chain starting from $x$.

The relaxation time $T_{\mathrm{rel}}$ of the chain is $\frac{1}{1-\lambda_{\ast}}$, where $\lambda_{\ast}$ is the second highest absolute value among the eigenvalues of $P$. For reversible chains it is known that the transition matrix has only real eigenvalues. From now on we will talk only about reversible chains.

The Dirichlet form of a function $f$ on $\Omega^{V}$ is defined as
$$
\mathcal{E}(f) : = \E_{\mu}[f (f-Pf)] = \frac{1}{2}\sum_{x, y \in \Omega^{V}}{\mu(x)P(x,y) (f(x)-f(y))^2}.
$$


\bl \label{l.eigenclue}
Let $V$ and $\Omega$ be finite sets, a probability measure $\p$ on $\Omega^V$, a random subset $\mathcal{U}\subseteq V$ independent of $\p$, and the corresponding heat-bath block dynamics $X^{\mathcal{U}}_t$ with transition matrix  $P_{\mathcal{U}}$. Then
$$
\max_{f \in L^2(\Omega^{V},\p)} {\E[\clue_{\p}(f \; |\;\mathcal{U})] =\lambda_{2}(P_{\mathcal{U}})} ,
$$
where  $\lambda_{2}$ denotes the second largest eigenvalue of $P_{\mathcal{U}}$.
\el

Since this lemma is of fundamental importance in that it connects maximal clue to the more familiar notion of the spectral gap $1-\lambda_2$, we are going to give two proofs.

\bpf[First proof.] It is well-known (see \cite[equations (13.3) and (13.4)]{LPW}) that, for any reversible Markov chain $(X_t)_{t\ge 0}$ started from the stationary measure $\p$, we have  
\beq \label{dirclue}
\gamma = \min_{f \in L^2(\p, \Omega^{V})}{\frac{\mathcal{E}(f)}{\Var(f)}}= \min_{f \in L^2(\p, \Omega^{V})}{\frac{\frac{1}{2}\E[(f(X_0) - f(X_1))^2]}{\Var(f(X_0))}},
\eeq
where $\gamma := 1 - \lambda_{2}$ is the spectral gap of the chain. Let $ f \in L^2(\p, \Omega^{V})$. We can write
$$
\frac{1}{2}\E[(f(X^{\mathcal{U}}_0) - f(X^{\mathcal{U}}_1))^2] 
= \sum_{S \subseteq V} \p[\mathcal{U} = S] \frac{1}{2}\E\left[(f(X^{\mathcal{U}}_0) - f(X^{\mathcal{U}}_1))^2 \,\big|\,\mathcal{U} = S\right].
$$
We observe that
\beq \label{eq.diricondvar}
\frac{1}{2}\E[(f(X^{\mathcal{U}}_0) - f(X^{\mathcal{U}}_1))^2 \,|\,\mathcal{U} = S] = \E\left[\Var(f(X) \,|\, \mathcal{F}_S )\right],
\eeq
where $X$ is distributed according to $\p$. Indeed, the spin systems $X^{\mathcal{U}}_0$ and $X^{\mathcal{U}}_1$ are both distributed according to $\p$. Moreover,  conditioned on the event $\{\mathcal{U} = S \}$,  we have    $X^{\mathcal{U}}_0(S) = X^{\mathcal{U}}_1(S)$ while $X^{\mathcal{U}}_0(V \setminus S)$ and $X^{\mathcal{U}}_1(V \setminus S)$ are independent conditionally on $\mathcal{F}_S$.

Thus, using the law of total variance, we conclude that
\begin{align*}
\frac{1}{2}\E[(f(X^{\mathcal{U}}_0) - f(X^{\mathcal{U}}_1))^2] 
&= \sum_{S \subseteq V} {\p[\mathcal{U} = S] \, \E\left[ \Var(f(X) \,|\, \mathcal{F}_S ) \right] } \\
&= \sum_{S \subseteq V} {\p[\mathcal{U} = S] \, \big( \Var(f(X)) -  \Var(\E[f(X) \,|\, \mathcal{F}_S])\big)  } \\
&=  \Var(f(X)) \left( 1 - \E[\clue_{\p}(f \; |\;\mathcal{U})]\right).
\end{align*}
So for arbitrary $f \in  L^2(\p, \Omega^{V})$ we have
$$
\frac{\frac{1}{2}\E[(f(X_0) - f(X_1))^2]}{\Var(f(X_0))} = 1 - \E[\clue_{\p}(f \; |\;\mathcal{U})],
$$
and taking the minimum over all $f \in  L^2(\p, \Omega^{V})$ we get the statement of the lemma. \epf

\bpf[Second proof.] The key ingredient is the following inequality, reminiscent to Lemma~\ref{corr_lemma}.

\bl \label{corr_lemma2} Let $X_V = \{ X_v : v \in V \}$ be a system of random variables on a common probability space. Let  $Y$ be  an $X_V$-measurable random variable with finite second moment and $\mathcal{U}$ a random subset of $V$, independent of $X_V$. Then
$$
\Corr^2(Y, P^{\mathcal{U}}[Y]) \cdot \E\left[ \Corr^2\big(P^{\mathcal{U}}[Y], \E[Y \,|\,\mathcal{F}_{\mathcal{U}}] \;\big|\;\mathcal{U}\big) \right] \geq \E\left[\Corr^2\big(Y,\E[Y\,|\, \mathcal{F}_{\mathcal{U}}] \;\big|\;\mathcal{U}\big)\right].
$$ 
\el

\bpf
For the first factor on the left hand side, observe that
\beq \label{eq.cov2_1}
\begin{aligned}
\Cov \left(Y, P^{\mathcal{U}}[Y]\right) 
&= \sum_{S \subseteq V} {\p[\mathcal{U}=S] \, \Cov(Y, \E[Y| \mathcal{F}_{S}])}\\
&= \E[\Cov(Y,\E[Y |\mathcal{F}_{\mathcal{U}}] \,|\,\mathcal{U})]
= \E\left[\Var(\E\left[Y| \mathcal{F}_{\mathcal{U}}\right]\,|\,\mathcal{U})\right].
\end{aligned}
\eeq
For 
the last equality we used that conditional expectation is an orthogonal projection.
 
Now we estimate the second factor of the left hand side. First, quite similarly to~\eqref{eq.cov2_1},
\begin{equation}\label{eq.corr2_1}
\begin{aligned}
\Cov\big(P^{\mathcal{U}}[Y], \E[Y| \mathcal{F}_{\mathcal{U}}] \,\big|\, \mathcal{U} =T\big) 
&= \sum_{S \subseteq V} {\p[\mathcal{U}=S] \, \Cov(\E[Y| \mathcal{F}_{S}], \E[Y| \mathcal{F}_{T}])}\\
&= \E\left[\Cov\left(\E[Y| \mathcal{F}_{\mathcal{U}_1}], \E[Y| \mathcal{F}_{T}]\,\big|\,\mathcal{U}_1\right)\right],
\end{aligned}
\end{equation}
and therefore
$$
\E\left[ \Corr^2\big(P^{\mathcal{U}}[Y], \E[Y \,|\,\mathcal{F}_{\mathcal{U}}] \;\big|\;\mathcal{U}\big) \right]
=\frac{1}{\Var(P^{\mathcal{U}}[Y])}\E \left[\frac{\E^2\left[\Cov\big(\E[Y |\mathcal{F}_{\mathcal{U}_1}], \E[Y| \mathcal{F}_{\mathcal{U}_2}]\,|\, \mathcal{U}_1,\mathcal{U}_2\big) \,\big|\, \mathcal{U}_2\right]}{\Var\big(\E[Y| \mathcal{F}_{\mathcal{U}_2}]\,\big|\,\mathcal{U}_2\big)}\right].
$$
Without loss of generality we may assume that $\E[Y] = 0$, so  we can write $\Var\left(\E[Y| \mathcal{F}_{\mathcal{U}_2}]\,|\,\mathcal{U}_2\right) = \E^2 \left[ \E[Y| \mathcal{F}_{\mathcal{U}_2}]\,|\,\mathcal{U}_2 \right]$.
Now we get the following lower bound:
\begin{equation} \label{eq.corr2_2}
\begin{aligned}
\E \left[\frac{\E^2\left[\Cov\big(\E[Y |\mathcal{F}_{\mathcal{U}_1}], \E[Y| \mathcal{F}_{\mathcal{U}_2}]\,\big|\, \mathcal{U}_1,\mathcal{U}_2\big) \,\big|\, \mathcal{U}_2\right]}{\Var\big(\E[Y| \mathcal{F}_{\mathcal{U}_2}]\,\big|\,\mathcal{U}_2\big)}\right]&\\
&\hskip -3cm \geq \frac{\E^2 \left[ \E\left[\Cov\big(\E[Y |\mathcal{F}_{\mathcal{U}_1}], \E[Y| \mathcal{F}_{\mathcal{U}_2}]\,|\, \mathcal{U}_1,\mathcal{U}_2\big)\,\big|\, \mathcal{U}_2\right]\right]}{\E\left[\Var\left(\E\left[Y| \mathcal{F}_{\mathcal{U}}\right]\,\big|\,\mathcal{U}\right)\right]}\\ 
&\hskip -3cm = \frac{\E^2 \left[\Cov\big(\E[Y| \mathcal{F}_{\mathcal{U}_1}], \E[Y| \mathcal{F}_{\mathcal{U}_2}]\,|\,{\mathcal{U}_1}, \mathcal{U}_2\big)\right]}{\E\left[\Var\big(\E\left[Y| \mathcal{F}_{\mathcal{U}}\right]\,\big|\,\mathcal{U}\big)\right]}\\
&\hskip -3cm = \frac{\Var(P^{\mathcal{U}}[Y])^2}{\E\left[\Var\big(\E\left[Y| \mathcal{F}_{\mathcal{U}}\right]\,\big|\,\mathcal{U}\big)\right]},
\end{aligned}
\end{equation}
where we first used that, by the Cauchy-Schwarz inequality, $\E[X^2/Y^2] \geq\E^2[X]/\E[Y^2]$, then we used the tower property, and for the final line we recalled~\eqref{eq.avgcov}, noting that $\mathcal{U}_1$  and $\mathcal{U}_2$ are independent.
So, we altogether get
$$
\E[\Corr^2(P^{\mathcal{U}}[Y], \E[Y |\mathcal{F}_{\mathcal{U}}]\,|\,\mathcal{U})] \geq \frac{\Var(P^{\mathcal{U}}[Y])}{\E\left[\Var(\E\left[Y| \mathcal{F}_{\mathcal{U}}\right]\,|\,\mathcal{U})\right]}.
$$
From this and~\eqref{eq.cov2_1}, we conclude that
$$
\Corr^2(Y, P^{\mathcal{U}}[Y]) \cdot \E\left[ \Corr^2\big(P^{\mathcal{U}}[Y], \E[Y \,|\,\mathcal{F}_{\mathcal{U}}] \;\big|\;\mathcal{U}\big) \right]
\geq  \frac{\E\left[\Var(\E\left[Y| \mathcal{F}_{\mathcal{U}}\right]\,|\,\mathcal{U})\right]}{ \Var(Y)}.
$$
Which is exactly what we wanted, since the right hand side of our claim can be written as:
$$
\E\left[\Corr^2\big(Y,\E[Y\,|\, \mathcal{F}_{\mathcal{U}}] \;\big|\;\mathcal{U}\big)\right] 
= \E\left[\frac{\Cov^2(Y,\E[Y |\mathcal{F}_{\mathcal{U}}] \,|\,\mathcal{U})}{\Var(Y) \Var(\E\left[Y| \mathcal{F}_{\mathcal{U}}\right]\,|\,\mathcal{U})}  \right] = \E\left[\frac{\Var(\E\left[Y| \mathcal{F}_{\mathcal{U}}\right]\,|\,\mathcal{U})}{\Var(Y) }  \right],
$$
finishing the proof.\epf

It follows directly from  Lemma \ref{corr_lemma2}
 that, for every random variable $Y$, we have
\beq \label{eq.PisGood}
\E[\clue(Y\,|\,\mathcal{U})] \leq \E[\clue( P^{\mathcal{U}}[Y]\,|\,\mathcal{U})].
\eeq
Moreover, when applying $P^{\mathcal{U}}$ to any random variable $Y$, the clue clearly increases whenever $\Corr(Y, P^{\mathcal{U}}[Y])<1$. So we have equality in \eqref{eq.PisGood} if and only if $Y$ is an eigenfunction of $ P^{\mathcal{U}}$.

As $P^{\mathcal{U}}$ is self--adjoint, by the spectral theorem there is an eigenbasis, and in particular  $\lim_{n \to \infty} (P^{\mathcal{U}})^n[Y] $ exists and is an eigenfunction of  $P^{\mathcal{U}}$. This shows that $\clue$ is maximised by some eigenfunction of $P^{\mathcal{U}}$.

 On the other hand, we can explicitly calculate $\clue_{\p}(f \,|\,U)$ of an  eigenfunction $f$ with eigenvalue $\lambda$. Let $Z = f(X_{V})$. Then
\begin{align*}
\E[\clue(f\,|\,\mathcal{U})] &= \frac{\E\left[\Cov(Z, \E[Z\,|\, \mathcal{F}_{\mathcal{U}}]\,|\,\mathcal{U}) \right]}{\Var (Z)}\\
&= \frac{\lambda^{-1}\E\left[\Cov( P^{\mathcal{U}} [ Z], \E[Z\,|\, \mathcal{F}_{\mathcal{U}}]\,|\,\mathcal{U})\right]}{\lambda^{-2}\Var ( P^{\mathcal{U}} [ Z])} = \lambda\,,
\end{align*}
using that 
$$
\Var ( P^{\mathcal{U}} [ Z]) = \sum_{S,T \subseteq V} {\p[\mathcal{U}=S]\p[\mathcal{U}=T] \, \Cov(Y_{S}, Y_{T})} = \E\left[\Cov( P^{\mathcal{U}} [ Z], \E[Z\,|\, \mathcal{F}_{\mathcal{U}}]\,|\,\mathcal{U})\right].
$$
This completes the second proof of Lemma~\ref{l.eigenclue}. \epf

\begin{rem}
If $\mathcal{H}$ is a random subset satisfying  $\mathcal{H} \preceq \mathcal{U}$, that is, there exists a coupling between  $\mathcal{H}$ and $\mathcal{U}$ such that $\mathcal{H} \subseteq \mathcal{U}$, almost surely, then we trivially have $\E[\clue(f \; |\;\mathcal{H})] \leq \E[\clue(f \; |\;\mathcal{U})] = \lambda_2(P_{\mathcal{U}})$ as well.
\end{rem}

\begin{rem}
Lemma \ref{l.eigenclue} can be generalised for a broader class of block dynamics. Let  $\mathcal{U}$ be an independent random subset and $Y^{\mathcal{U}}_t$ a Markov chain with stationary distribution $\p$ that at each step samples a subset from $\mathcal{U}$ and it updates the configuration  $\Omega^{V\setminus \mathcal{U}}$ according to the measure $\p^{ \mathcal{U}}$ on $\Omega^{V\setminus \mathcal{U}}$ induced by the boundary condition  $\Omega^{ \mathcal{U}}$ (but possibly not independently from the previous values on $V\setminus \mathcal{U}$ --- this is where it differs from the heat bath block dynamics above).

Whenever $Q_{\mathcal{U}}$ is the transition matrix of such a dynamics, it holds that: 
$$
\max_{f \in L^2(\Omega^{V},\p)} {\E[\clue_{\p}(f \; |\;\mathcal{U})] }\leq \lambda_{2}(Q_{\mathcal{U}}).
$$
This can be proven in the same way as Lemma \ref{l.eigenclue}, the only difference being that equation \eqref{eq.diricondvar} becomes an inequality: $\frac{1}{2}\E[(f(X^{\mathcal{U}}_0) - f(X^{\mathcal{U}}_1))^2 \,|\,\mathcal{U} = S] \leq\E\left[\Var(f(X) \,|\, \mathcal{F}_S ))\right]$.
\end{rem}

\medskip

As we are about to see, in our applications instead of $\lambda_{2}$ we typically bound the eigenvalue with the second largest absolute value $\lambda_{\ast}$,  by showing that the relaxation time  of $X^{\mathcal{U}}_t$ is small. 

First we give an interesting consequence of Lemma \ref{l.eigenclue}. We need the following concept. 

Let $P$ be the transition matrix of an aperiodic irreducible Markov chain on $\Omega^{V}$ with stationary measure $\mu$ and let $E \subset \Omega^{V}$ such that $\mu(E) \leq 1/2$. Then the bottleneck ratio of $E$ is defined as
$$
\Phi(E) := \frac{\sum_{x  \in E, y \notin E}{\mu(x)P(x,y)}}{\mu(E)} = \frac{\cE(\1_E)}{\mu(E)}\,,
$$
and the quantity
$$
\Phi_{\ast} := \min_{E \subset \Omega^{V} :\; \mu(E) \leq 1/2 }{\Phi(E)}
$$
is called the expansion (bottleneck ratio) of the chain $P$.

\bc \label{c.Wsr_sr}
Let $\Omega$ be finite  and consider  a sequence of 
 probability spaces $(\Omega^{V_n}, \mu_n)$. There is random sparse reconstruction for such a sequence if and only if there is a sequence of Boolean functions $f_n: \Omega^{V_n} \lora \{-1,1\}$ such that there is full random sparse reconstruction for $f_n$.

\ec
\bpf Suppose there is sparse reconstruction for $(\Omega^{V_n}, \mu_n)$. By Theorem \ref{clueto1} we can find a sequence of random subsets $\mathcal{U}_n$ and real functions $f_n$ such that
 $$\E[\clue_{\mu_n}(f_n \; |\;\mathcal{U}_n)] \to 1.$$ 
 For any fixed $\eps>0$ let us choose $n$ large enough so that $\E[\clue_{\mu_n}(f_n \; |\;\mathcal{U}_n)]\geq 1- {\eps^2}/{8}$.

Let $P_{\mathcal{U}_n}$ be the transition matrix of the heat bath block dynamics generated by $\mathcal{U}_n$. 
Let  $E \subset \Omega^{V_n}$ and $p : = \mu_n(E) \leq 1/2$. Let $\Phi(E)$ be the bottleneck ratio with respect to $P_{\mathcal{U}_n}$.  Then we have
$$
2\, \Phi(E) \ge \frac{\Phi(E)}{1-p} = \frac{\mathcal{E}( \1_{E})}{\Var(\1_{E})} = 1 - \E\left[\clue_{\mu_n}(\1_{E} \; |\;\mathcal{U}_n)\right],
$$
where the second equality used the proof of Lemma \ref{l.eigenclue}. So we have the bound
$$
1 - 2 \Phi(E) \leq \E\left[\clue_{\mu_n}(\1_{E} \; |\;\mathcal{U}_n)\right]\,.
$$
It was proved in \cite{SJ} that, whenever the transition matrix is reversible, we have
$$
\frac{\Phi_{\ast}^2}{2} \leq  \gamma = 1- \lambda_2.
$$
Let $E_{\ast} \subset \Omega^{V_n}$ satisfy $ \mu_n(E_{\ast}) \leq 1/2$ and $\frac{\Phi(E_{\ast})^2}{2} \leq  \gamma$.  Then we have
$$
1-2\sqrt{2(1- \lambda_2) } \leq 1 - 2\Phi(E_{\ast}) \leq  \E[\clue_{\mu_n}(\1_{E_{\ast}} \; |\;\mathcal{U}_n)].
$$
Now, by Lemma~\ref{l.eigenclue}, the maximal achievable clue is $\lambda_2$, hence  $\lambda_2 \geq 1- {\eps^2}/{8}$, by our choice of $n$. Therefore,
$$
\E[\clue_{\mu_n}(\1_{E_{\ast}} \; |\;\mathcal{U}_n)]>1-\eps,
$$
showing full RSR for $f_n = 2 \1_{E_{\ast}}-1$, as  $\clue$ is invariant under linear transformations.
\epf

\subsection{SSM  implies no RSR} \label{ss.SSM}

The aim of this section is to prove Theorem~\ref{t.norSR}, using Lemma~\ref{l.eigenclue}. We start by introducing a condition on spatial correlation decay that ensures rapid mixing of the Glauber block dynamics. We start by introducing some notation.  Let $G(V, E)$ be a (possibly infinite) graph and for a $v \in V$ and $R \in \N$ let $B(v,  R)$ denote the ball of radius $R$ around $v$ with respect to the graph distance. For $T\subset V$ we denote by $\partial T$  the external vertex boundary of $T$, the collection of vertices outside of $T$ which are connected to some $v \in T$. We call a configuration $\psi \in \Omega^{\partial T}$ a boundary condition, and for  $u \in \partial T$ and $j \in \Omega$ we denote by  $\psi^{j}_{u}$ the boundary condition that satisfies   $\psi^{j}_{u} (u) =j$, while at every other vertex,  $\psi^{j}_{u}$ equals to $\psi$.

For the following definition and throughout this section, we restrict ourselves to the case of $\Z_n^d$. 

\bde[Strong Spatial Mixing on cubes] \label{SSM}
Let $\Lambda \subseteq \Z_n^d$ be a $d$-dimensional cube and let $\psi$ be a boundary condition on $\partial \Lambda$. For $u \in \partial \Lambda$, let $\psi_{u}$ be any boundary condition that is equal to $\psi$ everywhere except for $u$. We say that a measure  $\p$ with the spatial Markov property satisfies the Strong Spatial Mixing (SSM) property if there exists $a, b > 0$ such that, for every pair of boundary conditions $\psi$ and $\psi_{u}$ as above and any $B\subseteq \Lambda$, we have 
$$
\|\p_{B}^{\psi} - \p_{B}^{\psi_{u}}\|_{TV} \leq b \exp{(-a\cdot d(u,B))},
$$
where $\p_{B}^{\psi}$ and $\p_{B}^{\psi_{u}}$ denote the measures on $\Omega^B$ induced by the boundary conditions $\psi$ and $\psi_{u}$, respectively, and $d(u,B)$ is the graph distance between a vertex and a subset of vertices.
\ede

The SSM criterion is known to hold for the subcritical (that is, for all  $\beta<\beta_c$) Ising model  on $\Z^d$. This was shown recently in \cite{DSS}for $d\geq 3$, using a novel correlation inequality. More generally, SSM is also known to hold for the Potts model for every $q \ge 2$ for low enough $\beta$ for $\Z^{d}$ with $d>2$, and for the whole subcritical regime in case $d=2$ \cite{MO1}.  

Let $\{ L_n: n \in \N \}$ be a sequence of integers. We shall now define a sequence of random subsets $\mathcal{N}(L_n)$ of $\Z^d$, as in \cite{BCSV}. Denote by $\Lambda_0$ the cube $[0,L_{n}-1]^d \subseteq \Z^d$. Let $e_i \in \Z^d$ be the $i$th unit vector and let $N_0 := \Z^d \setminus \bigcup_{k \in \Z, \; i \in [d]} \left\{\Lambda_0 + k(L_{n}+3) e_i \right\}$. For $w \in \Z^d$ we define  $N_w : = N_0 + w$. It is clear that $N_0$ is invariant under translations with $k(L_{n}+3)x$, where $k \in \Z$ and $x \in \Z^d$. So in fact it is enough to consider translates by the group $\Z^d/(L_{n}+3)\Z^d$, which has $(L_{n}+3)^d$ elements. We can now define the random subset $\mathcal{N}(L_n)$ by picking a uniformly random $w \in \Z^d/(L_{n}+3)\Z^d$ and $N_w$ accordingly. That is, $\p[\mathcal{N}(L_n) =N_w ] = \frac{1}{(L_{n}+3)^d}$.  

The following lemma is a straightforward sharpening of Lemma 3.1 in \cite{BCSV}.

\bl \label{tilerelax}
Let $\p$ be a measure on $\Omega^{\Z^d}$ with the spatial Markov property that satisfies the SSM condition of Definition~\ref{SSM}.  Let $V_n=[0,n]^d$  and $\mathcal{C}_n := \mathcal{N}(L_n) \cap V_n$ and $X^{\mathcal{C}_n}_t$  the corresponding Glauber block chain on $\Omega^{V_n}$, as defined at the beginning of Subsection~\ref{ss.gap}. Then 
$\lambda_{\ast}(P_{\mathcal{C}_n}) \leq C \frac{\log^{d} L_n}{L_n}$.
\el

\bpf
Let $X_t$ and $Y_t$ be two copies of the heat bath block dynamics $X^{\mathcal{C}_n}_t$ on the box $V_n$ that differ at a single vertex $v \in V_n$ at time $t=1$. We will couple the two chains in a way that $\eta_{t}$, the expected number of vertices at which  $X_t$ and $Y_t$ differ after step $t$, tends to decrease.

So let the same $N_w$ be chosen uniformly at random in the two chains. Let $B_w := V_n \setminus N_w$. Note that $B_w$ is the union of disjoint boxes, mostly cubes (except possibly near the boundary of $V_n$).
Now we have three different cases:
\begin{itemize}
    \item $v \in B_w$, in which case   $X_t$ and $Y_t$ are the same on $V_n \cap N_v$, and the configuration on $B_w$ can be updated in a coupled way so that $X_{t+1}= Y_{t+1}$ and $\eta_{t+1}=0$.
    \item $v \in N_w \setminus \partial B_w$. Then at time $t+1$ the disagreement at $v$ remains. Since the two chains agree on $\partial B_w$, the boundary conditions are the same, and again the two chains can be updated in the same way, and thus $\eta_{t+1}=1$.
    \item $v \in \partial B_w$, in which case it is possible that the disagreement from $v$ propagates into $B_w$ because of the different boundary conditions. 
\end{itemize}
We now focus on the third case and bound $\E[\eta_{t+1} \, |\, v \in \partial B_w ]$. Set $R: = C \log L_n $ for a  constant $C$ large enough so that $b \exp{\{-a R\}} \leq \frac{1}{(L_n)^{d}}$ holds.

Let $\Lambda$ be the component of $\B_w$ the boundary of which $v$ lies. Obviously, $v$ might generate disagreement between $X_t$ and $Y_t$ only at vertices of $\Lambda$. Let $F : = \{ u \in \Lambda \, : \, d(u,v) > R \}$. Then by the SSM property we have 
$$
d_{TV}(\mu^{\phi}_{F}, \mu^{\phi^{v}}_{F}) \leq b \exp{\{-a R\}},
$$
where $\mu^{\phi}_{F}$ denotes the measure restricted to $F$ with boundary conditions according to the configuration of $X_t$, while $\mu^{\phi^{v}}_{F}$ is that with boundary conditions according to $Y_t$. The above inequality implies that there is a coupling of configurations $\sigma_{F}$ and $\Tilde{\sigma}_{F}$ distributed according to $\mu^{\phi}_{F}$ and $\mu^{\phi^{v}}_{F}$, respectively, in such a way that
$$
\p[\sigma_{F} \neq \Tilde{\sigma}_{F} ] \leq b \exp{\{-a R\}} \leq \frac{1}{(L_n)^{d}}.   
$$
Using this coupling for updating the two chains we get that
$$
\E[\eta_{t+1} \; |\;v \in \partial B_v ] \leq (2R)^d + \frac{(L_n)^d}{(L_n)^{d}} = C^d\log^{d} L_n + 1.
$$
Now we can conclude, using that $\p[v \in \partial B_w] =  \frac{2d(L_n)^{d-1}}{(L_n+3)^{d}} \leq \frac{2d}{L_n} $ and  $\p[v \in N_w \setminus \partial B_w] \leq \frac{K}{L_n}$ as well, for some absolute constant $K$, that
$$
\E[\eta_{t+1}]\leq \frac{K + 2d (C^d\log^{d} L_n +1)}{L_n} \leq C\frac{\log^{d} L_n}{L_n}, 
$$
for every large enough $n$.

Thus, using that the diameter of $V_n$ is $n$, the path coupling method (see \cite{BD} or \cite[Section 14.2]{LPW}) implies that
$$
\max_{\sigma \in \Omega^{V_n}} {d_{TV}(P_{\mathcal{C}_n}^t(\sigma, \cdot), \p)} \leq n \left(\frac{\log^{d} L_n}{L_n}\right)^{t},
$$
where $P_{\mathcal{C}_n}^t(\sigma, \cdot)$ denotes the distribution of $X^{\mathcal{C}_n}_t$ when started from $\sigma$.
Using that by Corollary 12.7 in \cite{LPW} we have 
$\lambda_{\ast} = \lim_{t \rightarrow \infty}{\left( \max_{\sigma \in \Omega^{V_n}} {d_{TV}(X^{\mathcal{C}_n}_t(\sigma), \p)} \right)^{1/t}}$, the lemma follows.
\epf

\bpf[Proof of Theorem \ref{t.norSR}]
Set $L_n : = \left(\frac{1}{\delta_n}\right)^{1/(d+1)}\log^{d/(d+1)} \left(\frac{1}{\delta_n}\right)$, where $\delta_n\to 0$ is the revealment of $\cU_n$. Consider the random set  $\mathcal{C}_n := \mathcal{N}(L_n) \cap V_n$ as in Lemma~\ref{tilerelax}. Note that $L_n \rightarrow \infty$.

Now consider the following dynamics. At every step $t$ we sample $\mathcal{U}_n$ and $\mathcal{C}_n$ independently, and update the configuration inside all the components  $\Lambda$ of $V_n\setminus \cN(L_n)$ that do not contain any element of $\mathcal{U}_n$. We will call these components ``good cubes''. These dynamics can be also described as a heat bath block chain induced by $\mathcal{H}_n = \left(\mathcal{N}(L_n) \cup \bigcup_{\Lambda_x \cap \mathcal{U}_n \neq \emptyset}{\Lambda_x}\right)$, where the $\Lambda_x$'s are the connected components (cubes) in the complement of $\mathcal{N}(L_n)$.   

With the above choice of $L_n$, for any cube $\Lambda$ outside of $N_v$ we have, by the union bound,
\beq\label{e.badcube}
\p[ \Lambda \cap \mathcal{U}_n \neq \emptyset ] \leq (L_n)^{d} \delta_n  =  \delta_n^{\frac{1}{d+1}} \log^{\frac{d^2}{d+1}} \left(\frac{1}{\delta_n}\right).
\eeq

We use again the path coupling method, and consider two copies of the above dynamics, $X_t$ and $Y_t$  that differ at a single vertex $v \in V_n$. Observe that:
\begin{itemize}
    \item if $v$ falls in a good cube, then we can again couple the two chains to agree from that step onwards;
    \item in case  $v$ is at the outer vertex boundary of a good cube, we can again use the same bound as in Lemma~\ref{tilerelax};
    \item in the remaining cases the number of disagreements remains $1$.
\end{itemize}
If $v \in B_w$, then the probability  that it is not in a good cube is at most~(\ref{e.badcube}), so we get that, for every large $n$,
\begin{align*}
 \E[\eta_{t+1}]&\leq \delta_n^{\frac{1}{d+1}} \log^{\frac{d^2}{d+1}} \left(\frac{1}{\delta_n}\right) +\frac{ 2d (C^d\log^{d} L_n +1)}{L_n} \\
  &\leq K \delta_n^{\frac{1}{d+1}} \log^{\frac{d^2}{d+1}} \left(\frac{1}{\delta_n}\right),    
\end{align*}
for some uniform constant $K<\infty$. So, finishing the argument of Lemma~\ref{tilerelax} for the chain $X^{\mathcal{H}_n}_t$, we get that
$$
\lambda_{\ast}(P_{\mathcal{U}_n}) \leq \delta_n^{\frac{1}{d+1}} \log^{\frac{d^2}{d+1}} \left(\frac{1}{\delta_n}\right).
$$
Now the theorem is immediate from  Lemma~\ref{l.eigenclue}, using that  $\lambda_{2}(P_{\mathcal{U}_n}) \leq \lambda_{\ast}(P_{\mathcal{U}_n})$.
\epf

\begin{rem}
The proof above also gives an explicit quantitative bound on the clue we can get from a random subset with revealment $\delta$. However, Corollary \ref{c.quantnoSR} gives the stronger (general) bound 
$\E[\clue(f_n\,|\,\mathcal{U}_n)]\leq C \sqrt{\delta}$.
\end{rem}

\subsection{ASSM implies no RSR}\label{ss.ASSM}

We shall  need another, stronger notion of spatial mixing, in a more general setting.

\bde[Aggregate Strong Spatial Mixing]\label{ASSM}
Let $G(V,E)$ be a finite graph, $\Omega$ a finite set, and $\p$ a probability measure on $\Omega^V$. For a ball $B(v,R)$ in $G$ and $u \in \partial B(v,R)$, define
$$
a_u(v,R) : = \sup_{\psi,\psi_u \in \Omega^{\partial  B(v,  R)}} \|\p_{v}^{\psi} - \p_{v}^{\psi_{u}}\|_{TV},
$$
where the supremum is over all pairs of  boundary conditions $\psi$ and $\psi_u$ on $\partial B(v,R)$ that differ only at $u$, and where $\p_{v}^{\psi}$ and $\p_{v}^{\psi_{u}}$ are the measures at $v$ induced by the boundary conditions $\psi$ and $\psi_{u}$, respectively.

We say that a measure  $\p$ with the spatial Markov property satisfies the Aggregate Strong Spatial Mixing (ASSM) property  with $R \in \N$ if, for all $v \in V$,
$$
\sum_{u \in \partial B(v,  R)}{a_u}(v,R) \leq \frac{1}{4}.
$$
\ede

We also introduce  the Tree Uniqueness Threshold $\beta_{c}(d)$, defined as
\beq \label{TUT}
\tanh{\beta_{c}(d)} = \frac{1}{d-1}.
\eeq
There is a unique Ising measure on the infinite $d$-regular tree exactly when $\beta \leq \beta_{c}(d)$, hence the name.

\bl[\cite{MS}, Lemma 3]\label{l.ASSMcond}
For all graphs $G$ with maximum degree $d$ and $\beta < \beta_{c}(d)$, there exists $R(\beta,d) \in \N$ such that ASSM holds for the Ising measure $\mu_{G, \beta}$.
\el

For the sake of completeness, we give the definition of censoring, which plays a key role for the results in \cite{BCV} that we are relying on.  This concept expresses that some moves of a (monotone, reversible) Markov chain can be ignored without speeding up the mixing of the chain. Here, monotonicity of a chain on $\{-1,1\}^V$ means the taking one step from a coordinate-wise larger configuration results in a stochastically larger configuration. 

\bde[Censoring] \label{d.cens}
Let $G$ be a graph on $n$ vertices with maximum degree $d$, and let $\mu_{G, \beta}$ be the Ising measure at inverse temperature $\beta$ on $G$. Let $P$ be the transition matrix of an ergodic monotone Markov chain reversible with respect to its stationary measure $\mu_{G, \beta}$. We say that $\{P_{A}\}_{A \subseteq V}$, a collection of monotone stochastic matrices reversible with respect to  $\mu_{G, \beta}$, is a censoring for $P$ if it has the following properties: 
\begin{enumerate}
    \item $P_{A}(\sigma, \Tilde{\sigma}) > 0$ only if $\sigma_{V\setminus A} = \Tilde{\sigma}_{V\setminus A}$;
    \item if $f,g: \{-1, 1\}^{V} \lora \R^{+} $ are increasing (that is, if $\sigma \leq \tau$ coordinate-wise, then $f(\sigma) \leq  f(\tau)$), then, for any $A \subseteq V$, we have  $(Pf, g)_{\mu_{G, \beta}} \leq (P_{A}f, g)_{\mu_{G, \beta}}$ (or equivalently $\mathcal{E}_{P}(f, g) \geq \mathcal{E}_{P_A}(f, g)$).
\end{enumerate}
\ede

For the proof of Theorem~\ref{t.norSRgen}, we will rely on the following result, which is a generalization of Theorem 3 from \cite{MS} to  ergodic monotone Markov chains:

\bth[\cite{BCV}, Theorem 11] \label{t.bcv}
Let $G$ be a graph on $n$ vertices with maximum degree $d$, and let $\mu_{G, \beta}$ the Ising measure at inverse temperature $\beta$ on $G$. Let $P$ be the transition matrix of an ergodic monotone Markov chain reversible with respect to its stationary measure $\mu_{G, \beta}$. 

Suppose $\{P_{A}\}_{A \subseteq V}$ is a censoring for $P$. If ASSM holds with some $R \in \N$, and for all $v \in V$ the mixing time satisfies
$$
T_{\mathrm{mix}}(P_{B(v,R)}) \leq T\,,
$$
then
$$
T_{\mathrm{rel}}(P) = O(T).
$$
\eth
The most common case of a collection $P_{A}$ from Definition~\ref{d.cens}  is when one takes some dynamics with transition matrix $P$ and ignores updates concerning any vertex outside $A$. In our case, for a random subset  $\mathcal{U} \subseteq V$ with corresponding Glauber block dynamics  $X^{\mathcal{U}}_t$ and an $A \subseteq V$, we define a censored dynamics $X^{\mathcal{U}, A}_t$ with transition matrix $P_{\mathcal{U}, A}$ as follows:
\begin{enumerate}
	\item We sample a $U \subseteq V$ according to the distribution of $\mathcal{U}$.
	\item We get $X^{\mathcal{U}, A}_{t+1}$ by resampling the configuration $X^{\mathcal{U}, A}_t$ on  $U^{c} \cap A$  according to the distribution conditioned on $U \cap A$ and independently from the previous values on  $U^{c} \cap A$.
\end{enumerate}

The following lemma is a trivial generalization of Lemmas 14 and 15 in \cite{BCV}. The main difference is that in  \cite{BCV} the authors only consider block dynamics associated with random sets which are distributed uniformly  on a partition of $V$. Its relevance is that it allows us to use Theorem \ref{t.bcv} for Glauber block dynamics.

\bl\label{l.moncen}
Let $\mathcal{U}$ be a random subset of the vertex set $V$ of a graph $G$, independent of everything, with revealment strictly smaller than $1$. Consider the Glauber block dynamics $X^{\mathcal{U}}_t$ with stationary measure  $\mu_{G, \beta}$  associated with $\mathcal{U}$. Then 
\begin{enumerate}
	\item $X^{\mathcal{U}}_t$ is monotone.
	\item For any  $A \subseteq V$, the collection of transition matrices $P_{\mathcal{U}, A}$ as above is a censoring for  $X^{\mathcal{U}}_t$.
\end{enumerate}	
\el

The following theorem is a simple modification of Theorem 3 in \cite{BCV}.

\bth[]\label{t.contrelt}
Let $G(V,E)$ be a finite graph with maximal degree $d$ and let  $\mu_{G, \beta}$ the Ising measure on $G$ with inverse temperature $\beta$ for some $\beta < \beta_{c}(d)$. Fix $0<\delta<1$, and  let  $\mathcal{U}$ be a random subset of $V$, independent of everything, with revealment  at most $\delta$.

Then there exists a constant $K(d,\beta,\delta)<\infty$  such that
$$
T_{\mathrm{rel}}(P_{\mathcal{U}}) \leq K(d, \beta, \delta).
$$
\eth

\bpf
In light of Lemma \ref{l.moncen} and Lemma \ref{l.ASSMcond}, the conditions of Theorem \ref{t.bcv} will be satisfied if we can show rapid mixing of the censored heat bath block chain  $X_t^{\mathcal{U} , B(v,R)}$ for all $v \in V$, with $R = R(d,\beta)$ given by Lemma~\ref{l.ASSMcond}. That is, at each step of the dynamics we only update the spin of vertices which are outside of $\mathcal{U}$ and are in $B(v,R)$.

In order to bound the mixing time, we shall use the standard monotone grand coupling of the single site Glauber dynamics, with the difference that in each step we update all the vertices in $\mathcal{U} ^{c} \cap B(v,R)$ in an arbitrary order. For the monotone grand coupling, see, for example, \cite[Example 22.2]{LPW}.

We first bound the probability that all $w \in B(v,R)$ are updated by some time $T$. The probability that at a step a particular vertex is not updated is at most $\delta$, and therefore by time $T$ the probability that there exists  $w \in B(v,R)$ which is not updated is at most $|B(v,R)| \, \delta^{T}$, using the union bound. This means that, if
\beq\label{e.largeT}
T \geq \frac{\log{(2|B(v,R)|)} }{\log(\delta^{-1}) }, 
\eeq
then at time $T$ all vertices in $B(v,R)$ have already been updated at least once with probability at least $1/2$.

Now observe that each time a $w \in B(v,R)$ is updated, with a uniform positive probability its spin is set to the same value for all configurations in the grand coupling. Indeed, irrespective of the spins of the neighbouring vertices, the spin at $w$ will be $-1$ with probability at least 
$$
\frac{e^{-\beta d}}{e^{\beta d}+e^{-\beta d} }, 
$$
and will be $+1$ with at least the same probability. Altogether, we get the same spin in all configurations with probability at least twice this, which is larger than $e^{-2\beta d}$.

Suppose that by time $T$ every $w \in B(v,R)$ is updated at least once. If, at the time of the last update before $T$, the above coupling event does happen, for every $w$, then all the spins will agree at time $T$. This happens with conditional probability at least $e^{-2\beta d |B(v,R)|}$.

Altogether, if $T$ satisfies~(\ref{e.largeT}), then with probability at least $\frac12 e^{-2\beta d |B(v,R)|}$ all spins agree at time $T$. Repeating this some $O_{d,\beta,R}(1)$ times, the probability of agreement will be close to 1, hence
$$
T_{\mathrm{mix}}(P_{\mathcal{U}, B(v,R) }) = O_{d,\beta,R}(T) = O_{d,\beta,R}\left(\frac{1 }{\log(\delta^{-1}) }\right).
$$
Therefore, by Theorem~\ref{t.bcv}, 
$$
T_{\mathrm{rel}}(P) = O_{d,\beta,R}\left(\frac{1 }{\log(\delta^{-1}) }\right),
$$
proving Theorem~\ref{t.contrelt}.\epf

\bpf[Proof of Theorem \ref{t.norSRgen}]
Aiming for contradiction, suppose that there is random sparse reconstruction for  $\{ \mu_n \}_{n \in \N}$. Then by Theorem~\ref{clueto1} there exists a sequence of random subsets $\mathcal{U}_n \subseteq V_n$ with $\delta(\mathcal{U}_n) \to 0$ and a respective sequence $\{f_n: (\{-1,1\}^{V_n},\,\mu_n) \longrightarrow \R \}_{n \in \N}$ such that
$$
 \E[\clue(f_n\,|\,\mathcal{U}_n)]  \to 1.
$$
As before, we denote by $P_{\mathcal{U}_n}$ the transition matrix of the block dynamics generated by the random subset   $\mathcal{U}_n$. As the revealment tends to $0$, for every $\eta>0$ there is an $N \in \N$ such that $\delta(\mathcal{U}_n)\leq \eta$ for every $n\geq N$. 

Thus, by Theorem \ref{t.contrelt}, there is a constant $ K(d,\beta,\eta)$ such that for every $n\geq N$ 
$$
T_{\mathrm{rel}}(P_{\mathcal{U}_n}) \leq  K(d,\beta, \eta).
$$
This implies, using Lemma~\ref{l.eigenclue}, that for any function $\{g_n: \{-1,1\}^{V_n} \longrightarrow \R \}_{n \in \N}$ we have
$$
 \E[\clue(g_n\,|\,\mathcal{U}_n)] = \lambda_2(P_{\mathcal{U}_n}) \leq \lambda_{\ast}(P_{\mathcal{U}_n}) < c(d,\beta, \eta),
$$
for some constant $c(d,\beta,\delta)<1$, whenever $n \geq N$. This contradicts the fact that there is a sequence of functions  $f_n$ with $ \E[\clue(f_n\,|\,\mathcal{U}_n)]   \to 1$.
\epf

\section{Divide-and-Color measures}\label{s.DaC}

\subsection{Definition}
We introduce a class of measures that can be described via a random graph on the set vertices, a subclass of  the Generalised Divide and Color measures, as defined and investigated in \cite{StT}. The reason for studying this subclass is that the Fourier-Walsh transform and some of its nice properties extend to these measures, while, most importantly, the Ising measure on any graph can be described in this way. 

Let $G(V, E)$ be a graph.  We are going to define a measure on the usual configuration space $ \left\{ -1,1 \right\}^V$. Let $N  \subseteq E$, and consider   the event $A_{N} = \big\{  \sigma_{u} \sigma_{v} = 1 \text{ for all }(u,v) \in N \big\}$.  The measure $\p_{N}$ is defined as the uniform measure conditioned on $A_{N}$. In other words, we consider the uniform measure on configurations which are constant on the connectivity clusters of $G(V, N)$, the  subgraph of $G$ spanned by the edge set $N$. 

In this section, we consider probability measures on $ \left\{ -1,1 \right\}^V$ that are induced by a random subset of the edge set $\mathcal{N}$ as follows.
\bde [Generalised Divide and Color measure \cite{StT}]
Let $G(V, E)$ be a graph and $\mathcal{N} \subseteq E$ be an edge percolation on $E$  (a random subset of the edge set $E$). Then $\mathcal{N}$ induces a probability measure on $ \left\{ -1,1 \right\}^V$ by
$$
\p[\sigma = x] := \sum_{N \subseteq E}{\p_\nu[\mathcal{N} = N]\, \p_{N}[\sigma = x]}  = \E_\nu \left[\p_{\mathcal{N}}[\sigma= x ] \right],
$$
where $\nu$ denotes the probability measure associated with $\mathcal{N}$.
\ede


We may think of $\p$ in the following way: take $N  \subseteq E$, then flip a fair coin independently for each cluster of $N$; that is, we assign this value to every vertex of the respective cluster. It is easy to verify that this  defines the uniform measure on $A_{N}$. In case $\mathcal{N} \subseteq E$ is random, we first sample according to the percolation $\mathcal{N}$, then assign spin values to clusters as above. The most important special case of this is, of course, the FK--representation of the Ising measure. 

Generalised Divide and Color measures, as defined in  \cite{StT}, are more general in that they allow the measure $\p_{N}$ to be a $p$-biased product measure rather than just the uniform measure $p=1/2$. However, we do not see how to extend our arguments to this more general case.

\subsection{Generalized Fourier-Walsh transform for DaC measures }

We now introduce a generalisation of the Fourier-Walsh transform for this family of measures.  First,  for $N  \subseteq E$, we introduce the notation 
$$
\langle {N} \rangle : = \{ S \subseteq V \; :\; \E [\chi_{S} | A_{N}] = 1 \}.
$$

It is straightforward to verify that these are subsets which can be written as a union of even subsets of the connectivity clusters of $N$. Equivalently, if one identifies subsets of $V$ with elements of  $\mathbb{F}_2^{V}$, it is the linear subspace of $\mathbb{F}_2^{V}$ generated by $N$. 

Let $f(\omega)=\sum_{T\subset V} \widehat{f}(T)\chi_T(\omega)$. As, for every $S\subseteq V$ and $L \in \langle {N} \rangle$,
$$
\p_N[\chi_S(\omega) =\chi_{S\triangle L}(\omega)] =1,
$$
 the subsets of $V$ can be divided into $\langle N\rangle$-equivalence classes of the form
 $$
 [S]_{ N} := \{ S \triangle L : L \in \langle N \rangle \},
 $$
 so that, for any $\omega \in A_{N}$, the value of $\chi_S(\omega)$ only depends on the  class of $S$. Note that the size of each equivalence class is $|\langle N \rangle|$. Moreover, if we let $\mathcal{C}^{N}_j$ for $j =1, \dots,  k(N)$ denote the connectivity clusters according to $N$, then an equivalence class $[S]_N$ is determined by the parity of $|S \cap \mathcal{C}^{N}_j|$ for every $j \in [k(N)]$. Indeed, it is exactly the collection of all subsets that share the   parity profile with $S$. 

We define, for  any $S \subseteq V$, the Fourier-Walsh coefficient of an $\langle N\rangle$-equivalence class as
\beq \label{Fourier_congr}
\widehat{f}([S]_N): = \sum_{T \in  [S]_{N}}{ \widehat{f}(T)}\,,
\eeq
and then we can write
\beq  \label{eq.DaCFourier_congr}
f(\omega)=\frac{1}{|\langle N \rangle|}\sum_{T\subseteq V} \widehat{f}([T]_N)\chi_T(\omega)\,, \quad \text{ for }\omega \in A_{N}.
\eeq
Using the  $\mathbb{F}_2^{V}$ subspace interpretation, it is immediate that $|\langle N \rangle|= 2^{\dim \langle N \rangle}$. At the same time, every cluster $C$ of $N$ can be generated by $|C|-1$ edges, and subspaces generated by edges from different connectivity clusters are linearly independent over $\mathbb{F}_2$. Thus we get that
\beq \label{eq.SizeN}
|\langle N \rangle|= 2^{|V|-k(N)},
\eeq
where $k(N)$ denotes the number of connectivity clusters of $N$. One can of course get this formula from direct counting, as well, without the $\mathbb{F}_2^{V}$  interpretation.






By definition, we have $\E[\chi_T|A_{N}] = 1$ if and only if $T \in \langle N \rangle$; otherwise, $\E[\chi_T|A_{N}] = 0$. Therefore, for any $N \subseteq E$, we have 
\beq \label{eq. Dacprod}
\E_N[\chi_T \chi_S] = \E[\chi_{T\triangle S} |A_{N}] = \left\{
	\begin{array}{ll}
1	 & \mbox{if } T \triangle S \in \langle N \rangle \\
		0 & \mbox{if } T \triangle S \notin \langle N \rangle.
	\end{array}
\right.
\eeq
Hence, for any two functions $f,g :\left\{ -1,1 \right\}^n \lora \R$, 
\beq \label{Cov_cond}
\begin{aligned}
 \E[f g| A_{N}] &= \frac{1}{2^{|V|-k(N)}}\frac{1}{2^{|V|-k(N)}}\sum_{\substack{T,S\subseteq V \\ T \triangle S \in \langle N \rangle}}{\widehat{f}([T]_N)\,\widehat{g}([S]_N)} \\
 &= \frac{1}{2^{|V|-k(N)}}\sum_{T \subseteq V}{\widehat{f}([T]_N)\,\widehat{g}([T]_N)},   
\end{aligned} 
\eeq
where we used the fact that $[T]_N = [S]_N$ whenever $T \triangle S \in \langle N \rangle$.
In particular,
$$
\E[f^2| A_{N}]  = \frac{1}{|\langle N \rangle|}\sum_{T\subseteq V}{\widehat{f}([T]_N)^2} ,
$$
and
\beq
\label{eq.fourempty}
\E[f| A_{N}] = \frac{1}{|\langle N \rangle|}\sum_{T\in \langle N \rangle } \widehat{f}([T]_N) = \widehat{f}([\emptyset]_N).
\eeq
We can now express the first and the second moment of functions with respect to a Divide and Color measure: 
%

\bpr \label{cov_DaC} 
Let $\mathcal{N}$ be a random edge set with distribution $\nu$, and let $\p$ be the corresponding Divide-and-Color measure measure. Then, for any two functions $f,g: \left\{ -1,1 \right\}^n \lora \mathbb{R}$,
\beq \label{eq.Dacexp}
\E[f] =\E_{\nu}\left[\widehat{f}([\emptyset]_\mathcal{N})\right]
\eeq
\beq \label{eq.Dacquadexp}
 \E[fg] =   \frac{1}{2^{|V|}}\sum_{T \subseteq V}{\E_\nu \left[2^{k(\mathcal{N})} \widehat{f}([T]_\mathcal{N})\,\widehat{g}([T]_\mathcal{N})\right]}. 
\eeq
\epr

\bpf
The first statement is immediate from~\eqref{eq.fourempty}, the second statement from~\eqref{Cov_cond}.
%
\epf


\comm{\subsection{Example of FK-Ising}

It is worth having a closer look at our most important example, the FK Ising model. Recall that the Ising measure can be expressed as a DaC measure where the respective percolation measure is the random cluster model with $q=2$, that is 
$$
\phi_{p,2}(\omega) = \frac{1}{Z}\prod_{e \in E}{p^{\omega(e)}(1-p)^{1-\omega(e)}2^{k(\omega)}}
$$
where $Z$ is a constant that makes $\phi_{p,2}$ a probability measure. Now let $\sigma$ be a DaC measure induced by $\omega$, that is $\sigma$ is distributed according to the Ising measure with parameter $\beta$. Now let $p=e^{-2 \beta J }$ and we denote by  $\mathcal{B}$  a Bernoulli($p$) edge percolation on $G(V,E)$. that is, each edge $e \in E$ is in $\mathcal{B}$ with probability $p$ and not in $\mathcal{B}$ with probability $1-p$, independently from any other edge. So we can write
\begin{align*}
 \E [f(\sigma)g(\sigma)] = &\frac{1}{2^{|V|}} \sum_{T \subseteq V}{\E[2^{k(\omega)} \widehat{f}_{\omega}(T)\widehat{g}_{\omega}(T)]} \\
 =& \frac{1}{2^{|V|}}\frac{1}{Z}\sum_{T \subseteq V}{\E [2^{2k(\mathcal{B})} \widehat{f}_{\mathcal{B}}(T)\widehat{g}_{\mathcal{B}}(T)]}\\
 =& \frac{2^{|V|}}{Z}\sum_{T \subseteq V}{\E [ \widehat{\overline{f}}_{\mathcal{B}}(T)\widehat{\overline{g}}_{\mathcal{B}}(T)]}. 
\end{align*}}

\begin{comm}{In such a way  one can express scalar products and second moments for the Ising measure via a Bernoulli percolation. In particular, when calculating $\clue$ the factor $\frac{2^{|V|}}{Z}$ cancels out-

The most important property of this is not so much that now the expectation is according to a product measure but the following fact. Let $\mathcal{N}$ be a linear subspace and $L \subseteq [n]$ arbitrary. Let $\Tilde{\mathcal{N}} =  \langle \mathcal{N}, L \rangle$. Then, obviously, $\Tilde{\mathcal{N}} = \mathcal{N} \cup (\mathcal{N} \oplus L) $. Now it is easy to see that for any  $T \subseteq [n]$, we have
$$
\widehat{\overline{f}}_{\Tilde{\mathcal{N}}}(T) = \frac{1}{2}\left(\widehat{\overline{f}}_{\mathcal{N}}(T) + \widehat{\overline{f}}_{\mathcal{N}}(T\oplus L) \right).
$$
Indeed, in case $L \notin \mathcal{N}$, then $\widehat{\overline{f}}_{\Tilde{\mathcal{N}}}(T) = \widehat{\overline{f}}_{\mathcal{N}}(T) =  \widehat{\overline{f}}_{\mathcal{N}}(T\oplus L)$, and the equality holds trivially. If  $L \in \mathcal{N}$,we only need to observe that $\dim \Tilde{\mathcal{N}} = \dim \mathcal{N} +1$ and that $\widehat{f}_{\Tilde{\mathcal{N}}}(T) = \widehat{f}_{\mathcal{N}}(T) + \widehat{f}_{\mathcal{N}}(T \oplus L)$. An important consequence is that
$$
\widehat{f}^2_{\Tilde{\mathcal{N}}}(T) + \widehat{f}^2_{\Tilde{\mathcal{N}}}(T \oplus L) \leq \widehat{f}^2_{\mathcal{N}}(T) + \widehat{f}^2_{\mathcal{N}}(T \oplus L),
$$
for any $L \in\Tilde{\mathcal{N}}$ and $T \subseteq [n]$.
Indeed, 
\begin{align*}
\widehat{f}^2_{\Tilde{\mathcal{N}}}(T) + \widehat{f}^2_{\Tilde{\mathcal{N}}}(T \oplus L) = 2 \widehat{f}^2_{\Tilde{\mathcal{N}}}(T) = \frac{1}{2}\left(\widehat{\overline{f}}_{\mathcal{N}}(T) + \widehat{\overline{f}}_{\mathcal{N}}(T\oplus L) \right)^2 &\leq
\\
\frac{1}{2}\left(\widehat{\overline{f}}_{\mathcal{N}}(T) + \widehat{\overline{f}}_{\mathcal{N}}(T\oplus L) \right)^2 + \frac{1}{2}\left(\widehat{\overline{f}}_{\mathcal{N}}(T) - \widehat{\overline{f}}_{\mathcal{N}}(T\oplus L) \right)^2 &,
\end{align*}
which, in particular, implies that $\E[f^2|\mathcal{N}]$ is monotone decreasing in $\mathcal{N}$.

In what way these observations help us? Our aim is to show that 
$$
\sum_{T \subseteq [n^{d}]} {\E\left[\beta_{N}(T) \widehat{\overline{f}}^2_{\mathcal{N}}(T)\right]} \leq C \E\left[\beta_{N}(X) \right]
$$}
\end{comm}

\subsection{A clue bound for DaC measures}

In contrast with the usual Fourier-Walsh transform, see~\eqref{VarFourierCond}, here in general there seems to be no efficient way to express $ \Var(\E[f \,|\, \mathcal{F}_{U}])$. Yet, for any fixed $N \subseteq E$, we have
\beq \label{eq.DaCFour_condcongr}
\E_N[ \chi_S\,|\, \mathcal{F}_{U}] = \left\{
	\begin{array}{ll}
\chi_S & \mbox{if exists } T \subseteq U \text{ such that } T   \in [S]_{N}\,, \\
		0 & \mbox{otherwise.}
	\end{array}
\right.
\eeq
To see this, first note that $\E_N[ \chi_T\,|\, \mathcal{F}_{U}]$ is constant for $T\in [S]_{N}$, because, for any $L \in \genN$, we have 
$$
\E_N[ \chi_S\,|\, \mathcal{F}_{U}] = \E_N[\chi_L \chi_S\,|\, \mathcal{F}_{U}]  = \E_N[\chi_{S \triangle L}\,|\, \mathcal{F}_{U}].
$$
The first line of \eqref{eq.DaCFour_condcongr} follows immediately. At the same time, in case for every $T \subseteq U$  we have $T \notin [S]_{N}$, then by~\eqref{eq. Dacprod} we also have $\E_N[\chi_T \chi_S] =0$, and, since $\sigma(\{\chi_T \,:\, T \subseteq U\}) = \mathcal{F}_{U}$, this proves the second line of~\eqref{eq.DaCFour_condcongr}. 

Thus, as in \eqref{eq.DaCFourier_congr}, we have  
$$
\E_N[f \,|\, \mathcal{F}_{U}  \}](\sigma) =\frac{2^{k(N)}}{2^{|V|}} \sum_{T : [T]_{N} {\subseteq} U} \widehat{{f}} ([T]_N)\chi_T(\sigma)\,,
$$
where we write $ [T]_{N} \subseteq U$  to mean that there exists an $S \in  [T]_{\langle N\rangle}$ such that $S \subseteq U$. Moreover, by~\eqref{Cov_cond},
$$
\E_N\left[\E^2_N\left[f \,|\, \mathcal{F}_{U}\right]\right] = \frac{2^{k(N)}}{2^{|V|}} \sum_{T : [T]_{N} {\subseteq} U} \widehat{f} ([T]_N)^2\,.
$$
Consequently, for a DaC measure, we have
$$
\E\big[\E_\cN^2\left[f \,|\,   \mathcal{F}_{U}\right] \,\big|\,\mathcal{N} \big] =\frac{2^{k(\cN)}}{2^{|V|}} \sum_{T : [T]_{\mathcal{N} \subseteq U}} \widehat{f}([T]_\mathcal{N})^2.
$$
Taking expectation on both sides w.r.t.~$\cN$, and assuming that $\E[f]= 0$, so that $\E[f^2]=\Var(f)$, give 
\beq\label{e.VarkN}
\begin{aligned}
\Var(\E[f \,|\, \mathcal{F}_{U}]) 
\leq
\Var\big(\E[f \,|\,  \mathcal{F}_{U}  \vee \sigma(\mathcal{N})]\big)
&=\E\big[\E^2[f \,|\,   \mathcal{F}_{U}  \vee \sigma(\mathcal{N})] \big] \\
&=  \E_\nu\left[\frac{2^{k(\mathcal{N})}}{2^{|V|}}\sum_{ [T]_{\mathcal{N}} \subseteq U } \widehat{f}([T]_\mathcal{N})^2\right]\,,  
\end{aligned}
\eeq
where the inequality holds because we take the $L^2$-norm of a projection to a larger $\sigma$-algebra. 

Similarly to Definition~\ref{d.Spec} that we had for the uniform and other product measures, we shall introduce a notion of spectral sample $\Spec^{\mathcal{N}}_{f}$ for DaC measures. In our setting, the value of the spectral sample is not a single subset of coordinates, but an equivalence class of a subset $T \subseteq V$ with respect to some $N \subseteq E(G)$. Namely, given the bond percolation configuration $\cN$ sampled from $\nu$,  we let, for any $T\subseteq V$,
\beq\label{e.SpecN}
\p\big[ \Spec^{\mathcal{N}}_{f} = [T]_{N} \,\big|\, \cN=N \big] := \frac{\widehat{f}([T]_{N})^2}{\E[f^2]}. 
\eeq


With this definition, using that $\Spec^\cN_f$ has $\gencN$ elements, we can rewrite~\eqref{e.VarkN} as follows, which can be thought of as a generalised version of~\eqref{clueFourier}. For every $f \in L^2( \{-1,1\}^{V}, \p)$ with $\E[f] =0$: 
$$
\clue(f\,|\,U) \leq \p\left[ \exists \, S \in \Spec^{\mathcal{N}}_{f} \,:\,S\subseteq U \right] = \p\left[  \Spec^{\mathcal{N}}_{f} \subseteq U \right].
$$ 
It is straightforward to verify that in case we substitute $U \subseteq V$ with a random subset $\mathcal{U}$, independent of the spin system, just the same way we get:
\beq \label{DaCclueFourier}
 \E\left[\clue(f\,|\,\mathcal{U})\right] \leq \p\left[ \exists \, S \in \Spec^{\mathcal{N}}_{f} \,:\,S\subseteq \mathcal{U} \right] = \p\left[  \Spec^{\mathcal{N}}_{f} \subseteq \mathcal{U} \right].
\eeq

As a first application of~\eqref{DaCclueFourier}, we bound the average clue of the magnetization in a Divide and Color model. In the sequel, for  one-element subsets, with a slight abuse of notation, we are going to write $[u]_N$ instead of $[\{u\}]_N$.

\bpr \label{DacClueMag} 
Let $\p$ be a  Divide and Color measure generated by the percolation $\mathcal{N}$, and let $M= \frac{1}{|V|}\sum_{v \in V} {\sigma_v}$.   For any random subset $\mathcal{U} \subseteq V$, independent of $\sigma$, with revealment $\delta$,
$$
 \E\big[\clue(M\,|\,\cU)\big] \leq \delta \, \frac{\sum_{v \in V}{\E_{\nu} |\mathcal{C}_{v}^{\mathcal{N}}|^2}}{\sum_{v \in V}{\E_{\nu} |\mathcal{C}_{v}^{\mathcal{N}}|}}, 
$$
where $|\mathcal{C}_{v}^{{N}}|$ denotes the cluster size of  $v \in V$ with respect to ${N}$.
\epr
\bpf
Denoting the spectral sample of the average by $\Spec^{\mathcal{N}}_{M}$, we have the following estimate:
\beq\label{e.sumu}
\p\big[  \Spec^{\mathcal{N}}_{M} \subseteq \mathcal{U} \big]
 =\p\big[ \exists u \in \mathcal{U} ,\; \{u\} \in \Spec^{\mathcal{N}}_{M} \big] 
 \leq \sum_{u \in V}{\p\left[   \Spec^{\mathcal{N}}_{M} = [u]_{\mathcal{N}} \right] \, \p\left[   u \in \mathcal{U} \right]}. 
\eeq

We can express the variance in terms of the cluster sizes:
$$
\Var(M) = \frac{1}{|V|^2}\sum_{u \in V}  \sum_{v \in V}{\Cov( \sigma_{u}, \sigma_{v})} = \frac{1}{|V|^2}\sum_{u \in V}  \sum_{v \in V}{\p[v \in \mathcal{C}^{\mathcal{N}}_{u}]} = \frac{1}{|V|^2}\sum_{u \in V} \E_{\nu} |\mathcal{C}_{u}^{\mathcal{N}}|. 
$$
Thus, using that $\widehat{M}( [u]_\mathcal{N}) = |\mathcal{C}_{u}^{\mathcal{N}}|/|V|$,  we obtain
$$
\p\left[   \Spec^{\mathcal{N}}_{M} = [u]_{\mathcal{N}} \right] = \sum_{N \subseteq E} \p_\nu[\mathcal{N} = N] \frac{\widehat{M}^2( [u]_\mathcal{N}) }{\Var(M)}=\frac{\E[\widehat{M}^2( [u]_\mathcal{N})}{\Var(M)} = \frac{\E_{\nu}\big[ |\mathcal{C}_{u}^{\mathcal{N}}|^2/|V|^2\big]}{\sum_{u \in V} \E_{\nu} \big[ |\mathcal{C}_{u}^{\mathcal{N}}|/|V|^2\big]}.
$$
Using~\eqref{e.sumu} and~\eqref{DaCclueFourier}, we obtain the proposition.\epf

Now we use \eqref{DaCclueFourier} to  bound the expected clue of a general  function. We introduce the notation $\clue(f \,|\, \mathcal{N}) = \frac{\Var(\E[f \,|\, \mathcal{N}])}{\Var(f)}$ as a natural extension of $\clue$.

\bth \label{DaCClueThm}
Let  $\p$ be a Divide and Color measure on $\{-1,1\}^{V}$. Let $f: \{-1,1\}^{V} \longrightarrow \mathbb{R}$ be any function satisfying $\E[f] = 0$, and $\mathcal{U}$ a random subset of $V$, independent of the spin values. Then
\beq
\E \left[\clue(f\,|\,\mathcal{U})\right] \leq \delta_{\mathcal{U}} \, \E_{\nu} \left[\max_{v \in V }{|\mathcal{C}^{\mathcal{N}}_{v}|}\right]  +  \clue(f\,|\,\mathcal{N}),
\eeq
 where $\mathcal{C}_{v}^{\mathcal{N}}$ denotes the connectivity cluster  of $v \in V$ with respect to the percolation $\mathcal{N}$. 
\eth 

\bpf 
First recall that $\E[f \,|\, \mathcal{N}] = \widehat{f}([\emptyset]_\mathcal{N})$ by \eqref{eq.fourempty}, and hence
$$
\clue(f\,|\,\mathcal{N})
 = \frac{\Var(\E[f \,|\, \mathcal{N}])}{\Var(f)}=  \frac{\E_{\nu}\left[\widehat{f}([\emptyset]_\mathcal{N})^2\right]}{\Var(f)}=\p\left[\Spec^{\mathcal{N}}_{f}   =[\emptyset]_\mathcal{N}\right].
$$ 
Therefore
\beq\label{splitspecclue} 
\p\left[  \Spec^{\mathcal{N}}_{f} \subseteq \mathcal{U}\right] = 
\p\left[  \Spec^{\mathcal{N}}_{f} \subseteq \mathcal{U},\; \; \Spec^{\mathcal{N}}_{f} \neq [\emptyset]_\mathcal{N}   \right] + \clue(f\,|\,\mathcal{N}).
\eeq
 We imitate the argument of \cite[Theorem 2.8]{GaPe} to bound the second term.  We first sample $\Spec^{\mathcal{N}}_{f}$ conditioned on $\{ \Spec^{\mathcal{N}}_{f} \neq [\emptyset]_\mathcal{N} \}$.
For a congruence class $ [T]_{ N}$, we define the set of coordinates
$$
\overline{[T]}_{ N}  := \{ v \in V\;: \; [v]_{N}{\subseteq} T  \}.
$$
Note that $\overline{[T]}_{ N}$ is the union of connectivity clusters $\mathcal{C}^{\mathcal{N}}_v$ such that  $T \cap \mathcal{C}^{\mathcal{N}}_v$ is odd.

Conditioned on $\Spec^{\mathcal{N}}_{f} \neq[\emptyset]_\mathcal{N}$, the set $\overline{\Spec_{f}}^{\mathcal{N}}
 $ is non-empty, and we can take a uniform random element of $\overline{\Spec_{f}}^{\mathcal{N}}
 $, denoted by $X$. For sake of completeness we define $X$ as a a uniform random element of $V$, in case $\overline{\Spec_{f}}^{\mathcal{N}}
 = \emptyset$.  Furthermore, if $\Spec^{\mathcal{N}}_{f} \subseteq \mathcal{U}$ holds, $X$ has to be in the same cluster as some $u \in \mathcal{U}$, or equivalently, $[X]_\mathcal{N}=[u]_\mathcal{N}$ must hold. So we can write:
\beq
\begin{aligned} \label{DaCspecest}
\p\left[  \Spec^{\mathcal{N}}_{f}  \subseteq \mathcal{U},\; \; \Spec^{\mathcal{N}}_{f} \neq [\emptyset]_\mathcal{N}   \right] 
&\leq \p\left[ \exists u \in \mathcal{U} \; \text{such that}\;  [X]_\mathcal{N}=[u]_\mathcal{N} \right]\\ 
 &\leq \sum _{u \in V}{\p\big[ [X]_\mathcal{N}=[u]_\mathcal{N}, \ u \in \mathcal{U}\big] }\\
&=\sum _{u \in V}{\p\big[ [X]_\mathcal{N}=[u]_\mathcal{N}\big]\p\left[ u \in \mathcal{U} \right]}\\
&\leq \delta_{\mathcal{U}}\sum _{u \in V}{\p\big[[X]_\mathcal{N}=[u]_\mathcal{N}\big]}.
\end{aligned}
\eeq
As for the sum of probabilities, we can write it in a more compact form as
$$
\sum _{u \in V}{\p\big[[X]_\mathcal{N}=[u]_\mathcal{N}\big]} = \sum _{u \in V}{\p\left[ u \in \mathcal{C}_{X}^{\mathcal{N}}\right]} = \E\left[|\mathcal{C}_{X}^{\mathcal{N}}|\right],
$$
where $\mathcal{C}_{X}^{\mathcal{N}}$ is the connectivity cluster of $X$ in $\mathcal{N}$.
Trivially, $\E|\mathcal{C}_{X}^{\mathcal{N}}|\leq \E \left[\max_{v \in V}{|\mathcal{C}^{\mathcal{N}}_v|}\right]$. Writing this back into~\eqref{splitspecclue}, we get the stated bound.
\epf

\bc \label{DacClueSym}
Besides the conditions of Theorem \ref{DaCClueThm} assume that $f: \{-1,1\}^{V} \longrightarrow \R$ is odd, that is $f(-\omega) = -f(\omega)$. Then
$$
\E\left[\clue(f\,|\,\cU)\right]\leq \delta_{\mathcal{U}} \, \E \left[\max_{v \in V}{|\mathcal{C}^{\mathcal{N}}_v|}\right]
$$
\ec

\bpf
It is a direct consequence of Theorem \ref{DaCClueThm}. We only have to note 
that in this case $\E[f \,|\, \mathcal{N}]  \equiv 0$ and thus $ \Var(\E[f \,|\, \mathcal{N}]) =0$ as well.
\epf

\bpf[Proof of Proposition~\ref{p.FK}] Using the FK random cluster representation of the Ising model as our Divide and Colour measure, (\ref{e.FKmagn}) follows directly from Proposition~\ref{DacClueMag}, while~(\ref{e.FKodd}) follows from Corollary~\ref{DacClueSym}.
\epf

\begin{rem}
It is worth pointing out that, in the proof of  Theorem \ref{DaCClueThm}, the definition of the random vertex $X$ was arbitrary and sub-optimal. First, observe that $X$ being uniform on the vertices of $\overline{\Spec_{f}}^{\mathcal{N}}
 $ is arbitrary, as the proof would work in the same way using \emph{any} distribution supported on 
$\overline{\Spec_{f}}^{\mathcal{N}}$. As we want to make  $\E|\mathcal{C}_{X}^{\mathcal{N}}|$ as small as possible, the best way to sample $X$ is to choose it uniformly from the smallest connectivity cluster of $\overline{\Spec_{f}}^{\mathcal{N}}$.

We can in fact express $\E|\mathcal{C}_{X}^{\mathcal{N}}|$ in terms of the generalised Fourier coefficients.  In case  $\Spec^{\mathcal{N}}_{f} = [T]_N$,  we sample $X$ from the smallest cluster of ${N}$ that has an odd intersection with $T$. Let us define, for $ [T]_N \neq  [\emptyset]_N$,
$$
\beta([T]_N) := \min_{v \in \overline{[T]}_N} \{ |\mathcal{C}_{v}^{N}|  \}.
$$ 
For $ [T]_N =  [\emptyset]_N$, set $\beta([T]_N) = 0$.

Now we can write
$$
\E\left[|\mathcal{C}_{X}^{\mathcal{N}}| \,\bigm|\,   \Spec^{N}_{f} = [T]_N \right] =  \beta([T]_N) ,
$$ 
and therefore for any non-constant function $f$, we have
\beq \label{betaeq}
\E|\mathcal{C}_{X}^{\mathcal{N}}|=   \frac{ \sum_{T \subseteq V } {\E\left[ 2^{k(\mathcal{N})} \widehat{f}([T]_{\mathcal{N}})^2 \beta([T]_\mathcal{N})  \right]}}{\sum_{T \subseteq V } {\E\left[ 2^{k(\mathcal{N})} \widehat{f}([T]_{\mathcal{N}})^2\right]}}.
\eeq
In the special case of the average magnetization $M= \frac{1}{|V|}\sum_{v \in V} {\sigma_v}$, the only non-zero Fourier coefficients are $\widehat{M}( [u]_\mathcal{N}) = |\mathcal{C}_{u}^{\mathcal{N}}|/|V|$, each cluster $\cC_j^\cN$ of $\cN$ appearing for $|\gencN|$ many subsets $T\subseteq V$; note here that $[u]_\cN=[v]_\cN$ for any $u,v \in \cC_j^\cN$. Furthermore, $\beta([v]_{\mathcal{N}}) = |\mathcal{C}_{v}^{\mathcal{N}}|$. Therefore,~\eqref{betaeq} translates to 
$$
\E|\mathcal{C}_{X}^{\mathcal{N}}| = \frac{ 
\E\left[\sum_{j=1}^{k(\cN)} |\cC_j^\cN|^3 \right]}
{\E\left[\sum_{j=1}^{k(\cN)} |\cC_j^\cN|^2 \right]}\,,
$$
which then exactly gives the statement of Proposition~\ref{DacClueMag}. 

In general, the heuristics are that larger subsets $T \subseteq V$ should typically have smaller $\beta([T]_{\mathcal{N}})$, since  a large $T$ seems likely to contain a vertex having  odd intersection with a small connectivity cluster. Nevertheless, $\beta([T]_{\mathcal{N}})$ and $ \widehat{f}^2([T]_{\mathcal{N}})$ may be  correlated.  Indeed, it is possible to give a sequence of DaC measures where magnetisation cannot be reconstructed, but some other odd function can be. \end{rem}

\begin{rem}
In case of product measures, it is immediate from the proof of Theorem \ref{t.cluegen}  (see \cite{GaPe}) that the magnetisation makes the bound sharp. This we cannot conclude for DaC measures in general, as $\Var(\E[f \,|\, \mathcal{N}]$ can be significant for some transitive functions, while it vanishes for magnetization. Still it is worth noting that in~\eqref{DaCspecest} we can write equality instead of inequality in the special case of magnetization:
$$
\p\left[  \Spec^{\mathcal{N}}_{f}  \subseteq \mathcal{U},\; \; \Spec^{\mathcal{N}}_{f} \neq [\emptyset]_\mathcal{N} \right]= \p\left[ \exists u \in \mathcal{U} ,\; \text{such that}\;  [X]_\mathcal{N}=[u]_\mathcal{N} \right]  .
$$
\end{rem}

\begin{rem} \label{rem.enough}
Suppose that, for some sequence of DaC measures $\p_n$, we can show on the one hand that the sequence $\E|\mathcal{C}_{X}^{\mathcal{N}}|$ is bounded by some constant, and on the other hand that $ \clue (f_n \,|\, \mathcal{N}_n)\leq 1-\eps$ for some $\eps>0$. Then there is no sparse reconstruction for $\p_n$.

Indeed, from Theorem~\ref{DaCClueThm} we would then have $\E[\clue(f_n\,|\,\mathcal{U}_n)]\leq 1-\eps/2$ for any sequence of functions, showing that there is no full reconstruction, and then Corollary~\ref{clueto1} says that there is no partial reconstruction, either.
\end{rem}

\begin{rem}\label{r.SW}
There is a connection between the \empha{Swendsen--Wang dynamics} and the question of estimating  $\clue(f_n \,|\, \mathcal{N}_n) = \frac{\Var(\E[f_n \,|\, \mathcal{N}_n])}{\Var(f_n)}$ for the FK-Ising model. 

The Swendsen--Wang algorithm is a non--local Markov chain on the configuration space $ \{-1,1\}^{V}$, where $V$ is the vertex set  of a graph $G(V, E)$, with stationary distribution being the Ising measure at some inverse temperature $\beta$. This chain takes a spin configuration $\sigma^0$, then samples a random cluster configuration $\cN$ which is consistent with  $\sigma^0$, according to the  distribution $\p_{\nu}(\cdot|\sigma^0)$, then generates a new spin configuration $\sigma^1$, chosen according to $\p_\cN$.

We claim that the Swendsen--Wang algorithm having a uniform spectral gap on a sequence of graphs is equivalent to  $\clue(f_n \,|\,\mathcal{N}_n)  < 1-\eps$, where $ \mathcal{N}_n$ is the sigma-algebra generated by the coupled random cluster configuration. Then, in light of Remark~\ref{rem.enough}, knowing that the Swendsen--Wang dynamics mixes fast would solve `half' of the problem of showing no SR --- it would remain to show that $\E|\mathcal{C}_{X}^{\mathcal{N}}|$ is bounded.

The claim above follows directly from a generalization of the  block dynamics $(X^{\mathcal{U}}_t)_{t\in \N}$ (as introduced at the beginning of Section~\ref{s.glauber}) and Lemma~\ref{l.eigenclue}. Consider {\it any} random variable $\cN$ coupled with the spin system $\sigma$. Then, $\cN$ induces a Markovian dynamics  $(X^{\cN}_t)_{t\in \N}$ on $\Omega^{V}$ in the exact ping-pong manner as the random cluster configuration $\cN$ does on the Ising spin configurations. 

We note that this chain is reversible and $\sigma$ is still its stationary distribution. Indeed, if the chain is started from the stationary measure $\p$, then the triple $(X^{\cN}_0, \cN, X^{\cN}_1)$ has the important feature that $(X^{\cN}_0, \cN)$ and $(X^{\cN}_1, \cN)$ are both distributed as $(\sigma, \cN)$, and are conditionally independent given $\cN$. 

We now show that Lemma~\ref{l.eigenclue} remains true in this setting. As in the proof there, we can recognize
\begin{align*}
  \frac{1}{2}\E[(f(X^{\cN}_0) - f(X^{\cN}_1))^2] 
&= \frac{1}{2}\E\left[\E[(f(X^{\cN}_0) - f(X^{\cN}_1))^2 \,|\, \cN]\right]\\
&= \E\left[  \Var(f(X) \,|\,\cN )\right],
\end{align*}
%
 and the remaining part of the proof is the same as in Lemma~\ref{l.eigenclue}.

There are hints that there are even stronger links between RSR and the mixing properties of the Swendsen--Wang algorithm. It is still an open question whether for high temperature Ising models on any reasonable graph sequence the Swendsen--Wang algorithm mixes fast (i.e., if there is a spectral gap). In \cite{BCSV} it is proven that for all high temperature Ising models on $\Z^d$ the relaxation time of the Swendsen--Wang algorithm is $O(1)$, while in \cite{BCV} the same result is shown for high enough temperature on any graph sequence of bounded degree. As large part of the techniques used in Subsections~\ref{ss.SSM} and~\ref{ss.ASSM} are borrowed from \cite{BCSV} and \cite{BCV}, it is natural to ask if there is a clear relationship between the mixing properties  of the Swendsen--Wang algorithm  and the lack or presence of sparse reconstruction (for a particular graph and value of $\beta$). See Question~\ref{q.SW}.
\end{rem}

\section{Sparse reconstruction in the Ising model} \label{s.Ising}


\subsection{Proofs of sparse reconstruction for the (super)critical Ising model, general case}\label{ss.supIsing}

It is well-known that, for $\beta > \beta_c(\Z^d)$,  the sequence $\{ \sigma^n \}$ on the torus $\Z^d_n$ converges to the Ising measure $\frac{1}{2} \mu_{+}+ \frac{1}{2} \mu_{-}$ on $\Z^d$, where $\mu_{+}$ and  $\mu_{-}$  are the weak limits of Ising measures along any exhaustion with $+$ and $-$ respective boundary conditions.  The limiting measure is a convex combination of two translation-invariant measures, thus non-ergodic. Nevertheless, this is not enough to conclude the existence of sparse reconstruction --- see the discussion before Proposition~\ref{SRnonerg}. And this route is clearly unavailable at the  critical temperature. So, we need to be more careful.

\bpf[Proof of Theorem~\ref{supercrit_is}] The first statement of part~(1), full sparse reconstruction for $M_n$, follows directly from Corollary~\ref{nonconstsusc}: we have $\Su(\sigma_n)\to\infty$, plus, by the $\pm$ symmetry of the system, we have $\Var[\sigma_n(v)]=1$.  

Now we prove the second statement of part~(1), sparse reconstruction of $\Maj_n$ on a general sequence of transitive graphs. We are going to show for the total magnetization $S_n$ that
\beq \label{magmoment}
\E[S_n^4]\leq 3\,\E[S_n^2]^2.  
\eeq
Let us first see how sparse reconstruction of $\Maj_n$ follows from this inequality. The Paley-Zygmund inequality implies that, for arbitrary $\eps>0$,
$$
 \p\left[ |S_n| > \eps \sqrt{\E[S^2_n]} \right] = \p\left[ S^2_n > \eps^2 \E[S^{2}_n] \right] \geq (1-\eps^2)^2 \frac{\E[S^{2}_n]^2}{\E[S_n^4]} \geq \frac{(1-\eps^2)^2}{3}. 
$$
By our assumptions, the susceptibility is unbounded, and thus
$$
\sqrt{\E[S^2_n]} = \sqrt{\Var(S_n)} \gg \sqrt{|V_n|}. 
$$
So, by Proposition \ref{SR_maj}, with high probability we have 
$$
\clue(\Maj_n\,|\,\mathcal{B}^{p_n})> \gamma(\mu=0,c=1/3 -\eps) > 0
$$
for every large $n$, with some $p_n\to 0$, giving non-trivial sparse reconstruction.

We are left to show that \eqref{magmoment} holds. Indeed, the so-called Lebowitz inequality (or Gaussian bound; see \cite{L}) says that, for any $x, y, u,v \in V_n$,
$$
\E[\sigma_{x} \sigma_{y} \sigma_{u} \sigma_{v} ]\leq \E[\sigma_{x} \sigma_{y}]\E[ \sigma_{u} \sigma_{v} ]  + \E[\sigma_{x} \sigma_{u}]\E[ \sigma_{y} \sigma_{v} ] + \E[\sigma_{x} \sigma_{v}]\E[ \sigma_{u} \sigma_{y} ].
$$
We can sum these inequalities for all possible quadruples $x, y, u,v \in V_n$, yielding~\eqref{magmoment}.
\medskip

For part~(2), let $\{\cC_i\}_{i\ge 1}$ be the set of FK-clusters, ordered by size: $|\cC_1| \ge |\cC_2| \ge \dots$. By assumption, the probability of the event $\big\{ |\cC_1| \ge \lambda|V_n|$ and $|\cC_2| \le |V_n|/L \big\}$ is at least $1-\eps/2$, for some $\lambda>0$ and any $\eps>0$ and $L>0$, if $n$ is large enough. Conditioning on this event, and all the clusters, the sum of the Ising spins in the clusters $\{\mathcal{C}_i\}_{i\ge 2}$ has conditional variance $\sum_{i \ge 2}{|\mathcal{C}_i|^2}$, which can be bounded as follows: 
$$ 
\sum_{i \ge 2}{|\mathcal{C}_i|^2} \leq \left(\sum_{i \ge 2}{|\mathcal{C}_i|}\right) \max_{i \ge 2}{|\mathcal{C}_i|} \leq (1-\lambda) |V_n|  \frac{|V_n|}{L} \leq  \frac{ |V_n|^2}{L} .
$$
 
If $L$ is chosen large enough, Chebyshev's inequality implies that the absolute value of the sum of the spins outside the giant cluster is smaller than  $\lambda |V_n|/2$ with probability at least $1-\eps/2$. Altogether, with probability at least $1-\eps$, for the total magnetization we have $|S_n| > \lambda |V_n|/2 $, and therefore, by Proposition~\ref{SR_maj}, $\Maj_n$ can be fully reconstructed from a sparse random set, of any density $p_n \gg 1/|V_n|$. 
\medskip

For the critical case of part~(3), we need that the susceptibility blows up. According to \cite{Si}, we have  $\lim_{n \to \infty}{\Su(\sigma^{f}_n)} = \infty $ at $\beta = \beta_c$ where $\sigma^{f}_n$ is the Ising model on $[-n,n]^d$ with free boundary conditions (that is, no magnetization on the boundaries of the square). In order to transfer this result to the torus, we need the fact that there is only one infinite Ising measure at critical temperature. This was shown in \cite{y} for $d =2$ first, in \cite{AF} for $d >3$, and recently in \cite{ADS} for $d =3$.  

This means that the sequence of Ising models on the sequences of tori converges weakly to the same measure as the sequence of  cubes with free boundary conditions, which satisfies $\Su(\sigma^{f}_n) \to \infty$. This implies that the respective covariances converge to the same value as well. Since all covariances are non-negative, we can apply the dominated convergence theorem, and thus we obtain that $\Su(\sigma_n) \to \infty$ as well. 

Regarding the super-critical case of part~(3), it was shown in \cite[Theorem 1.1]{Pi96} that the unique giant condition for the FK random cluster representation holds for boxes and tori of $\Z^d$. (In fact, this was at the time proven for $p$ above the so-called slab threshold for the FK model for  $d>2$, but later it was shown in \cite{Bo}  that this coincides with the critical density $p_c(2)$ that corresponds to the Ising $\beta_c$.)
\medskip

Part (4) could also be proved using the FK representation, studied in \cite{BGJ}. However, there are some statements in that paper that seem wrong to us, hence we prefer to give a proof relying on the direct analytic calculations of \cite{EllNew}. This will be done in Subsection~\ref{ss.CW}.
\epf

\subsection{Proofs for the critical planar Ising model}\label{ss.planar}

Originally, our idea for proving Theorem \ref{critical_is} was to show that the magnetization field scaling limit of \cite{CGN} can be constructed and approximated not only from all the spins, but also from a sparse sample, and hence sparse magnetization must be close to the full magnetization. The following simpler approach, using the scaling limit of correlation functions, was shown to us by Christophe Garban.

\bpf[Proof of Theorem \ref{critical_is}] The main ingredient of the proof is the following special case of \cite[Theorem 1.1]{CHI}: if $x,y \in [0,1]^2$ are two distinct points, and $x_n$ and $y_n$ are lattice points in $Q_n:=\{0,1,\dots,n-1\}^2$ such that $x_n/n \to x$ and $y_n/n \to y$, then 
\beq\label{e.rho2}
n^{1/4} \, \E[\sigma_{x_n}\sigma_{y_n}] \to \rho_2(x,y),
\eeq
where $\rho_2$ is the conformally covariant 2-point function given explicitly in \cite[(1.1), (1.2) and (1.3)]{CHI}.

We consider the sublattice
$$
H_n = H^s_n : = \left\{0, s, 2s \dots, \left \lfloor {\frac{n-1}{s}}\right \rfloor s  \right\}^2,
$$
where $s=s_n \in \Z_+$ with $1 \leq s_n \leq n$ will be fixed later. We will use the notation $M_n^{H_n} := \frac{1}{|H_n|} \sum_{x \in H_n }{\sigma_{x}}$ for the magnetization of $H_n$.

We first want to show that, as $n\to\infty$,
 \beq\label{e.CorrH}
\Corr(M_n, M_n^{H_n}) = \frac{\Cov(M_n, M_n^{H_n})}{\sqrt{ \Var (M_n) \Var (M_n^{H_n})  }} \to 1\,.
\eeq

Using~\eqref{e.rho2}, it was proved in \cite[Chapter 3]{CGN} that
\beq\label{e.convMn}
n^{1/4} \Var(M_n) = \frac{1}{n^{15/4}} \sum_{x,y \in Q_n} \E[\sigma_{x}\sigma_{y}] \to \iint_{[0,1]^2\times [0,1]^2} \rho_2(x,y) \, dx\, dy\,.
\eeq
The only technical difficulty here was to handle pairs $(x,y)$ that are close to each other and hence the convergence in~\eqref{e.rho2} does not hold uniformly for them. To this end, it was proved  in \cite[Proposition 3.5]{CGN} that
\beq\label{e.neardiag}
\frac{1}{n^{15/4}}  \sum_{\substack{x,y \in Q_n \\ \|x-y\| \leq \eps n}}  \E[\sigma_{x}\sigma_{y}] \leq C \, \eps^{7/4}\,;
\eeq
that is, the contribution from nearby pairs is negligible. We need to note here that~\eqref{e.neardiag} is stated in \cite{CGN} with $+$ boundary conditions, but that dominates the expectation with free boundary, hence the inequality remains true in the latter case. This domination holds because in the FK-representation with a wired boundary it is more likely that $x$ and $y$ are connected than with a free boundary (by the FKG inequality for the FK random cluster model), while the conditional expectation of $\sigma_{x}\sigma_{y}$ in the Edwards-Sokal coupling, under the condition that $x$ and $y$ are not in the same FK-cluster, is zero by symmetry. Let us also note that if we looked at the covariance instead of the expectation of the product, the domination between $+$ and free boundaries would be in the reversed direction, by the recent DSS correlation inequality \cite{DSS}.

Now, with the same proof as~(\ref{e.neardiag}), we also have
\beq\label{e.neardiagHn}
\frac{n^{1/4}}{|H_n|^2}  \sum_{\substack{x,y \in H_n \\ \|x-y\| \leq \eps n}}  \E[\sigma_{x}\sigma_{y}] \leq C_1 \frac{n^{1/4}}{|H_n|} + C_2\, \eps^{7/4}\,,
\eeq
where the first term comes from the pairs where the distance between $x$ and $y$ is of constant order (including the diagonal terms $x=y$). This term was also there in~(\ref{e.neardiag}), but there it was automatically negligible compared to the second term, for large $n$. However, now, in order to conclude that
\beq\label{e.convHn}
n^{1/4} \Var(M_n^{H_n}) = \frac{n^{1/4}}{|H_n|^2} \sum_{x,y \in H_n} \E[\sigma_{x}\sigma_{y}] \to \iint_{[0,1]^2\times [0,1]^2} \rho_2(x,y) \, dx\, dy\,,
\eeq
we need to make sure that $\frac{n^{1/4}}{|H_n|}\to 0$. This is indeed the case if the mesh size is $s_n=o(n^{7/8})$. 

Under the condition $\frac{n^{1/4}}{|H_n|}\to 0$, we also have
\beq\label{e.covMnHn}
n^{1/4} \Cov(M_n,M_n^{H_n}) = \frac{n^{1/4}}{|H_n| n^2} \sum_{\substack{x\in Q_n \\ y\in H_n}} \E[\sigma_{x}\sigma_{y}] \to \iint_{[0,1]^2\times [0,1]^2} \rho_2(x,y) \, dx\, dy\,.
\eeq
Combining~(\ref{e.convMn}), (\ref{e.convHn}) and~(\ref{e.covMnHn}) proves that $\Corr(M_n,M_n^{H_n}) \to 1$. 
\medskip

Thus we have full reconstruction for the average total magnetization $M_n$, from the sparse average $M_n^{H_n}$. We will obtain full reconstruction for the majority $\Maj_n = \sign\, M_n$ from $\sign\,M_n^{H_n}$. For this, we use more refined results of Camia, Garban and Newman on the scaling limit of $M_n$. Namely, it was proved in \cite{CGN} that $n^{1/8} M_n$ has a distributional limit, and in \cite{CGN2} that the limit has a smooth density. This implies that, for any $\delta>0$ there is some $\eps>0$ such that, for all $n$ large enough, 
\beq\label{e.sparseMaj}
\begin{aligned}
\p\left[\sign \, M_n \not= \sign \, M_n^{H_n}\right] &\leq \p\big[ |n^{1/8}M_n|<\eps \big] + \p\left[n^{1/8}\big|M_n-M_n^{H_n}\big| > \eps \right] \\
&\leq \delta + \frac{1}{\eps^2} \, \E\left[n^{1/4}\big|M_n-M_n^{H_n}\big|^2 \right]\,,
\end{aligned}
\eeq
where we used Markov's inequality for the second term on the RHS. Now (\ref{e.convMn}), (\ref{e.convHn}) and~(\ref{e.covMnHn}) together imply that $\E\left[n^{1/4}\big|M_n-M_n^{H_n}\big|^2 \right] \to 0$. Hence the RHS of~\eqref{e.sparseMaj} can be made arbitrarily small by first taking $\delta$ small then taking $n$ large. This finishes the proof of Theorem~\ref{critical_is}.
\epf

Let us note that the above proof does not presently work for critical Ising on the tori $\Z_n^2$, despite the results of \cite{IKT} --- the scaling limit of spin-spin correlations is not proved there.

\bpf[Proof of Corollary~\ref{cor.FK}~(2)] In critical FK-Ising on $\Z^2$, the one-arm event $\{ o \leftrightarrow \partial\, [-n,n]^2 \}$ has probability $\asymp n^{-1/8}$ by the classical result of \cite{Onsager}, independently of the boundary conditions, by the RSW-technology developed for critical FK-Ising in \cite{DCHN}. See, specifically, \cite[Lemma 26]{DCHN}. This RSW-technology also implies, as in \cite[Proposition 27]{DCHN}, that, for the cluster $\cC_u$ of a fixed vertex $u$ on the torus $\Z_n^2$, we have 
$$
\E[|\cC_u|] = \sum_{x \in \Z_n^2}\p[u \leftrightarrow x]  \asymp n^2 \, n^{-2/8}
$$ 
and 
$$
\E[|\mathcal{C}_{u}|^2] = \sum_{x,y \in \Z_n^2}\p[x \leftrightarrow u \leftrightarrow y] \asymp n^4 \, n^{-3/8}.
$$
(More precise asymptotics for the free or the plus boundary condition, instead of the torus, are given in \cite[Theorem 1.2]{CHI}.) Thus, Proposition~\ref{DacClueMag} gives that  
$$
\clue(M_n \,|\,U) \leq O(1) \frac{n^{4-3/8}}{n^{2-2/8}} \frac{|U|}{n^2} = O(n^{-1/8}) |U|.
$$

To get the same bound for a general odd transitive function $f$, we use Corollary~\ref{DacClueSym}, together with the fact that, if $\{\cC_k\}_{k\ge 1}$ are the FK-clusters on $\Z_n^2$, then 
$$
\E \left[\max_{k\ge 1}{|\mathcal{C}_k|}\right] \asymp n^{2-1/8}.
$$
For critical percolation, the analogous result was proved in \cite{BCKS}; for the critical FK-Ising model, the proof is the same, using the RSW-technology \cite{DCHN}. \epf

Note that, for magnetization in the critical case, Proposition~\ref{DacClueMag} and Corollary~\ref{DacClueSym} give the same bound. This is in contrast with the high temperature case, where $\max_{k}{|\mathcal{C}_k|}$ would yield an extra factor of $\log n$, which we know to be not sharp by Theorem~\ref{t.norSR}, using Strong Spatial Mixing.

\subsection{Proofs for the Curie-Weiss model} \label{ss.CW}

In this section, we prove Theorem~\ref{supercrit_is}~(4) and Theorem~\ref{t.CWnoSR}, describing exactly what kind of reconstruction is possible for the Curie-Weiss model, the Ising model on the complete graph $K_n$, at different temperatures.

The Curie-Weiss model $\sigma^{\beta,h}_{[n]}$ is defined through the following Hamiltonian:
$$
H(\sigma) = -\frac{1}{n}\sum_{(x,y)} {\sigma_{x}\sigma_{y}} - h  \sum_{x \in [n]}{\sigma_{x}}.
$$
The first sum is over the ${n\choose 2}$ edges of $K_n$. We have the normalisation term $\frac{1}{n}$ to compensate for the fact that the vertex degree is growing linearly with $n$.

For zero external field, $h=0$, there is a phase transition for $\sigma^\beta_{[n]}$ at $\beta=1$. Specifically, we will need the following classical results of Ellis and Newman \cite{EllNew}, also found in \cite{Ellis}: part~(1) is in \cite[Theorem V.9.4]{Ellis}, part~(2) is  in \cite[Theorem V.9.5]{Ellis}, and part~(3) is in \cite[Theorem IV.4.1]{Ellis}. 

\bth[CW magnetization limit theorems \cite{EllNew}]\label{t.EN} 
Consider the Curie-Weiss model $\sigma^\beta_{[n]}$ at $\beta\leq \beta_c=1$, and let $M_n(\sigma^\beta_{[n]})$ be the average magnetization.  
\bit{
\item[{\bf (1)}] For $\beta < 1$, the average magnetization satisfies $\Var[M_n]\asymp 1/n$, and in fact the Central Limit Theorem $\sqrt{n} \, M_n  \Rightarrow \mathsf{N}\left(0,\frac{1}{1-\beta}\right)$ holds, with $\E\left[nM_n^2\right]\to 1/(1-\beta)$.
\item[{\bf (2)}] For $\beta = 1$, we have $\Var[M_n]\asymp 1/\sqrt{n}$, and in fact there is a scaling limit for $n^{1/4} \, M_n$, with $\E\left[\sqrt{n} M_n^2\right]\to c_1 \in (0,\infty)$.
\item[{\bf (3)}] For $\beta > 1$, we have  $\Var[M_n]\asymp 1$, and in fact $M_n \Rightarrow (\delta_{m(\beta)}+\delta_{-m(\beta)})/2$ for some $m(\beta)>0$.
}
\eth

Using these theorems and our general results from Section~\ref{General_SR}, establishing reconstruction for magnetization and majority is easy:

\bpf[Proof of Theorem~\ref{supercrit_is}~(4)] Because of the complete symmetry of the spins and the functions $M_n$ and $\Maj_n$ w.r.t.~permutations, only the size of the set $U_n$ matters, not its distribution.

For $\beta>1$, from Theorem~\ref{t.EN}~(3) we know that $\Corr(\sigma_i,\sigma_j) \asymp 1$ for every $i,j$, and thus  Corollary~\ref{nonconstsusc} tells us that $\clue(M_n \,|\, U_n)\to 1$ for any $|U_n|\to\infty$. Similarly, from Proposition~\ref{SR_maj} we get that  $\clue(\Maj_n \,|\, U_n)\to 1$. 

For $\beta=1$, Theorem~\ref{t.EN}~(2) coupled with Corollary~\ref{nonconstsusc} and Proposition~\ref{SR_maj} give reconstruction for $|U_n| \gg \sqrt{n}$.
\epf 

We now turn to our small clue results. We will use the following information theoretic version of $\clue$, introduced in \cite{GaPe}.

\bde[I-clue] \label{d.Infoclue}
Let $X_V:=\{X_v :v \in V\}$ be a finite family of discrete random variables on some probability space $(\Omega^V,\,\p)$, and for some  $f: \Omega^V \lora \R$ consider the random variable $Z = f(X_V)$. The \empha{information theoretic clue (I-clue)} of $f$ with respect to $U \subseteq V$  is
$$
\clue^{I}(f\,|\,U)= \frac{I(Z : X_U )}{I(Z : X_V )} = \frac{I(Z : X_U )}{\HH(Z)},
$$
where $I(Z : X) := \HH(Z)+\HH(X)-\HH(Z,X)$ is the mutual information. We will also use the conditional entropy $\HH(X \,|\, Z) := \HH(X,Z) - \HH(Z)$, which is also the expectation over $Z$ of the entropy of the conditional distribution of $X$ given $Z$.
\ede

It was proved in \cite[Proposition 3.3]{GaPe} that, for non-degenerate sequences of Boolean functions on any sequence $\big(\{-1,+1\}^{V_n},\p_n\big)$ of probability spaces, sparse reconstruction with respect to the clue according to Definition~\ref{clue} and according to Definition~\ref{d.Infoclue} are equivalent. In the remaining part of this section, we are going to use the I-clue.

\bpf[Proof of Theorem \ref{t.CWnoSR}] 
The proof consists of two main steps. The first one comes from an idea of~\cite{GaPe}.

\bl[Large entropy implies no non-deg RSR]\label{largeentSR} 
Let $\sigma_V$ be a system of random variables such that  
$$
\HH(\sigma_V) \geq  \sum_{i \in V}{\HH(\sigma_{i})} -C\,.
$$
Let $\mathcal{U}$ be a random subset of $V$, independent of $\sigma_V$. As usual, $\delta_{\mathcal{U}} = \max_{i \in V}\p[i \in \mathcal{U}]$. Then, for any  $f: \sigma_V \lora \R$, we have 
$$
\E \left[ \clue^{I}(f\,|\,\cU) \right] \leq \delta_{\mathcal{U}} \left(1 +  \frac{C}{\HH(Z)}\right).
$$
%
\el

\bpf
The proof is very similar to that of \cite[Theorem 4.1]{GaPe}. 
 As in there, we will use the following well-known inequality. For a proof, see, for example, \cite[Theorem 6.28]{LP}, where the statement is formulated for a deterministic family of subsets of $[n]$, but probabilities and expectations over a random subset $\cU$ can be approximated by averages over a deterministic family (allowing for sets appearing with multiplicities), hence the two formulations are equivalent.
 

\bth[Shearer's inequality \cite{Shearer}] \label{shearer}
Let $\sigma_1, \sigma_2, \dots, \sigma_n$ be random variables defined on the same probability space. Let $\mathcal{U}$ be a random subset of $[n]$ independent of $\sigma$. Then
$$
\HH(\sigma_{[n]}) \, \min_{i \in V}\p[i \in \mathcal{U}]  \leq \E[ \HH(\sigma_{\mathcal{U}})].
$$
\eth

We will show that, if the spin system $\sigma_{[n]}$ and the random set $\mathcal{U}$ satisfy the conditions of the lemma, then 
\beq \label{revshearer}
\E[I(Z : \sigma_{\mathcal{U}} )]\leq \delta_{\mathcal{U}} \, \big( I( Z: \sigma_V) +C \big) =  \delta_{\mathcal{U}} \, \big( \HH(Z) +C \big).
\eeq
Dividing by $\HH(Z)$, this gives  Lemma~\ref{largeentSR}.

Without loss of generality we can assume that $\p[i \in \mathcal{U}] = \delta_{\mathcal{U}}$ for every $i \in [n]$. Indeed, if
this is not the case, we can always blow up $\mathcal{U}$ to fullfill this. It is clear that, by adding more spins to  $\mathcal{U}$, the right hand side of the inequality does not change, while the left hand side can only increase.

Observe now that 
\beq \label{almostindep}
 \E\left[ \HH\left(\sigma_{\mathcal{U}}\right)\right] \leq  \E\left[ \sum_{i \in \mathcal{U}} \HH(\sigma_{i}) \right] =  \delta_{\mathcal{U}} {\sum_{i \in V}{\HH(\sigma_{i}) }} \leq \delta_{\mathcal{U}} (\HH(\sigma_V)+ C),
\eeq 
where, for the last inequality, we used  the condition of the lemma. In turn, Shearer's inequality says that
\beq \label{searer}
- \E[\HH(\sigma_{\mathcal{U}} \,|\, Z)] \leq -\delta_{\mathcal{U}} \,  \HH(\sigma_V \,|\, Z).
\eeq
Summing \eqref{almostindep} and \eqref{searer}, we get \eqref{revshearer}, finishing the proof of Lemma~\ref{largeentSR}.
%
%
%
%
\epf

The second ingredient is that Lemma~\ref{largeentSR} is satisfied by the subcritical Curie-Weiss model with $C=O(1)$ and by the critical one with $C=O(\sqrt{n})$. This is an observation of some independent interest. Note that the Ising model on any sequence of bounded degree finite graphs at any $\beta>0$ does not satisfy this, for the following reason. In any exploration process of the spins, $V_n=\{v_1,v_2,\dots,v_n\}$, if the already explored neighbours of the next vertex in the exploration do not happen to have an equal number of $+$ and $-$ spins (which certainly happens with a positive frequency during the exploration), then the conditional distribution of the next spin is a coin flip that is biased by a uniformly positive amount, hence its $\log_2$-entropy is strictly bounded away from 1. Thus, using the chain rule for entropy,
\beq \label{chainrule}
\HH(\sigma_{V_n}) = \sum_{k=0}^{n-1}{\HH\big(\sigma_{v_{k+1}}\,|\,\sigma_{\{v_1,\dots,v_k\}} \big)},
\eeq
the total entropy in the system can be linear with at most a constant factor strictly less than 1. 

These large entropy results for the Curie-Weiss model are related to the propagation of chaos phenomenon \cite{chaos1,chaos2}, which says, e.g., in the subcritical case, that if we look at $o(n)$ of the spins, the joint distribution is very close to being IID.

\bpr[The entropy of (sub)critical Curie-Weiss]\label{CWlargeent} Consider the CW model $\sigma^\beta_{[n]}$ on $K_n$ at external field $h=0$ and $\beta\leq \beta_c=1$.
\bit{
\item[{\bf (1)}] For $\beta < \beta_c$, there exists a positive constant $C_\beta$ such that, for all large enough $n$,
$$
\HH(\sigma^\beta_{[n]}) \geq n - C_\beta.
$$
\item[{\bf (2)}] For $\beta=\beta_c$, we have 
$$
\HH(\sigma^{\beta_c}_{[n]}) \geq n - C \sqrt{n},
$$
with an absolute constant $C<\infty$.
}
\epr

Originally, our proof was based on the chain rule~\eqref{chainrule}, by noting that the already observed spins act on the next spin as an external field, hence the distribution of that next spin can be estimated quite effectively. Then Amir Dembo suggested the following much more elegant approach, which we work out here with his kind permission.

\bpf
We will use the notion of {\it relative entropy} between two random variables $X$ and $Y$ (or more precisely, between their distributions), also called the {\it Kullback-Leibler divergence}:
$$
D(X\,\|\,Y)=\sum_{x\in S} \p[X=x]\,\log\frac{\p[X=x]}{\p[Y=x]}\,,
$$
with $\log$ still denoting $\log_2$.  First notice that, since $\sigma_{0}[n]$ is just the uniform distribution on $|\Omega|=2^n$ configurations,
\beq\label{e.KLent}
D\left(\sigma^{\beta}_{[n]} \,\big\|\,\sigma^{0}_{[n]} \right)=\sum_{\omega \in \Omega} \p\big[\sigma^{\beta}_{[n]} =\omega\big] \left\{ \log\p\big[\sigma^{\beta}_{[n]} =\omega\big]+\log(2^n) \right\} = n-\HH(\sigma^{\beta}_{[n]}),
\eeq
hence our task is exactly the computation of this relative entropy. For this, the main idea is that the distribution of $\sigma^{\beta}_{[n]}$, conditioned on any possible value $M_n(\sigma)=x \in [-1,1]$ for the average magnetization, is just the uniform distribution on the ${n\choose n (1+x)/2 }$ configurations with $n (1+x)/2$ plus and $n (1-x)/2$ minus spins, regardless of $\beta$. Therefore,  the relative entropy decomposition \cite[Theorem D.13]{DZ} says that all the relative entropy comes from the difference in the  distributions of the magnetization:
\beq\label{e.KLCW}
D\left(\sigma^{\beta}_{[n]} \,\big\|\,\sigma^{0}_{[n]} \right) = D\left(M_n(\sigma^{\beta}_{[n]}) \,\big\|\, M_n(\sigma^{0}_{[n]}) \right).
\eeq
Now, since $H(\sigma)=nM_n(\sigma)^2/2$, we have
$$
\p\left[M_n(\sigma^{\beta}_{[n]})=x\right]= \frac{1}{Z_{n,\beta}} {n \choose n(1+x)/2} \exp( \beta n x^2/2),
$$
thus 
$$
\log\frac{\p\left[M_n(\sigma^{\beta}_{[n]})=x\right]}{\p\left[M_n(\sigma^{0}_{[n]})=x\right]}= \log\frac{Z_{n,0}}{Z_{n,\beta}} + \log(e) \beta n x^2/2,
$$
and
\beq\label{e.KLmagn}
\begin{aligned}
D\left(M_n(\sigma^{\beta}_{[n]}) \,\big\|\, M_n(\sigma^{0}_{[n]}) \right) &=  \frac{\log(e) \beta}{2} n \, \E\left[ M_n(\sigma^{\beta}_{[n]})^2  \right] - \log\frac{Z_{n,\beta}}{Z_{n,0}} \\
&\leq  \frac{\log(e) \beta}{2} n \, \E\left[ M_n(\sigma^{\beta}_{[n]})^2  \right],
\end{aligned}
\eeq
where the inequality is from the trivial observation $Z_{n,\beta} \ge Z_{n,0}$.

%
%
%

Now substituting the limit results of Theorem~\ref{t.EN}~(1) and (2) into~(\ref{e.KLmagn}), we get 
$$
D\left(M_n(\sigma^{\beta}_{[n]}) \,\big\|\, M_n(\sigma^{0}_{[n]}) \right) \leq 
\begin{cases}
\frac{\log(e) \beta}{2(1-\beta)} &\text{ for }\beta<1,\\
\frac{\log(e) c_1}{2} \sqrt{n} &\text{ for }\beta=1.
\end{cases}
$$
Using~(\ref{e.KLCW}) and~(\ref{e.KLent}), this finishes the proof of Proposition~\ref{CWlargeent}. \epf

Lemma~\ref{largeentSR} (large entropy implies small clue) and Proposition~\ref{CWlargeent} (the Curie-Weiss model has large entropy) together imply Theorem~\ref{t.CWnoSR} immediately.
\epf

\section{Open problems}\label{s.open}

In Example~\ref{SRvsRSR}, the limiting spin system (say, in the Benjamini-Schramm sense) on $\Z\times\Z_2$ is ergodic, but not totally ergodic (there are non-trivial functions that are invariant under the infinite subgroup $\Z$), not weakly mixing, not a factor of IID measure, not a Markov random field. This leads to the following question:

\begin{quest}[SR vs RSR]\label{q.SRvsRSR}
Suppose that for a transitive, Benjamini-Schramm convergent sequence of spin systems $\p_n$, there is RSR. Under what condition does this imply that there is also SR for $\p_n$?  
\end{quest}

Another possible distinction between SR and RSR comes from the different quantitative bounds obtained in Corollary~\ref{c.quantnoSR}:

\begin{quest}[Quantitative SR vs RSR]\label{q.RSRsqrt}
Is the bound $\E[\clue(f_n\,|\,\mathcal{U}_n)]\leq C \sqrt{\delta(\mathcal{U}_n)}$  of Corollary~\ref{c.quantnoSR} for RSR ever sharp?
\end{quest}

Our next question concerns whether sparse reconstruction can be a property of the limiting measure, despite simple  counterexamples like Examples~\ref{prodlim} and~\ref{nonerglim}: 

\begin{quest}[Robust Sparse Reconstruction]\label{q.robust}
Is there an ergodic $\mathsf{Aut}(G)$-invariant spin system $(G, \p)$ such that whenever $(G_n, \p_n)$ converges  to $(G, \p)$ in the Benjamini-Schramm sense, there is (R)SR for $(G_n, \p_n)$?
\end{quest}

An interesting example is the critical Ising model on $\Z^2$, which is ergodic, but the magnetization $M_n$ on tori admits full sparse reconstruction, by Theorem~\ref{supercrit_is}~(3). At the same time, one can take the sequence of critical Ising measures conditioned on $M_n=0$. It can be shown that this spin system also converges to the critical Ising on $\Z^2$, but the magnetization clearly cannot be reconstructed. In this model, however, we suspect that other transitive functions can still be reconstructed. More generally, the critical planar Ising model could provide a positive answer to Question~\ref{q.robust}.

A similar question is the following, using the notion of asymptotic entropy, Definition~\ref{d.asympent}.

\begin{quest}[Zero entropy]\label{q.zero}
Let $(G_n, \p_n)$ converge  to a non-trivial $\mathsf{Aut}(G)$-ergodic spin system $(G, \p)$. If $\mathcal{H}(\{ \p_n \})=0$, then is it true that there is (R)SR for $(G_n, \p_n)$?
\end{quest}

Another positive answer to Question~\ref{q.robust} could come from the following important example, which we are planning to address in future work:

\begin{quest}[2d Ising plus phase]\label{q.plus}
Consider the Ising model on $\Z_n^2$, at $\beta > \beta_c$, conditioned to have positive total magnetization. Is there (non-degenerate, random) sparse reconstruction?
\end{quest}

We suspect a positive answer, despite the model having bounded susceptibility. The paper \cite{BM} shows the vanishing of the spectral gap using an observable constructed from the relatively large minus islands that can exist due to the subexponential large deviations in the model. Similarly, we think that the volume of the largest minus island, mean of order $\log^2 n$, standard deviation of order $\log n$, could be a function admitting sparse reconstruction. (For comparison: in a subcritical percolation model, where the volume of the largest cluster has mean of order $\log n$, standard deviation of order 1, detecting whether it is larger than its median looks definitely impossible using a sparse sample.) 
\medskip

Notwithstanding the above possible examples, it seems that, for FFIID measures, finite or infinite susceptibility might indeed be a decisive factor, which brings us to possibly the main open question of this paper. Let us remark here that the low temperature Ising plus phase is known to be not a FFIID measure exactly from the above-mentioned property of having subexponential large deviations for magnetization \cite[Theorem 2.1]{BS}. 

\begin{quest}[Finite expected volume FFIID]\label{q.adam}
Suppose that  a sequence of FFIID spin systems $(G_n, \p_n)$ converges to the FFIID $(G, \p)$ regularly (see Definition~\ref{d.Adam}), where $(\p, G)$ has finite expected coding volume. Is it true that there is no sparse reconstruction on  $(G_n, \p_n)$?
\end{quest}

Despite all our results, there are still some open problems even in the most classical case, the critical Ising model on the $\Z^d$ lattices.

\begin{quest}[Full sparse reconstruction for critical Ising magnetization]\label{q.IsingMgen}
Prove the analogues of Theorem~\ref{critical_is} for the boxes $[0,n]^d$ with $d\ge 3$, and tori $\Z_n^d$ with $d\ge 2$. In particular, prove that not only does the susceptibility blow up at criticality, but also, the average magnetization is much larger than $n^{-d/2}$ with high probability.
\end{quest}

The fact that Proposition~\ref{DacClueMag} or  Corollary~\ref{DacClueSym} do not give the sharp bounds of Theorem~\ref{t.CWnoSR} for the Curie-Weiss model suggests that, also in the planar case, the quantitative bound of Theorem~\ref{critical_is} might be sharp, as opposed to the bound of Corollary~\ref{cor.FK}~(2):

\begin{quest}[Quantitative bounds for 2d critical Ising]\label{q.IsingMplane}
For critical Ising on two-dimensional boxes $[0,n]^2$ or tori $\Z_n^2$, is reconstruction possible from $U_n$ when $n^{1/8} \leq |U_n| \leq n^{1/4}$,
\bit{
\item[{\bf (a)}] for average magnetisation?
\item[{\bf (b)}] for any transitive function?
     }
\end{quest}

The next question is about the Ising model beyond $\Z^d$. On the $d$-regular tree $\T_d$, $d\ge 3$, there are three main regimes for the Ising model. Letting $\eps=\eps(\beta)=(1+e^{2\beta})^{-1}$, the FK representation of the model on any finite tree, or on an infinite tree with the free boundary condition, simplifies to first performing Bernoulli$(1-2\eps)$ bond percolation, then colouring each cluster $\pm$ according to an independent fair coin. The high temperature regime on $\T_d$ is $1-2\eps \leq \frac{1}{d-1}$, where every free FK cluster is finite, and one can use this to show that there is only one extremal Ising Gibbs measure. The intermediate regime is $\frac{1}{d-1} < 1-2\eps \leq \frac{1}{\sqrt{d-1}}$, with three extremal translation-invariant Gibbs measures: the $+$, the $-$, and the free measure. In the low temperature regime, $1-2\eps > \frac{1}{\sqrt{d-1}}$, only the $+$ and the $-$ measures are the extremal translation-invariant Gibbs measures; even though the free measure is still ergodic, it is a convex combination of non-translation-invariant Gibbs measures. See \cite{MP} for an overview.

\begin{quest}[Free Ising approximations]\label{q.free} Let $\p$ be the free Ising measure on $\T_d$ for some $\eps$ in the intermediate regime $\frac{1}{d-1} < 1-2\eps \leq \frac{1}{\sqrt{d-1}}$. 
\bit{
\item[{\bf (a)}] Is it true that along every exhaustion $V_n\nearrow \T_d$ the projected measures $\p_n$ have no RSR? 
\item[{\bf (b)}] Is it true that for every Benjamini-Schramm approximation of $\p$ there is RSR?\\
 For approximations that have non-negative two-point correlations, one can show (similarly to Remark~\ref{r.positive}) that $\E[\Var(M_n)] \gg 1/n$, hence the average magnetization $M_n$ has RSR by Corollary~\ref{nonconstsusc}.
}
\end{quest}

For part (a), our main reason for guessing no RSR is that if $V_n$ is a large ball with $n$ vertices, then $\Var[M_n] \asymp 1/n$. This is also the reason why census reconstruction from the boundary to the root does not work, and it turns out that for this measure this implies (in a highly non-trivial manner) that no reconstruction works --- i.e., the measure is tail-trivial, hence extremal. 

For part (b), note that the Ising model on the $d$-regular random graph $G_{n,d}$ does not converge to the free measure $\p$, because for $1-2\eps > \frac{1}{d-1}$ there is a unique giant cluster in the FK random cluster model on $G_{n,d}$. (Unfortunately, we have not managed to locate a reference, but we are sure the claim is correct. But even if it is false, that would just be good news for finding more approximations.) Nevertheless, there do exist Benjamini-Schramm approximations: from the general results of \cite{Elek,BLS}, or, for a large part of the intermediate regime, from the FIID construction of \cite{NSZ}. It is not clear though how to get a Benjamini-Schramm approximation with positive two-point correlations. 
\medskip

After these special Ising cases, two general bold questions are the following:

\begin{quest}[Swendsen-Wang spectral gap vs RSR]\label{q.SW}
In the spirit of Remark~\ref{r.SW}, is a uniform spectral gap for the Swendsen-Wang dynamics of the high temperature Ising model on some sequence of graphs equivalent to the impossibility of RSR? 
\end{quest}

\begin{quest}[Glauber spectral gap vs RSR]\label{q.gap} Consider any spin system on any sequence of finite transitive graphs $G_n$. Is it true that the continuous time Glauber dynamics, updating the status of each vertex according to an independent Poisson clock of rate 1, has a uniformly positive spectral gap if{f} there is no random sparse reconstruction?
\end{quest}

Strictly speaking, the spin systems in the latter question should have the finite energy property, so that it is clear what we mean by Glauber dynamics, resampling a single spin. However, one can also consider natural extensions, such as dynamics resampling a bounded portion of the configuration; see Question~\ref{q.pm}~(c) for a natural example.

Note that \cite{MO1} deduced a uniform spectral gap (and more) for the Glauber dynamics from a weak form of spatial mixing, which was already known to hold at that time for all subcritical Ising models. In this sense, Theorem~\ref{t.norSR} is a special case of Question~\ref{q.gap}. Furthermore, by \cite{BM}, Question~\ref{q.plus} is also a special case of this question, just in the other direction. 

It was shown in  \cite{BKMP} that a uniformly positive spectral gap for the Glauber dynamics of the Ising model on finite balls on regular trees, with an arbitrary boundary condition at the leaves, holds precisely when the free Ising measure is tail trivial. That is,  Question~\ref{q.free}~(a) is also special case of Question~\ref{q.gap}. Furthermore, for the \empha{hard-core model} on finite balls of regular trees, with an arbitrary boundary condition at the leaves, \cite{RSVVY} showed that the spectral gap is almost constant in the tail-trivial regime, while, in the opposite direction, any reconstruction algorithm gives rise to a bottleneck for the dynamics, giving a significantly smaller spectral gap. In all these examples it is of course completely unclear why non-reconstruction for local functions (tail-triviality of the limiting measure along exhaustions) would imply no (random) sparse reconstruction for all functions; however, the fact that we have an (almost) uniform spectral gap in the tail-trivial regime seems to suggest that somehow a global relaxation is also implied, giving some support to Question~\ref{q.gap}.

Recall that we proved two special cases of the following conjecture in Subsection~\ref{FFIID.tt}:

\bcj[Bal\'azs Szegedy]\label{c.SzyB} Any factor of IID measure on any transitive infinite graph has trivial sparse tail: there exist no sparse (Bernoulli) reconstructable functions (see Definitions~\ref{d.stber} and~\ref{d.stgen}).
\ecj

Here is an example, quite different from Ising measures, to test some of the above general conjectures.

\begin{quest}[Perfect matchings on large girth graphs]\label{q.pm} Consider the unique automorphism-invariant perfect matching measure $\p$ on the tree $\T_3$, which is a FIID measure by \cite{LyNaz}.
\bit{
\item[{\bf (a)}] Show that it has a trivial sparse tail, as predicted by Conjecture~\ref{c.SzyB}.
\item[{\bf (b)}] Give a sofic approximation $\{G_n,\p_n\}$ to $(\T_3,\p)$ that does admit random sparse reconstruction.
\item[{\bf (c)}] For any local resampling Markov chain for the systems $(G_n,\p_n)$ of the previous item, show that there is no uniform spectral gap.
\item[{\bf (d)}] Show that $\p$ is not a FFIID measure with finite expected coding volume.
}
It is of course possible that (b) is correct, while (c) and/or (d) is not, in which case we would get a negative answer to Question~\ref{q.gap} and/or Question~\ref{q.adam}, respectively.
\end{quest}

As a final example, consider uniform random proper 3-colorings of the boxes $[n]^2$, a critical system in many senses, a version of the 6-vertex model. See, e.g., \cite{PelSpi} and \cite[Section 9]{GlaP}. To be precise, take a uniform $h: [n]^2\lora \Z$ height function with $|h(x)-h(y)|=1$ for every edge $(x,y)$, for $n$ even, with boundary condition being 0 on every second vertex of the boundary, then take the colouring $\sigma(x):=h(x)$ (mod 3). This height function has variance $h(x)\asymp \log n$, suggesting critical 2-dimensional behaviour. Moreover, not only $h(x)$, but also $|h(x)|$ satisfies the FKG-inequality, hence two-point correlations are positive. All these properties suggest that parts (a), (b) of following question might actually be quite straightforward; however, writing up the proofs would still be worthwhile.

\begin{quest}[Proper 3-colorings of 2d boxes]\label{q.3col}
Consider the uniform height functions $h$ on boxes $[n]^2$ with boundary conditions as above, and the corresponding (mod 3) coloring. 
\bit{
\item[{\bf (a)}] Show that $h$ has sparse reconstruction.
\item[{\bf (b)}] Show that the 3-coloring itself has sparse reconstruction.
\item[{\bf (c)}] What happens for the uniform proper 3-coloring of the torus $\Z_n^2$?
}
\end{quest}


\ \\
{\bf P\'al Galicza}\\
\texttt{galicza[at]gmail.com}\\
\ \\
{\bf G\'abor Pete}\\
HUN-REN Alfr\'ed R\'enyi Institute of Mathematics, Re\'altanoda u. 13-15, Budapest 1053 Hungary, and\\
Department of Stochastics, Institute of Mathematics, Budapest University of Technology and Economics, M\H{u}egyetem rkp.~3., Budapest 1111 Hungary\\
\texttt{gabor.pete[at]renyi.hu}, \url{http://www.math.bme.hu/~gabor}\\
\end{document}